%% file: LERW160718.tex
\documentclass[11pt]{article}
\usepackage{amssymb,latexsym}
\usepackage[dvips]{graphicx} 
\usepackage{url}
\usepackage{picture}
\usepackage{color}

\newtheorem{prp}{Proposition}
\newtheorem{lem}[prp]{Lemma}\newtheorem{thm}[prp]{Theorem}

\newtheorem{cor}[prp]{Corollary}
\newenvironment{prf}{\begin{trivlist}\item[\emph{Proof.}]}{\end{trivlist}
  \medskip\par}
 
\newenvironment{rem}{\begin{trivlist}\item[\emph{Remarks.}]}{\end{trivlist}  \medskip\par}
\def\prpb{\begin{prp}}\def\prpe{\end{prp}}
\def\lemb{\begin{lem}}\def\leme{\end{lem}}
\def\thmb{\begin{thm}}\def\thme{\end{thm}}
\def\corb{\begin{cor}}\def\core{\end{cor}}
\def\prfb{\begin{prf}}\def\prfe{\end{prf}}
\def\remb{\begin{rem}}\def\reme{\end{rem}}
\def\prpa#1{\label{p:#1}}\def\prpu#1{Proposition~\ref{p:#1}}
\def\lema#1{\label{l:#1}}\def\lemu#1{Lemma~\ref{l:#1}}
\def\thma#1{\label{t:#1}}\def\thmu#1{Theorem~\ref{t:#1}}

\def\seca#1{\label{s:#1}}\def\secu#1{Section~\ref{s:#1}}

\def\itmb{\begin{enumerate}}\def\itme{\end{enumerate}}
\def\itdb{\begin{itemize}}\def\itde{\end{itemize}}
\def\ittb{\begin{description}}\def\itte{\end{description}}
\def\eqnb{\begin{equation}}\def\eqne{\end{equation}}
\def\arrb#1{\begin{array}{#1}}\def\arre{\end{array}}
\def\tabb#1{\par\noindent\begin{tabular}{#1}}
\def\tabe{\end{tabular}\par\noindent}
\def\eqna#1{\label{e:#1}}\def\eqnu#1{(\ref{e:#1})}

\def\QED{\relax\ifmmode\let\@tempa\relax\ifcase\@eqcnt\def\@tempa{& & &}\or
  \def\@tempa{& &}\else\def\@tempa{&}\fi\@tempa $\Box$ \else\hfill $\Box$ \fi}
\def\DDD{\relax\ifmmode\let\@tempa\relax\ifcase\@eqcnt\def\@tempa{& & &}\or
 \def\@tempa{& &}\else\def\@tempa{&}\fi\@tempa $\Diamond$
 \else\hfill $\Diamond$ \fi}

\def\Rom#1{\uppercase\expandafter{\romannumeral#1}}
\def\dsp{\displaystyle}

\def\liminf{\displaystyle \mathop{\underline{\lim}}\limits}
\def\limsup{\displaystyle \mathop{\overline{\lim}}\limits}

\def\Ccomb#1#2{\setbox0=\hbox{$\displaystyle\mathrm{C}$}\setbox1=\hbox{%
$\scriptstyle #1$}\kern \wd1{\mathrm{C}}_{\kern -1.05\wd0\kern -0.99\wd1{#1}
 \kern 1.15\wd0{#2}}}

\def\clvec#1#2#3{\def\clvecone{#3}\left(\arrb{c} \dsp #1\\ \dsp #2
 \ifx\clvecone\empty\else\\ \dsp #3\fi\arre\right)}

 \def\leq{\leqq} \def\geq{\geqq}

\def\pintegers{{\mathbb Z}_+}
\def\nintegers{{\mathbb N}}

\def\prb#1{\def\prbone{#1}
  \ifx\prbone\empty{\mathrm{P}}\else{\mathrm{P[\;}}#1{\mathrm{\;]}}\fi}
\def\prbseq#1#2{\def\prbseqone{#2}
  \ifx\prbseqone\empty{\mathrm{P}}_{#1}\ignorespaces
  \else{\mathrm{P}}_{#1}{\mathrm{[\;}}#2{\mathrm{\;]}}\fi}

\def\EEseq#1#2{\def\EEseqone{#2}
  \ifx\EEseqone\empty{\mathrm{E}}_{#1}\else
 {\mathrm{E}}_{#1}{\dsp\mathrm{[\;}}#2{\mathrm{\;]}}\fi}

\def\VVseq#1#2{\def\VVseqone{#2}
  \ifx\VVseqone\empty{\matrm{V}}_{#1}\else
 {\mathrm{V}}_{#1}{\dsp\mathrm{[\;}}#2{\mathrm{\;]}}\fi}

\def\sg{Sierpi\'{n}ski gasket}
\def\parr{\par\noindent}
\def\eqsb{\begin{eqnarray*}}\def\eqse{\end{eqnarray*}}

\makeatletter
\@addtoreset{equation}{section}
\makeatother

\title{
Displacement exponent for loop-erased random walk on the Sierpi\'nski gasket
}
\author{
Kumiko Hattori
}

\date{\today}
\begin{document}
\maketitle

\footnotetext{Department of Mathematics and Information Sciences,
 Tokyo Metropolitan University, Hachioji, Tokyo 192-0397, Japan.}


\begin{center}
ABSTRACT
\end{center}
We prove that loop-erased random walks on finite pre-Sierpi\'{n}ski gaskets can be extended to the infinite pre-Sierpi\'{n}ski gasket by virtue of the 
  `erasing-larger-loops-first'  method, and 
obtain the asymptotic behavior of the walk as the number of steps increases, in particular, 
the displacement exponent and a  
law of the iterated logarithm.

\parr

\vspace*{1in}\par

\bigskip\par
\noindent\textit{MSC2010 Subject Classifications:}
60F99, 60G17, 28A80, 37F25, 37F35 
\bigskip\par
\noindent\textit{Corresponding author:} 
Kumiko Hattori, 
khattori@tmu.ac.jp
\parr
Department of Mathematics and Information Sciences,
 Tokyo Metropolitan University, Hachioji, Tokyo 192-0397, Japan.
\parr tel: +81 42 677 2475


\section{Introduction}
\seca{Intro}

Loop-erased random walk (LERW) is a 
process obtained by erasing loops from a simple random 
walk in chronological order (as soon as each loop is made).  
LERW was originally considered on ${\mathbb Z}^d$ and    
the existence of the scaling limit has been proved for all $d$.
The asymptotic behavior of the walk has been studied in terms of 
the growth exponent (expected to be the reciprocal of the displacement 
exponent).  For the growth exponents for LERW on  ${\mathbb Z}^d$, see, for example, \cite{Lawler}, \cite{Lawler2}, \cite{Lawler3}, \cite{Kenyon} and \cite{Shiraishi}.
  

In this paper, we consider LERW on the Sierpi\'nski gasket and prove the following Theorems 1--3.
\par\noindent
\thmb
\thma{infinite}
Loop-erased random walks on the finite  Sierpi\'nski gaskets  
can be extended to a
loop-erased random walk on the infinite \sg .
\thme
Let $\lambda =(20+\sqrt{205})/15$ and  $\nu= \log 2/\log \lambda $. 

\par\noindent
\thmb
\thma{exponent}
For any $s >0$, there exist positive constants $C_1(s)$ and $C_2(s)$ such that   
\[C_1(s) n^{s\nu } \leq E[|X(n)|^s] \leq C_2(s)n^{s\nu } ,\]
where $X(n)$ denotes the location of the LERW starting at the origin after $n$ steps and  $|\ \cdot \ | $ the Euclidean distance. 

\thme
$\nu$ is called the displacement exponent. 
\thmb
\thma{loglog}
There are positive constants $C_{3}$ and $C_{4}$  such that 
\[C_{3} \leq \limsup _{n\to \infty } \frac{|X(n)|}{\psi (n)}\leq C_{4},  \mbox{ a.s.},\]
where $\psi (n)=n^{\nu}(\log \log n)^{1-\nu}$.
\thme

Our main tool for the proof is  the `erasing-larger-loops-first' (ELLF) method, which was introduced to study the scaling limit (the limit as the edge length tends to $0$).
The scaling limit for LERW on the Sierpinski gasket was obtained by two groups independently, using different methods. 
For the `standard' LERW on general graphs, 
the uniform spanning tree proves to be a powerful tool (\cite{STW}). 
By `standard', we mean the loops are erased chronologically from a simple random walk  as first introduced by G. Lawler  (\cite{Lawler}).  On the other hand, 
\cite{HM} constructed a LERW on the Sierpi\'{n}ski gasket by ELLF, that is, 
by erasing loops in 
descending order of size of loops and proved that the resulting LERW has the same 
distribution as that of the `standard' LERW.  Futhermore, 
in \cite{HOO}, it is proved that ELLF does work not only for simple random walks, but also for other kinds of random walks on some fractals, in particular, for self-repelling walks on the Sierpi\'nski gasket introduced in \cite{HHH}. An important reason for this flexibility is that the ELLF method is based on self-similarity of the Sierpi\'nski gasket.

Another advantage of the ELLF method is facilitate the extension of   
LERW to the infinite Sierpi\'nski gasket  by providing us with a natural definition of 
two series of probability measures on sets of loopless paths.
The extension is not trivial, for the simple random walk on the infinite  Sierpi\'nski gasket is recurrent.  
The exact value of the displacement exponent has been known by a scaling argument (\cite{DD}).  As for the 
proof of the existence,  
the authors erroneously wrote in \cite{HM} that 
\thmu{exponent} has been proved in \cite{STW}, however,  \cite{STW}
deals with the scaling limit, not LERW on the infinite \sg , and proves the short-time behavior of the limit process $\overline{X}(t)$: 
\thmb (Theorem 7.10 in \cite{STW})
\thma{shorttime}
For any $p>0$, there exist constants $C_5(p)$,  $C_{6}(p)>0$ such that  for all $t \in [0,1]$,    
\[C_5(p) t^{p \nu }\leq E[|\overline{X}(t)|^p] \leq C_{6}(p)  t^{p \nu },\]
where $|\overline{X}(t)|$ denotes the Euclidean distance from the starting point at time $t$ and 
$\nu= \log 2/\log \lambda $, $\lambda =(20+\sqrt{205})/15$. 
\thme

It is expected that the same exponent also rules the long-time behavior of the 
walk, but the method of proof is different, for one has to look into how the scaled 
number of steps converges, not only the limit distribution.  
Thus, the author corrects her error and proves \thmu{exponent} 
in this paper.

The first mathematical result on the displacement exponent for a non-Markov random  walk on the  Sierpi\'nski gasket was obtained in \cite{HK},  dealing with the `standard' self-avoiding walk, which is defined by the uniform measure on self-avoiding paths of a given length.  They showed the existence of the exponent in the form of 
\eqnb\eqna{nuSAW} \lim_{n \to \infty} \frac{\log E_n[|X'(n)|^s]}{\log n}=s\nu_{SAW}, \ \ s>0\eqne
where $|X'(n)|$ denotes the end-to-end distance of an $n$-step self-avoiding path, 
 and  $\nu_{SAW}=\log 2/\log (\frac{7-\sqrt{5}}{2})$.  Since the exponent $\nu_{SAW}$ is different from $\nu$ in \thmu{exponent}, the LERW is in a different universality class from the 
self-avoiding walk.  Note that self-avoiding walk cannot be extended to infinite length, for 
the consistency condition is not satisfied because of culs-de-sac, thus the expectation is taken over the uniform measure on the $n$-step self-avoiding paths.  
Note also that we have a sharper result in \eqnu{nuSAW}, which comes from the refinement in the analysis.  

The structure of the paper is as follows. 
In \secu{SRW}, we define our notation and 
in \secu{LE}, we describe the ELLF method of loop-erasing. 
\secu{Steps} deals with the asymptotics of the exit times from  
a series of triangles, which is used in \secu{Proof}.
In \secu{Extention} we extend the walk to the infinite Sierpi\'{n}ski gasket and finally, 
in \secu{Proof} we prove Theorems 2 and 3. 

\section{Random walk on the pre-\sg s}
\seca{SRW}

\subsection{The pre-{\sg}s}

Let us recall the definition of the pre-{\sg}: denote  
$O=(0,0),\,
a_0=({1 \over 2},{\sqrt{3} \over 2}),\,
b_0=(1,0)\, $, 
$a_N=2^Na_0$ and $b_N=2^Nb_0$ for  $N \in 
{\mathbb N}$. 
Let $F'_{0}$ be the graph that consists of the three vertices and three edges of $\triangle Oa_0b_0 $
and define a recursive sequence of graphs 
$\{F'_{N}\}_{N=0}^{\infty }$ by
\[F'_{N+1}=F'_{N} \cup (F'_{N}+a_N) \cup (F'_{N}+b_N)
, \ \  N \in {\mathbb Z}_+ =\{0,1,2, \ldots \}\, ,
\]
where $A+a=\{x+a\ :\  x\in A\}$ and $kA=\{kx\ :\  x\in A\}$.  
$F'_0$, $F'_1$ and $F'_2$ are shown in Fig. 1. 
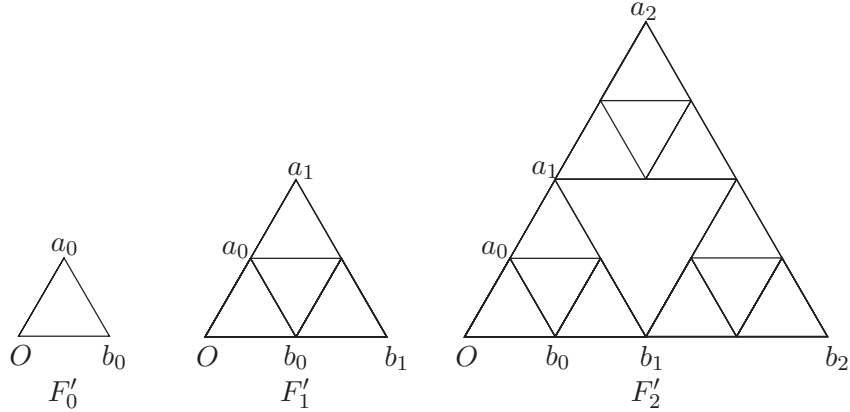
\begin{figure}[htb]
\begin{center}
\input{Fig1.tex}
\\[1\baselineskip]
\caption{$F'_0$, $F'_1$ and $F'_2$.}
\end{center}
\end{figure}

Finally, we let $F'^R_N$ be the reflection od $F'_N$ with respect to the $y$-axis, 
and denote $\dsp F_{0} =\bigcup _{N=1}^{\infty}(F'_N \cup F'^R_N)$;
the graph $F_{0}$ is called the (infinite) {\bf pre-{\sg}}.  $F_0$ is shown in Fig. 2.
\begin{figure}[htb]
\begin{center}
\input{Fig2.tex}
\\[1\baselineskip]
\caption{The pre-{\sg} $F_0$.}
\end{center}
\end{figure}

Furthermore, by letting 
$G_0$ and $E_0$ denote the set of vertices and the set of edges of $F_0$, respectively,
we see that, for each $N \in {\mathbb Z}_+$, $F_{N}=2^N F_{0}$ 
can be regarded as a coarse graph with vertices $G_N=\{2^Nx\ :\ x\in G_0 \}$
and edges $E_N=\{ 2^N (x, y) \ :\ (x, y) \in E_0\}$.
We call 
an upward (closed and filled) triangle which is a translation of $
 \triangle Oa_Mb_M$ and whose vertices are in $G_M$ 
 a {\bf $2^{M}$-triangle}.

\vspace{0.5cm}\parr
\subsection{Paths on the pre-{\sg}s}

Let us denote the set of finite paths on $F_0$ starting at $O$ by
\[W = \{\ w=(w(0), w(1), \cdots , w(n)):
\  w(0) =O, \ (w(i-1), w(i)) \in G_0, \ 1 \leq i \leq n,  \ n \in {\mathbb N}\  \}.\]
This gives the natural definition for 
the length $\ell $ of a path  $w=(w(0), w(1),$ $ \cdots , w(n))\in W$; namely, $\ell (w)=n$. 
 
For a path $w\in W$ and $A \subset G_0$,
 we define the  hitting time of $A$ 
 by
\[T_A(w)=\inf  \{j \geq 0 :\ w(j) \in A\},\]
where we set $\inf \emptyset =\infty $. 
By taking 
$w\in W$ and $M \in \pintegers$, we shall define a recursive sequence 
$\{T_i^M(w)\}_{i=0}^{m}$ of
{\bf hitting times of $G_M$} as follows:
Let 
 $T_{0}^{M}(w)=0$, and for $i\geq 1$, let
\[T_{i}^{M}(w)=\inf \{j>T_{i-1}^{M}(w) :\  w(j)\in G_{M}\setminus
\{w(T_{i-1}^{M}(w))\}\};\]
here we take  $m$ to be the smallest integer such that 
$T_{m+1}^{M}(w)=\infty $. Then 
$T_{i}^{M}(w)$ can be interpreted as being the time (steps) taken for the path $w$ to 
hit vertices in $G_{M}$ for the $(i+1)$-st time, 
under the condition that if $w$ hits the same vertex in $G_{M}$ 
more than once in a row, we count it only once.

Now, we consider two sequences of subsets of $W$ as follows: 
for each $N \in \pintegers$, let the set of paths from $O$ to $a_N$, which  
do not hit any other vertices in $G_N$ on the way, be 
\[W_N=\{ w =(w(0),w(1),\cdots,w(n)) \in W : \  w(T_1^N(w))=a_N, \ n=T_1^N(w) \},\]  
and let the set of paths from from $O$ to $a_N$ that 
hit $b_N$ `once' on the way (subject to the counting rule explained above) be 
\[V_N=\{w= (w(0),w(1),\cdots,w(n))\in W : \   w(T_1^N(w))=b_N,\  w(T_2^N(w))=a_N, 
\ n=T_2^N(w)\}.\]

Then, for a path $w\in W$ and each $M\in {\mathbb N}$, we 
 define the  {\bf coarse-graining map}
$Q_{M}$ by \[(Q_{M}w)(i)=w(T_{i}^{M}(w)), \ \ \mbox{ for } i=0,1,2,\ldots, m,\]
where $m$ is the smallest integer such that 
$T_{m+1}^{M}(w)=\infty $ as above. Thus,
\[Q_M w=(w(T_0^M(w)), w(T_1^M(w)), \ldots , 
w(T_m^M(w)) )\]
is  a path on a coarser graph $F_M$. For 
$w\in W_N\cup V_N$ and $M \leq N$, the end point of the coarse-grained path is $w(T_m^M(w))=a_N$, 
and if we write $(2^{-M}Q_{M}w)(i)=2^{-M}w(T_{i}^{M}(w))$,
then $2^{-M}Q_{M}w$ is a path in $W_{N-M}\cup V_{N-M}$ and $\ell (2^{-M}Q_{M}w)=m$.
In the following, we often write $w(T^M_i)$ instead of $w(T^M_i(w))$.

Define a family of probability measures $P_N$ on $W_N$, $N=1,2,\cdots$ by 
assigning  each $w\in W_N$, 
\[P_N[w]=\biggl(\frac{1}{4}\biggr)^{\ell (w)-1}.\]
$(W_N, P_N)$ defines a family of fixed-end random walks $Z_N$ on $F_N$ 
such that 
\eqnb
\eqna{Z}
Z_N(w)(i)=w(i),\ \  i=0, \cdots , \ell (w), \ \ w\in W_N.\eqne
This is a simple random walk on $F_0$ starting at $O$ and stopped at the first hitting time of $a_N$ conditioned that the walk does not hit any vertices in $G_N \setminus \{O\}$ on the way.  The factor $(1/4)^{-1}$ comes from this conditioning.

Define another family of probability measures $P'_N$ on $V_N$, $N=1,2,\cdots .$ by assigning each $w\in V_N$, 
\[P'_N[w]=\biggl(\frac{1}{4}\biggr)^{\ell (w)-2}.\]
$(V_N, P'_N)$ defines a family of fixed-end random walks $Z'_N$ on $F_0$ 
such that 
\eqnb
\eqna{Zprime}
Z'_N(w)(i)=w(i),\ \  i=0, \cdots , \ell (w), \ \ w\in V_N.\eqne
This is a simple random walk on $F_0$ starting at $O$ and stopped at the first hitting time of $a_N$ conditioned that the walk hits $b_N$ `once' on the way.  

Note that a coarse grained simple random walk is again a simple random walk on a 
coarse graph, that is, for $M<N$, $P_N\circ Q^{-1}_M=P_{N-M}$ and 
$P'_N\circ Q^{-1}_M=P'_{N-M}$.


\section{Loop erasure by the erasing-larger-loops-first rule}
\seca{LE}

For  $(w(0), w(1), \cdots , w(n))\in W_N\cup V_N$, 
if  there are $c \in G_0$,  $i$ and $ j$, $0\leq i < j \leq n$ such that $w(i)=w(j)=c \ $ and 
$w(k)\neq c$ for any $i<k<j$,  
we call the path segment $[w(i),w(i+1), \ldots , w(j)]$ a {\bf loop formed at} ${\bf c}$ and 
define its {\bf diameter} by $d=\max _{i\leq k_1<k_2 \leq j} |w(k_1)-w(k_2)|$, where 
$|\ \cdot \ |$ denotes the Euclidean distance. 
Note that a loop can be a part of another larger loop formed at some other 
vertex. 
By definition, the paths in $W_N\cup V_N$ do not have any loops 
with diameter greater than $2^{N-1}$. 
For each $N \in {\mathbb Z}_+$, let $\Gamma _N$ be the set of loopless paths from $O$ to $a_N$:
\[\Gamma _N= \{\ (w(0), w(1), \cdots , w(n)) \in W_N \cup V_N:
\  w(i) \neq w(j), \ 0\leq  i< j\leq n,\ n\in \nintegers \  \} .
\]
Note that any loopless path in  $\Gamma _N$ is confined in $\triangle Oa_Nb_N$.

We shall now describe the  loop-erasing procedure in a more organized manner than \cite{HM}.
We start by erasing loops from paths in  $W_1\cup V_1$.
\vspace{0.3cm}\parr
{\bf Loop erasure for $W_1\cup V_1$}

\itmb
\item[(i)] Erase all the loops formed at $O$;
\item[(ii)] Progress one step forward along the path, and 
erase all the loops at the new position;
\item[(iii)] Iterate this process, taking another step forward along the path and erasing the 
loops there, until reaching $a_1$.
\itme

Denote the resulting path $Lw$, where $L: W_1\cup V_1 \to \Gamma _1$ is the loop-erasing operator.
Fig. 3 shows all the possible loopless paths from $O$ to $a_1$ on $F_1$. 
Here only the parts in $\triangle Oa_1b_1$ are shown, for 
any path cannot go into the other triangles without making a loop.   
Note that $w \in W_1$ implies $Lw \in W_1\cap \Gamma _1$, but that 
 $w \in V_1$ can result in  $Lw \in W_1 \cap \Gamma _1$, with $b_1$ being erased 
 together with a loop.
So far, our loop-erasing procedure is the same as the chronological method defined for paths on 
${\mathbb Z}^d $ in \cite{Lawler}.

\begin{figure}[htb]
\begin{center}
\input{Fig4.tex}
\\[1\baselineskip]
\caption{Loopless paths from $O$ to $a_1$ on $F_1$.}
\end{center}
\end{figure}
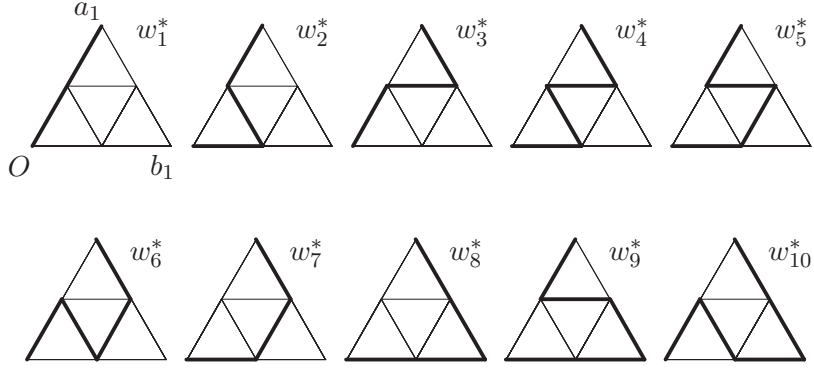

For a general $N$, 
we erase loops from the largest-scale loops down, repeatedly applying the 
loop-erasing procedure for $W_1 \cup V_1$.  To describe the procedure,  
we introduce a `step-based' decomposition of a path based on the self-similarity and the symmetries of the pre-\sg s.  
Assume $w\in W_N\cup V_N$ and $0\leq M< N$. 
Note that the pair of adjacent $2^{M}$--triangles including $(Q_Mw)(i-1)$, $(Q_Mw)(i)$
and $(Q_Mw)(i+1)$ is similar to $F_{0} \cap (\triangle Oa_Mb_M \cup \triangle Oa^R_Mb^R_M )  $,  
where $\triangle Oa^R_Mb^R_M$ is the reflection of $\triangle Oa_Mb_M$ with regard to the $y$--axis. 
This leads to a unique decomposition: 
\eqnb \eqna{decomposition1}
(\tilde{w}; w_1, \cdots , w_{\ell (\tilde{w})}), \ \tilde{w} \in W_{N-M}\cup V_{N-M}, \ w_i \in 
W_{M}, \ i=1, \cdots , \ell (\tilde{w})
\eqne
such that $\tilde{w}$ is similar to $Q_Mw$ and that 
the path segment  $(w(T_{i-1}^M(w)), w(T_{i-1}^M(w)+1)), \cdots , $ $w(T_{i}^M(w)))$ of 
$w$ is identified with $w_i \in W_{M} $ by appropriate  
rotation, translation and reflection so that $w(T_{i-1}^M(w))$ is identified with 
$O$ and $w(T_{i}^M(w))$ with $a_{M}$. 
We shall use this kind of identification throughout the paper.
We illustrate a simple example of the decomposition for $N=2$ and $M=1$ in Fig. 4.

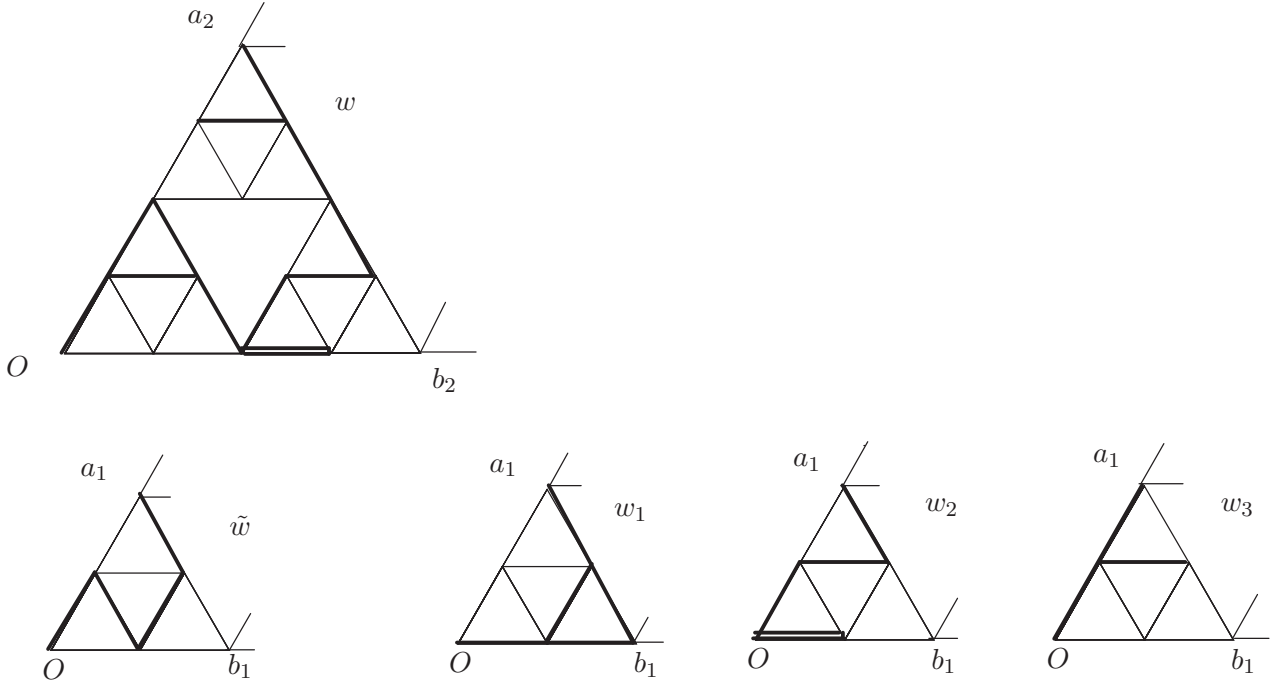
\begin{figure}[htb]
\begin{center}
\input{Fig3.tex}
\\[1\baselineskip]
\caption{$w, \tilde{w}, w_1, w_2, w_3$.}
\end{center}
\end{figure}

\vspace{0.3cm}\parr
{\bf Erasure of the largest loops}

\itmb
\item[(1)]
Decompose a  path $w\in W_N \cup V_N$ into $(\tilde{w}; w_1, \cdots , w_{\ell (\tilde{w})})$, $\tilde{w}=2^{-(N-1)}Q_{N-1}w \in  W_1\cup V_1$, 
$w_i \in W_{N-1} $ 
$i=1, \cdots , \ell(\tilde{w})$ as in \eqnu{decomposition1} with $M=N-1$.  Fig. 5(a) shows the original $w$ and 
Fig. 5(b) shows $Q_{N-1}w$.  
\item[(2)]
Erase all the loops 
from $\tilde{w}$ following the loop-erasure for $W_1\cup V_1$ 
to obtain
$L\tilde{w} \in \Gamma_1$. 
Denote the coarse, loopless path $2^{(N-1)}L\tilde{w} $ on $F_{N-1}$ by $\hat{Q}_{N-1}w$ (Fig. 5(c)).
\item[(3)]
Restore the original fine structures to the remaining parts as shown in Fig. 5(d) to 
obtain a path $w' \in W_N\cup V_N$.  To 
be more precise, if  we write  $\hat{Q}_{N-1}w=(w(T_0^{N-1}), w(T_{s_1}^{N-1}), \cdots ,$ $w(T_{s_n}^{N-1}) )$ , then for each $i$, between 
$w(T_{s_i}^{N-1})$ and $w(T_{s_{i+1}}^{N-1})$, insert the path segment 
$w_{s_{i}+1} =(w(T_{s_i}^{N-1}), w(T_{s_i}^{N-1}+1), \cdots, w(T_{s_{i}+1}^{N-1}))$
chosen from the original decomposition in Step (1). Note that 
$Q_{N-1}w'=\hat{Q}_{N-1}w$ holds.
\itme

In this stage all the loops with diameter greater than $2^{N-2}$ have been erased.
We repeat Procedure (1)--(3) within each $2^{N-1}$--triangle to erase all the 
loops with diameter greater than $2^{N-3}$, and then within 
each $2^{N-2}$--triangle, and so on, until there remain no loops.

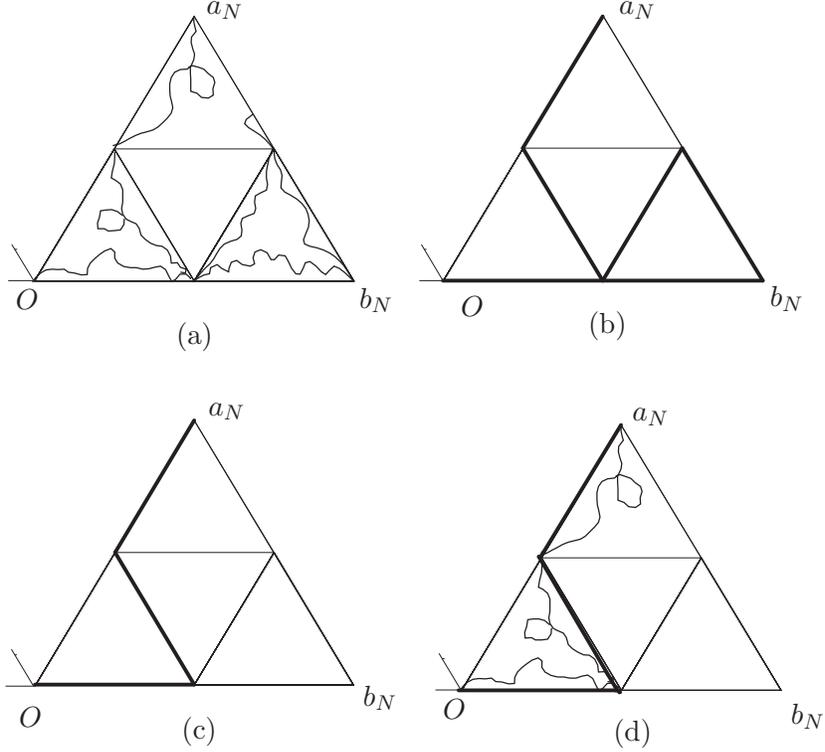
\begin{figure}[htb]
\begin{center}
\input{Fig5.tex}
\\[1\baselineskip]
\caption{The loop-erasing procedure: (a) $w$, (b) $Q_{N-1}w$, (c) $\hat{Q}_{N-1}w$, (d) 
fine structures restored.}
\end{center}
\end{figure}


To describe the procedure more precisely, we prepare another kind of decomposition, 
a `triangle-based' decomposition.
For $w \in  W_N$ and $0\leq M \leq N$, we shall define 
the  sequence $(\Delta_1, \ldots , \Delta_k)$ of  the  $2^{M}$--triangles 
$w$ `passes through',  and 
their exit times $\{T_i^{ex,M}(w)\}_{i=1}^{k}$ 
as a subsequence of  $\{T_i^{M}(w)\}_{i=1}^{m}$ as follows:
Let  $T_0^{ex, M}(w)=0$. 
There is a unique  $2^{M}$--triangle that contains $w(T_0^M)$ and $w(T_1^M)$, 
which we denote by  $\Delta _1$. 
For $i\geq 1$, define 
\[J(i)=\min \{j \geq 0\ : \ j<m,\ T_j^M(w)>T_{i-1}^{ex, M}(w),\ 
w(T_{j+1}^M(w))\not\in \Delta _i\},\] 
if the minimum exists, 
otherwise $J(i)=m$.   
Then define $ T_i^{ex, M}=T_i^{ex, M}(w)=T_{J(i)}^{M}(w)$, and let $\Delta _{i+1}$ be
the unique $2^{M}$--triangle that contains both $w(T_i^{ex, M})$ and $w(T_{J(i)+1}^{M})$. 
By definition, we see that  $\Delta _{i} \cap 
\Delta _{i+1 }$ is a one-point set $\{w(T_i^{ex, M})\}$, for $i=1, \ldots , k-1$. 
We denote the sequence of these triangles by 
$\sigma _M(w)=( \Delta _1, \ldots , \Delta _k)$, 
and call it the ${\bf 2^{M}}${\bf --skeleton} of $w$.
We call the sequence  
$\{T_i^{ex, M}(w)\}_{i=0}^{k}$ 
{\bf exit times} from the triangles in the skeleton.  
For each $i$, there is an $n=n(i)$ such that 
$T_{i-1}^{ex, M}(w)=T^M_{n}(w)$.  If $T_{i}^{ex, M}(w)=T^M_{n+1} (w)$,
we say that $\Delta _i \in \sigma _M(w)$ is  {\bf Type 1},  
and if  $T_{i}^{ex, M}(w)=T^M_{n+2}(w)$,  {\bf Type 2}.
For $w\in W_N \cup V_N$ and $M<N$, if $Q_Mw$ is similar to a path in $ \Gamma_{N-M}$, namely, 
 $2^{-M}Q_Mw \in \Gamma_{N-M}$, then 
its $2^{M}$--skeleton is a collection of distinct $2^{M}$--triangles 
and each of them is either Type 1 or Type 2.  

Assume $w\in W_N\cup V_N$ and $M\leq N$.
For each $\Delta $ in $\sigma _M(w)$, the {\bf path segment of} ${\bf w} $ {\bf in}
 ${\bf \Delta }$ is 
defined by 
\eqnb
\eqna{segment}
w|_{\Delta }=[w(n), \ T^{ex, M}_{i-1}(w) \leq n \leq T^{ex, M}_{i}(w)].
\eqne 
Note that the definition of $T^{ex, M}_i (w)$ allows a path segment $w|_{\Delta }$ to leak into the 
neighboring 
$2^{M}$--triangles.  
If $Q_Mw $ is similar to a path in $\Gamma _{N-M}$,  then 
$w|_{\Delta } \in W_{M}$ or $w|_{\Delta } \in V_{M}$ (identification implied),  
according to the type of $\Delta \in \sigma_M(w)$, 
where the entrance to $\Delta $ is identified with $O$ and the exit with $a_M$. 
This means that each $w $ such that 
$Q_Mw $ is similar to a path in $\Gamma_{N-M}$ 
can be decomposed uniquely to 
\eqnb 
\eqna{decomposition2}
(\sigma_M(w); \ w|_{\Delta_1}, \cdots , w|_{\Delta_k}), \ \ 
w|_{\Delta _i } \in W_{M}\cup V_{M},\ i=1, \cdots , k.
\eqne
We call a loop $[w(i), w(i+1), \cdots , w(i+i_0)]$ a ${\bf 2^{M}}${\bf -scale loop} whenever 
there exists an $M \in \pintegers $ such that  
\[\max \{N' : w(i)=w(i+i_0)\in G_{N'}\}=M,\ d\geq 2^{M},\]
where $d$ is the diameter of the loop. 

\vspace{0.3cm}\parr
{\bf Induction step of loop erasure}

Let $w\in W_N \cup V_N$ and 
$1\leq M\leq N$.  
Assume that all of the $2^{N-1}$ to $2^{N-M}$--scale loops have been erased from  
$w$, and denote the path obtained at this stage by $w' \in W_N\cup V_N$.  
Note that $Q_{N-M}w' $ is similar to a path in $\Gamma _M$.

\itmb
\item[1)]
Decompose $w'$ to obtain  
 $(\sigma _{N-M}(w') ; \ w'_1, \cdots w'_k )$,
$w'_i \in W_{N-M} \cup V_{N-M}$ as
 given in \eqnu{decomposition2}.

\item[2)] 
From each $w_i' $, 
erase $2^{N-M-1} $--scale 
loops (largest-scale loops) according to the base step procedure (1)--(3) above to obtain $\tilde{w}_i' \in W_{N-M}\cup V_{N-M}$.

\item[3)]
Assemble $(\sigma _{N-M}(w'); \ \tilde{w}'_1, \cdots , \tilde{w}'_k)$ 
to obtain $w'' \in W_N \cup V_N$, which is determined uniquely.
$w''$ has no 
$2^{N-1}$ to $2^{N-M-1}$--scale loops. 
\itme
\QED

We repeat 1)--3) until we have no loops and denote the resulting loopless path $Lw \in \Gamma_N$.
In this way, the loop erasing operator $L$, first defined for $W_1\cup V_1$, has been extended to 
$L:\bigcup _{N=1}^{\infty } (W_N\cup V_N) \to \bigcup _{N=1}^{\infty }\Gamma _N$ 
with $L(W_N\cup V_N)=\Gamma _N$.
Note that the operation described above is essentially a repetition of loop-erasing  for $W_1\cup V_1$.

We induce measures $\hat{P}_N=P_N\circ L^{-1}$ and ${\hat{P}}'_N={P}'_N\circ L^{-1}$, which satisfy 
$\hat{P}_N[\Gamma_N]={\hat{P}}'_N [\Gamma _N]=1$.  
For $w^*_1, \cdots , w^*_{10}$ shown in Fig. 3, 
denote 
\[
p_i= \hat{P}_1 [w^*_i] =P_1[w: Lw=w^*_i], \ \ 
q_i= {\hat{P}}'_1 [w^*_i] ={P}'_1[w: Lw=w^*_i].
\] 
They were obtained in \cite{HM} by direct calculation: 
  
\eqnb \eqna{pvalue} p_1=1/2, \ \ 
p_2=p_3=p_7=2/15, \ \ p_4=p_5=p_6=1/30, \ \ p_8=p_9=p_{10}=0, \eqne 
\eqnb \eqna{qvalue} q_1=1/9,\ \ q_2=q_3=11/90, \ \ q_4=q_5=q_6=2/45, \ \ q_7=8/45, \ \ q_8=2/9, \ \  
q_9=q_{10}=1/18.\eqne

 $\hat{P}_N$ and  ${\hat{P}}'_N$ define two kinds of walks $Y_N=LZ_N$ 
and $Y'_N=LZ'_N$ on $F_0 \cap \triangle Oa_Nb_N$ obtained 
by erasing loops from the simple random walks $Z_N$ and  ${Z}'_N$, respectively. 
We remark that $\displaystyle \frac{2}{3}
\hat{P}_N+\frac{1}{3}{\hat{P}}'_N$ equals to the `standard' LERW studied 
in \cite{STW}.


For $w \in W_N\cup V_N$, we defined $\hat{Q}_{N-1}w$ in Step (2) for the erasure  of the largest-scale loops.  For 
later use we define $\hat{Q}_{N-K}w$ on $F_{N-K}$ for all $K=0,1, \cdots , N$.
Repeat the induction step 1)--3) $K$ times to have down to $2^{N-K}$--scale loops erased and denote the resulting path $w'$. Let   $\hat{Q}_{N-K}w=
Q_{N-K}w'$, namely the coarse path before restoring fine structures.  In particular, $\hat{Q}_Nw =Q_Nw$ and $\hat{Q}_0w=Lw$. By construction,  the distributions of 
$2^{-(N-K)}\hat{Q}_{N-K}Z_N$ and  $2^{-(N-K)}\hat{Q}_{N-K}Z'_N$ equal to  $\hat{P}_K$ and $\hat{P}'_K$, respectively.

\vspace{0.5cm}\parr
\section{Asymptotic behavior of the exit times}
\seca{Steps}

In this section, we look into the asymptotics of exit times $T^{ex, N}_1(Y_N)$ and 
$T^{ex, N}_1(Y'_N)$ as $N \to \infty$, 
which will be used in \secu{Proof}. 
 
For $w\in \Gamma _N$, let us denote  the number of $2^0$-- triangles of Type 1
(the path passes two of the vertices) and
those of Type 2 (the path passes all three vertices) in $\sigma _0(w)$  by $s_1(w)$ and $s_2(w)$, respectively.
Note that $T_1^{ex, N}(w)=\ell (w)=s_1 (w)+2 s_2(w)$. 
Define two sequences, $\{\Phi ^{(1)}_N\}_{N \in {\mathbb N}}$ and $\{\Phi^{(2)}_N\}_{N \in {\mathbb N}}$, of 
generating functions by:  
\[\Phi^{(1)}_N (x,y)=\sum _{w\in \Gamma _N }\hat{P}_N(w)x^{s_1(w)}y^{s_2(w)},\]
\[\Phi^{(2)}_N (x,y)=\sum _{w\in \Gamma _N }{\hat{P}}'_N(w)x^{s_1(w)}y^{s_2(w)},\ \ \ x,y\geq 0.\]
For simplicity, we shall denote $\Phi^{(1)}_1 (x,y)$ and $\Phi^{(2)}_1 (x,y)$ by 
$\Phi^{(1)} (x, y)$ and $\Phi^{(2)}(x,y)$. 
A crucial observation is that in the process of erasing loops from 
$Z_{N+1}$, if we stop at the point where we have obtained
$\hat{Q}_1Z_{N+1}$ after erasing down to $2^1$--scale loops, it is nothing but the procedure for obtaining $LZ_N$ from
$Z_N$, namely, the distribution of $2^{-1}\hat{Q}_1Z_{N+1}$ equals to  $\hat{P}_N$.  The same holds for ${Z}'_{N+1}$ as well.  This 
combined with \eqnu{pvalue} and \eqnu{qvalue} leads to the recursion relations for the generating functions given below:

\prpb (Proposition 3 in \cite{HM})
 \prpa{recursion}

The above generating functions satisfy the following recursion relations for all $N \in {\mathbb N}$ :
\[\Phi ^{(1)}(x,y)= \frac{1}{30}(15 x^2 +8 xy + y^2 +2 x^2y + 4 x^3),
\]
\[\Phi ^{(2)} (x,y)=\frac{1}{45}(5 x^2 + 11 xy + 2y^2 + 14 x^2y + 8 x^3+ 5x y^2) ;\]
\[\Phi ^{(i)}_{N+1}(x,y)=\Phi ^{(i)}_{N}(\Phi  ^{(1)} (x,y), \Phi^{(2)} (x,y)),\ \ i=1,2 .\]
\prpe

\vspace{0.5cm}\par
Define the mean matrix by
\eqnb
\eqna{matrix}
{\bf M}=\left [
\begin{array}{cc}
\displaystyle \frac{\partial}{\partial x}\Phi^{(1)} (1,1) 
& \displaystyle \frac{\partial}{\partial y}\Phi ^{(1)}(1,1) 
\\
\displaystyle \frac{\partial}{\partial x}\Phi^{(2)} (1,1) 
& \displaystyle \frac{\partial}{\partial y}\Phi^{(2)} (1,1)
\end{array}
\right ]
=\left [
\begin{array}{cc}
\displaystyle \frac{9}{5} & \displaystyle\frac{2}{5} 
\\
\\
\displaystyle \frac{26}{15} & \displaystyle\frac{13}{15}
\end{array}
\right ].
\eqne
It is a strictly positive matrix, and the 
larger eigenvalue is given by 
 $\lambda =(20+\sqrt{205})/15 =2.2878 \ldots $.
The following is a restatement of Proposition 9 in \cite{HM}.


\prpb
\prpa{supbr}
\itmb
\item[(1)]
Let $G_N^{(1)}(t)$ and $G^{(2)}_N(t)$ be the Laplace transforms of  $\lambda ^{-N}T_1^{ex, N}(Y_N)$ 
and $\lambda ^{-N}T_1^{ex, N}(Y'_N)$, respectively, that is,     
\[G^{(1)} _N(t)=\hat{E}_N[\exp ( -t\lambda ^{-N}T_1^{ex, N}(w))],\]
\[G^{(2)}_N(t)=\hat{E}'_N[\exp ( -t\lambda ^{-N}T_1^{ex, N}(w))],\ \ \ t\in {\mathbb C}\]
where $\hat{E}_N$ and $\hat{E}'_N$ are expectations with regard to 
$\hat{P}_N$ and $\hat{P}'_N$, respectively. Then they are expressed in terms of 
the generating functions as 
\eqnb \eqna{GPhi}
G^{(i)}_N(t)=\Phi ^{(i)}_N(e^{-\lambda ^{-N}t}, \ e^{-2\lambda ^{-N}t } )\ \ i=1,2.
\eqne

\item[(2)]
$\lambda ^{-N}T_1^{ex, N}(Y_N)$ and $\lambda ^{-N}T_1^{ex, N}(Y'_N)$ 
converge  in law to some integrable random variables $T^*_1$ and $T^*_2$, respectively, as 
$N \to \infty$.
$T^*_1$ and $T^*_2$ have strictly positive probability density functions on $(0, \infty )$.

\item[(3)]
Let $g_i(t)$ be the Laplace transform of $T^*_i$.
For each $i$, $G_N^{(i)}(t)$ converges to $g_i(t) $ uniformly on any compact set in ${\mathbb C}$ as $N \to \infty$. 
$g_1(t)$ and $g_2(t)$  
are entire functions on ${\mathbb C}$ and the unique solution to 
\[g _1 (\lambda  t)=\Phi ^{(1)}( g _1( t), g _2 ( t)),\  \  
g _2 (\lambda  t)=\Phi ^{(2)} ( g _1( t), g _2( t)),\
g _1(0)=g _2(0)=1.
\]

\itme
\prpe


To obtain the left tail behavior of the scaled exit times, the following Tauberian theorem has a most suitable form. 

\thmb
\thma{Tauberian} (Theorem 5.9 in \cite{tets})

Let
$\mu_N$, $N\in {\mathbb N}$ be a family of  probability measures on $[0,\infty)$ and let $\displaystyle G_N(s)=\int _0^{\infty} e^{-sx} \mu _N(dx) $, $s>0$ be their Laplace transforms.  
If there exist positive constants $C_{4.1}$-- $C_{4.4}$,  $s_0>0$, $s_1\in {\mathbb R}$ and $0<\nu <1$ 
 such that 
\[C_{4.1}\exp (-C_{4.2}s^{\nu})\leq G_N(s) \leq C_{4.3}\exp (-C_{4.4}s^{\nu}), \]
holds for all  $s>s_0$ and $N> s_1+\frac{\nu}{\log 2} \log s $, 
then the following holds.

\itmb
\item[(1)]
There exist positive constants $C_{4.5}$ and $C_{4.6}$ such that for any 
positive sequence satisfying 
\parr
$\displaystyle \lim_{N \to \infty } 2^{N(1-\nu)/\nu}\alpha _N
 =\infty$ and $\displaystyle \lim_{N \to \infty}\alpha _N=0$, the following holds:  

\begin{eqnarray*}
-C_{4.5} &\leq& \liminf _{N\to \infty} \alpha _N^{\nu /(1-\nu) } \log \mu _N ([0, \alpha _N])\\ 
&\leq&
\limsup _{N\to \infty} \alpha _N^{\nu /(1-\nu) } \log \mu _N ([0, \alpha _N]) \leq -C_{4.6}.
\end{eqnarray*}

\item[(2)]
There exist positive constants $C_{4.7}$ -- $C_{4.9}$
such that for any $\xi  >0$ and $N \in {\mathbb N} $ satisfying $\displaystyle (2^{\frac{1}{\nu} -1})^N\xi \geq C_{4.7}$,
\[\mu _N ([0, \xi ]) \leq C_{4.8} e^{- C_{4.9} \xi ^{-\nu /(1-\nu) }} \]
holds.
\itme

\thme

(1) is a kind of restatement of a Tauberian theorem of exponential type given in \cite{Kasahara} and \cite{KasaharaKosugi}, and (2) is the combination of Chebyshev's inequality and (1).

\prpb
\prpa{G}
For $t>0$, 
$G^{(1)} _N(t)$ and $G^{(2)} _N(t)$ satisfy the condition for \thmu{Tauberian} 
with $\nu =\log 2/\log \lambda$. 
\prpe

\prfb

Using \eqnu{GPhi}, we rewrite the recursion as 
\eqnb \eqna{GPhiN}  G^{(i)}_{N+1}(t)=\Phi^{(i)}(G^{(1)}_N(t/\lambda), G^{(2)}_N(t/\lambda)), \ \ i=1,2.\eqne
From the explicit form of $\Phi^{(i)}$ in \prpu{recursion}, we have for $0<x, y<1$,
\[ q_1(x \wedge y)^2\leq \Phi^{(i)}(x,y) \leq (x\vee y)^2,\ \ i=1, 2,\]
where $q_1=1/9$. 
Repeating this $M$ times, we have 
\eqnb \eqna{2M} 
\{q_1 (x \wedge y)\}^{2^M} \leq \Phi^{(i)}_M(x,y) \leq (x\vee y)^{2^M},\ \ i=1, 2.\eqne
This combined with  \eqnu{GPhiN} gives  
\eqnb \eqna{q1}\{q_1(G^{(1)}_N(t/\lambda^M) \wedge G^{(2)}_N(t/\lambda^M)) \}^{2^M}\leq G_{N+M}^{(i)}
(t) \leq\{G^{(1)}_N(t/\lambda^M) \vee G^{(2)}_N(t/\lambda^M)\}^{2^M}
\eqne
Fix $t_0>0$ arbitrarily.  Since $\{G^{(1)}_N(t_0) \vee G^{(2)}_N(t_0)\}_{N=1}^{\infty}$ 
and $\{(G^{(1)}_N(\lambda t_0) \wedge G^{(2)}_N(\lambda t_0) \}_{N=1}^{\infty}$ 
are positive convergent sequences by \prpu{supbr} 
(3), there exist constants $c_1, c_2 \in (0,1)$ 
such that 
\eqnb \eqna{G1G2} q_1(G_N^{(1)}(\lambda t_0) \wedge G_N^{(2)}(\lambda t_0))>c_1,\ \ 
G_N^{(1)}(t_0)\vee G_N^{(2)}(t_0) <c_2, \eqne
for all $N \in {\mathbb N}$.
For any $t>t_0$, choose $M \in {\mathbb Z}_+$ such that 
\eqnb \eqna{chooseM} \lambda ^M\leq \frac{t}{t_0} <\lambda ^{M+1}.\eqne
Then, the monotonicity of $G_N^{(i)}$ combined with \eqnu{q1}, \eqnu{G1G2} and \eqnu{chooseM} gives  
\[\ c_1^{2^M}\leq G^{(i)}_{N+M}(t) \leq c_2^{2^M},\ \ \ i=1,2. \]
This further leads to
\[\exp (-C_{4.2} t^{\nu}) \leq  G^{(i)}_{N}(t) \leq \exp (-C_{4.4} t^{\nu}), \ \ i=1,2 \]
for all $t>t_0$ and $N>\log _{\lambda}(t/t_0)$, where we 
put $\dsp C_{4.2}=-\frac{\log c_1}{t_0^{\nu}}$ and $\dsp C_{4.4}=-\frac{\log c_2}{2t_0^{\nu}}$.

\QED
\prfe


\vspace{0.5cm}\parr
\section{Extention to the infinite \sg}
\seca{Extention}

In this section, we show that the loop-erased random 
walks defined in \secu{LE} can be extended to a loop-erased random walk on the infinite \sg .
For this purpose, we need walks from $O$ to $b_N$ as well as those from $O$ to $a_N$.
For each $N \in {\mathbb Z}_+$, let 
\[W_N^{b}=\{ w =(w(0),w(1),\cdots,w(n)) \in W : \  w(T_1^n(w))=b_N, \ n=T_1^N(w) \},\]  
\[V_N^{b}=\{w= (w(0),w(1),\cdots,w(n))\in W : \ 
 w(T_1^N(w))=a_N,\  w(T_2^N(w))=b_N, \ n=T_2^N(w)\}.\] 
and probability measures $P^{(2)}_N$ on $W_N^{b}$ and 
$P^{(4)}_N$ on $V_N^{b}$ by 
\[P_N^{(2)}[w]=\biggl(\frac{1}{4}\biggr)^{\ell (w)-1}, \ \ \ w \in W_N^{b},\]
\[P_N^{(4)}[w]=\biggl(\frac{1}{4}\biggr)^{\ell (w)-2}, \ \ \ w \in V_N^{b}.\]
Let $U_N=W_N\cup V_N \cup W_N^b \cup V_N^b$
and extend the loop-erasing operator $L$ to $\bigcup _{N=1}^{\infty} U_N$.  
Denote $P^{(1)}_N=P_N$,  $P^{(3)}_N=P'_N$ and 
$\hat{P}^{(i)}_{N}=  P^{(i)}_N\circ L^{-1}$, for $i=1,2,3,4$. 
In the rest of the paper,  we use the same notation $\Gamma_N$ for loopless paths in $U_N$.
Let 
\[\Omega =\{  \omega=(\omega_0, \omega_1, \omega_2, \cdots )\ :\ \omega _0 \in \Gamma _0, \
\omega _N \in \Gamma _N, \ 
\omega_N |_{N-1} =\omega_{N-1}, \ N \in {\mathbb N} \}, \]
where $\omega_N |_{N-1}$ denotes the path $\omega _N$  stopped at $T_1^{ex, N-1}(\omega _N)$ and 
${\cal B}$ the $\sigma $-algebra on $\Omega $ generated by cylinder sets.
Define the projection onto the first $N+1$ elements by 
\[\pi _N \omega =(\omega_0, \omega_1, \ldots , \omega_N). \]
and a probability measure $\tilde{P}_N$ on $\pi _N\Omega $ by 
\eqnb
\eqna{kumiawase}
\tilde{P}_N[(\omega_1, \ldots , \omega_N)]=\frac{11}{28}( \hat{P}_N^{(1)}[\omega_N] +\hat{P}_N^{(2)}[\omega_N]) +\frac{3}{28}( \hat{P}_N^{(3)}[\omega _N] +\hat{P}^{(4)}_N[\omega _N]).
\eqne

\prpb
\prpa{consistencycond}
The sequence $\{\tilde{P}_N\}$, $N\in {\mathbb Z}_+$ defined in \eqnu{kumiawase}  satisfies:
\eqnb
\eqna{consistency}
\tilde{P}_N [(\omega_0, \omega_1, \ldots , \omega_N) ]=\sum _{\omega '} \tilde{P}_{N+1} [(\omega_0, \omega_1, \ldots , \omega_N, \omega ' ) ],
\eqne 
where the sum is taken over all possible $\omega ' \in \Gamma _{N+1}$ such that $\omega '|_N =\omega _N$.
\prpe

\prfb
Assume $u \in U_{N+1}$. 
Recall that in Step (2) of erasing the largest-scale loops, namely,  
$2^N$--scale loops, from $u$, we obtain  $\hat{Q}_Nu$, which  satisfies $2^{-N}\hat{Q}_Nu \in\Gamma_1$ and whose law under $\hat{P}^{(i)}_{N+1}$ is equal to  $\hat{P}^{(i)}_1$.
Let $\Delta _0 =\triangle Oa_0b_0$ and denote the path segment of 
$2^{-N}\hat{Q}_Nu$ 
in $\Delta _0$ by $u_1:=(2^{-N}\hat{Q}_Nu )|_{\Delta _0}$.  Then $u_1 \in \Gamma _0 =\{(O, a_0), (O, b_0), (O, b_0, a_0), (O, a_0, b_0)\}$.  Denote $v_1^*=(O, a_0), v_2^*=(O, b_0), v_3^*=(O, b_0, a_0)$,  
$v_4^*=(O, a_0, b_0),$ and $\Delta = \triangle Oa_Nb_N$.  
For $\hat{w}\in \Gamma_N$, we classify the event $\{ u\in U_{N+1}: Lu|_{\Delta}=\hat{w}\}$ by $u_1$.  Note that under the condition that  $u_1=v_j^*$, 
the distribution of $ Lu|_{\Delta }$ is equal to $\hat{P}_{N}^{(j)}$.
Thus, for $i=1, 3$, 
\[ \begin{array}{l} \dsp
\hat{P}_{N+1}^{(i)}[ \ w \in \Gamma _{N+1} \ :\ w|_{\Delta }=\hat{w}\ ] =
P_{N+1}^{(i)}[\ u \in U _{N+1}  \ :\ Lu|_{\Delta }=\hat{w}\ ] \\ 
\dsp\phantom{\hat{P}_{N+1}^{(i)}}
=\sum_{j=1}^{4} 
P_{N+1}^{(i)} [\ Lu|_{\Delta }=\hat{w} \mid u_1=v^*_j \ ] \ P_{N+1}^{(i)}[\  u_1=v^*_j \ ]\\ 
\dsp\phantom{\hat{P}_{N+1}^{(i)}}
=\sum_{j=1}^{4} 
\hat{P}_{N}^{(j)}[\hat{w} ] \ \hat{P}_{1}^{(i)} [\ v \in \Gamma_1\ :\ v|_{\Delta_0} =v_j^* \ ] \\ 
\dsp\phantom{\hat{P}_{N+1}^{(i)}}
=\hat{P}_{N}^{(1)}[\hat{w} ] \ \hat{P}_{1}^{(i)}[\{w_1^*, w_3^*\}]+
\hat{P}_{N}^{(2)}[\hat{w} ] \ \hat{P}_{1}^{(i)}[\{w_5^*, w_7^*, w_8^*, w_9^*\}]\\
\\
\dsp\phantom{\hat{P}_{N+1}^{(i)}=}
+\hat{P}_{N}^{(3)}[\hat{w} ] \ \hat{P}_{1}^{(i)}[\{w_2^*, w_4^*\}]+
\hat{P}_{N}^{(4)}[\hat{w} ] \ \hat{P}_{1}^{(i)}[\{w_6^*, w_{10}^*\}].
\end{array}\]
Thus, we have 
\[\hat{P}_{N+1}^{(1)}[\ w \in \Gamma _{N+1} \ :\ w|_{\Delta }=\hat{w}\ ] =\frac{19}{30}\hat{P}_{N}^{(1)}[\hat{w} ] +
\frac{1}{6}\hat{P}_{N}^{(2)}[\hat{w} ] +\frac{1}{6}\hat{P}_{N}^{(3)}[\hat{w} ] +\frac{1}{30}\hat{P}_{N}^{(4)}[\hat{w} ] ,
\] 
\[\hat{P}_{N+1}^{(3)}[\ w \in \Gamma _{N+1} \ :\ w|_{\Delta }=\hat{w}\ ] =
\frac{7}{30}\hat{P}_{N}^{(1)}[\hat{w} ] +
\frac{1}{2}\hat{P}_{N}^{(2)}[\hat{w} ] +\frac{1}{6}\hat{P}_{N}^{(3)}[\hat{w} ] +\frac{1}{10}\hat{P}_{N}^{(4)}[\hat{w} ] ,
\] 

For $i=2$, let $\hat{w}^R$ and $v_i^{*R}$ be the paths obtained by 
reflection of $\hat{w}$ and $v_i^{*}$ with regard to  the line $y=x$, respectively. Then   
we have

\[ \begin{array}{l} \dsp
\hat{P}_{N+1}^{(2)}[\ w \in \Gamma _{N+1} \ :\ w|_{\Delta }=\hat{w}\ ] = 
P_{N+1}^{(2)}[\ u \in U _{N+1}  \ :\ Lu|_{\Delta }=\hat{w}\ ] \\ 
\dsp\phantom{\hat{P}_{N+1}^{(2)}}
=\sum_{j=1}^{4} 
P_{N+1}^{(2)}[\ Lu|_{\Delta }=\hat{w} \mid u_1=v^*_j \ ] \ P_{N+1}^{(2)}[\  u_1=v^*_j \ ]\\ 
\dsp\phantom{\hat{P}_{N+1}^{(2)}}
=\sum_{j=1}^{4} 
P_{N+1}^{(1)}[ \ Lu|_{\Delta }=\hat{w}^R \mid u_1=v^{*R}_j \  ] \ P_{N+1}^{(1)}[\  u_1=v^{*R}_j \ ]\\  
\\
\dsp\phantom{\hat{P}_{N+1}^{(2)}}
=\hat{P}_{N}^{(2)}[\hat{w}^R ] \ \hat{P}_{1}^{(1)} [\ v \in \Gamma_1\ :\ v|_{\Delta_0} =v_2^* \ ]
+\hat{P}_{N}^{(1)}[\hat{w}^R ] \ \hat{P}_{1}^{(1)} [\ v \in \Gamma_1\ :\ v|_{\Delta_0} =v_1^* \ ]\\
\\
\dsp\phantom{\hat{P}_{N+1}^{(2)}=}
+\hat{P}_{N}^{(4)}[\hat{w}^R ] \ \hat{P}_{1}^{(1)} [\ v \in \Gamma_1\ :\ v|_{\Delta_0} =v_4^* \ ]
+\hat{P}_{N}^{(3)}[\hat{w}^R ] \ \hat{P}_{1}^{(1)} [\ v \in \Gamma_1\ :\ v|_{\Delta_0} =v_3^* \ ]\\
\\
\dsp\phantom{\hat{P}_{N+1}^{(2)}}
=\hat{P}_{N}^{(1)}[\hat{w} ] \ \hat{P}_{1}^{(i)}[\{w_5^*, w_7^*, w_8^*, w_9^*\}]+
\hat{P}_{N}^{(2)}[\hat{w} ] \ \hat{P}_{1}^{(i)}[\{w_1^*, w_3^*\}]\\
\\
\dsp\phantom{\hat{P}_{N+1}^{(2)}=}
+\hat{P}_{N}^{(3)}[\hat{w} ] \ \hat{P}_{1}^{(i)}[\{w_6^*, w_{10}^*\}]+
\hat{P}_{N}^{(4)}[\hat{w} ] \ \hat{P}_{1}^{(i)}[\{w_2^*, w_4^*\}].
\end{array}\]
Thus,  
\[\hat{P}_{N+1}^{(2)}[\ w \in \Gamma _{N+1} \ :\ w|_{\Delta }=\hat{w}\ ] =\frac{1}{6}\hat{P}_{N}^{(1)}[\hat{w} ] +
\frac{19}{30}\hat{P}_{N}^{(2)}[\hat{w} ] +\frac{1}{30}\hat{P}_{N}^{(3)}[\hat{w} ] +\frac{1}{6}\hat{P}_{N}^{(4)}[\hat{w} ] ,
\] 
Similarly, we have 
\[\hat{P}_{N+1}^{(4)}[\ w \in \Gamma _{N+1} \ :\ w|_{\Delta }=\hat{w}\ ] =
\frac{1}{2}\hat{P}_{N}^{(1)}[\hat{w} ] +
\frac{7}{30}\hat{P}_{N}^{(2)}[\hat{w} ] +\frac{1}{10}\hat{P}_{N}^{(3)}[\hat{w} ] +\frac{1}{6}\hat{P}_{N}^{(4)}[\hat{w} ] .
\] 
Thus,  we see that  
\[(\alpha_1, \alpha _2 , \alpha _3, \alpha _4) =\biggl(\frac{11}{28}, \frac{11}{28}, \frac{3}{28},  
\frac{3}{28}\biggr) \]
is the unique choice that satisfies  
\[\sum _{i=1}^4 \alpha_i \hat{P}_{N+1}^{(i)}[w|_{\Delta} =\hat{w}] =
\sum _{i=1}^4 \alpha_i \hat{P}_N^{(i)}[\hat{w}]\]
for every $\hat{w} \in \Gamma _N$, $N \in {\mathbb N}$. 
\QED
\prfe

\prpu{consistencycond} provides a consisitency condition for Kolmogorov's extension theorem, and we have   
the unique probability measure $P$ on $(\Omega , {\cal B})$, 
 such that 
\[P\circ \pi _N^{-1}=\tilde{P}_N.\]

For any $n \in {\mathbb N}$, take an $N$ satisfying $n \leq 2^N$, then 
the distribution of 
the first $n$ steps of the path, $\omega _N|_n$ is uniquely determined 
independently of $N$.
$(\Omega , {\cal B}, P)$ defines a loop-erased random walk $X$ on $F_0$ such that 
for each  $\omega =(\omega _1, \omega _2, \cdots )$ and $i \in {\mathbb Z}_+$, 
\[X(\omega )(i)=\omega _N(i), \ \ i\leq 2^N.\]
This completes the proof of \thmu{infinite}.

\vspace{0.3cm}\parr
{\bf Remark}

For $N \in {\mathbb N}$ and $w \in U_N$, let $u_M=(2^{-M}\hat{Q}_Mw)|_{\Delta _0}, 
M=0, 1, \cdots , N$, where $\hat{Q}_M$ is defined at the end of \secu{LE} and $\Delta _0 =\triangle Oa_0b_0$.  
Note that $u_M \in \Gamma _0$.
For any $M\leq N-1$ and any $x_k \in \Gamma _0$, $k=M, M+1, \cdots , N$, 
\eqnb \eqna{Markov2}
P_N^{(i)}[\ u_M =x_M\mid u_k=x_k, \ k=M+1, M+2, \cdots  ,N\ ]=
 P_N^{(i)}[\ u_M =x_M\mid u_{M+1} =x_{M+1}\ ].
\eqne
Thus, 
$P^{(i)}_N, i=1,2,3,4, N \in {\mathbb N}$ define a family of backward Markov chains on the state space $\Gamma _0 = \{v_1^*, v_2^*, v_3^*, v_4^* \}$ such that 
\[ P_N^{(i)}[\ u_N=v^*_i\ ]=1,\]
and for $M \leq N-1$, 
\[ P_N^{(i)}[\ u_M =v_j^* \mid u_{M+1} =v_k^*\ ]=P_{kj},\]
where $P_{kj}$ denotes the $(k, j)$- element of the transition probability matrix 
\[
{\bf P}= \frac{1}{30}\left [
\begin{array}{cccc}
19 & 5 & 5 & 1 \\
5 & 19 & 1 & 5\\
7 & 15 & 5 & 3\\
15 & 7 & 3 & 5
\end{array}
\right ]
.\]
$\displaystyle \alpha = \frac{1}{28}(11,11,3,3)$ is the unique invariant probability vector, 
that is, the unique solution to 
\[\alpha =\alpha P.\]
Moreover, for any probability vector $a$, it holds that 
\[\lim_{n \to \infty} aP^n =\alpha.\]
In terms of the loop-erased walk measures, the above fact can be expressed 
as 
\[  \hat{P}^{(i)}_{N+K}[\ w|_K \in A_K\ ]=\sum _{j=1}^{4} (P^N)_{ij}\hat{P}_K^{(j)}[A_K],\]
where for $w \in \Gamma _{N+K}$,  $w|_K$ denotes the path $w$ stopped at 
$T_1^{ex, K}(w)$ and $A_K \subset \Gamma _K$.
Thus, for any probability vector $a$, we have as $N \to \infty$, 
\[\sum _{i=1}^{4} a_i \hat{P}^{(i)}_{N+K}[\ w|_K \in A_K\ ]
\to \sum_{i=1}^4 \alpha _i 
\hat{P}_K^{(i)}[A_K].\]
In particular, $\displaystyle \frac{1}{6}(2,2,1,1)$ represents the `standard' LERW 
studied in \cite{STW}.


\vspace{0.5cm}\parr
\section{Proof of the theorems}
\seca{Proof}
Let X be the loop-erased random walk defined in  \secu{Extention} and 
let 
\eqnb
\eqna{phitilde}
\tilde{\Phi}_N(x,y)=\frac{11}{14}\Phi ^{(1)}_N(x,y) +\frac{3}{14}\Phi ^{(2)}_N(x,y),
\eqne
where $\Phi ^{(i)}_N(x,y)$, $i=1,2$ are defined in \secu{Steps}. 
The laplace transform of $\lambda ^{-N}T_1^{ex, N}(X)$ is given by 
\eqnb
\eqna{LTtilde}
\tilde{g}_N(t):=\tilde{\Phi}_N(e^{-t\lambda ^{-N} }, e^{-2t\lambda ^{-N} }).
\eqne

Define for each $n \in {\mathbb N}$, 
\[D_n (X)=\min \{M\geq 0\\ :\ |X(i)| \leq 2^M,\ 0\leq i \leq n \},\]
and let $K=K(n)$ be the positive integer such that 
\eqnb \eqna{enke} \lambda ^K\leq n<\lambda ^{K+1}\eqne
holds. 

\prpb (short-path estimate)
\prpa{shortpath}
There exist positive constants $C_{6.1}$ and $C_{6.2}$ such that 
\[P[\ D_n (X)<K(n)-M\ ]\leq C_{6.1} e ^{-C_{6.2}\lambda ^M}\]
holds for any $n, M \in {\mathbb N}$ satisfying  $K(n)>M$.
\prpe

\prfb

Take $C_{6.2}>0$ arbitrarily. Since \prpu{supbr} (3) implies that $\{\tilde{g}_N(t)\}$ is a convergent 
sequence  for any $t \in {\mathbb C}$,  
we can take $C_{6.1}>0$ such that 
$\tilde{g}_N(-C_{6.2}) <C_{6.1}$ for all $N \in {\mathbb N}$.
By Chebyshev's inequality, we have 
\[ \tilde{P}_N[\ \lambda ^{-N} T_1^{ex, N}(X)\geq  \lambda ^M\ ] \leq 
\tilde{g}_N(-C_{6.2}) \ e^{-C_{6.2}\lambda ^M}<C_{6.1} \ e^{-C_{6.2}\lambda ^M}.\]
This leads to 
\begin{eqnarray*}
P[\ D_n (X)<K(n)-M\ ] &\leq& P [\ T_1^{ex, K-M}(X) >n\ ] \\ 
&=&\ 
\tilde{P}_{K-M}[\  T_1^{ex, K-M}(w) >n\ ] \\ 
&\leq&
\tilde{P}_{K-M}[\ \lambda ^{-(K-M)}T_1^{ex, K-M} (w) >\lambda ^M\ ] \\ 
&\leq &
C_{6.1} e^{-C_{6.2} \lambda ^M}.
\end{eqnarray*}
\QED\prfe


\prpb (long-path estimate)
\prpa{longpath}
There exist $C_{6.3}, C_{6.4}>0$ and $N_0\in {\mathbb N}$ such that 
\[P[\ D_n (X)> K(n)+M\ ] \leq C_{6.3} e^{-C_{6.4} 2^M}\]
for any $n$ satisfying $K(n)\geq N_0$ and any  $M\in {\mathbb N}$.
\prpe

\prfb
First note that
\begin{eqnarray*}
P[\ D_n (X)> K(n)+M\ ] &\leq& P [\ T_1^{ex, K+M}(X) <n\ ] \\ 
&=&\ 
\tilde{P}_{K+M}[\  T_1^{ex, K+M}(w) <n\ ] \\ 
&\leq&
\tilde{P}_{K+M}[\ T_1^{ex, K+M} (w) \leq \lambda ^{K+1}\ ] . 
\end{eqnarray*}
Fix $0< \delta <1$ arbitrarily, then 
\begin{eqnarray*}
\tilde{P}_{K+M}[\ T_1^{ex, K+M}(w) \leq \lambda ^{K+1}\ ] &=& \sum _{w \in \Gamma _{K+M},\ \ell (w) \leq \lambda ^{K+1}}\tilde{P}_{K+M}[w] \\ 
&\leq&
\delta ^{-1} \sum _{w \in \Gamma _{K+M}, \ \ell (w) \leq \lambda ^{K+1}}\tilde{P}_{K+M}[w]  \ \delta ^{\ell (w) \lambda ^{-(K+1)}} \\ 
&\leq& \delta ^{-1} \tilde{\Phi} _{K+M} (\delta ^{\lambda ^{-(K+1)}}  ,\ \delta ^{2\lambda ^{-(K+1)} }).
\end{eqnarray*}
Let $t'=-\lambda ^{-1} \log \delta >0$.  Since \prpu{supbr} (3) implies that $
\displaystyle \tilde{\Phi} _{N} (\delta ^{\lambda ^{-(N+1)}},  \ \delta ^{2\lambda ^{-(N+1)} })=\tilde{g}_N(t')  $ converges 
as $N \to \infty$ to a limit strictly smaller than $1$.  
we can choose $0<r<1$ and $N_0 \in {\mathbb N}$ such that 
\eqnb \eqna{r} 
\Phi ^{(i)}_{N}(\delta ^{\lambda ^{-(N+1)}}, \ \delta ^{2\lambda ^{-(N+1)} })<r ,\ \ i=1,2
\eqne
for all $N \geq N_0$. 
Thus if $K \geq N_0$, 
\[ \tilde{\Phi} _{K+M} (\delta ^{\lambda ^{-(K+1)}}  , \ \delta ^{2\lambda ^{-(K+1)} })
<\tilde{\Phi} _{M}(r,r) \leq   r ^{2^M} =  e ^{-C_{6.4}2^M},\]
where we used \eqnu{2M} in the last inequality and set $C_{6.4}=-\log r$.
Taking $C_{6.3}=\delta ^{-1}$ completes the proof.

\QED
\prfe

To obtain the displacement exponent, 
we shall use the following inequality that holds for any ${\mathbb N}$--valued 
random variable $Y$ and $s>0$:
\eqnb \eqna{Y} s \ C_{6.5}(s) \sum_{k=1}^{\infty } k^{s-1}  P [\ Y\geq k\ ] \leq  E[Y^s]
\leq s \sum _{k=1}^{\infty } k^{s-1}  P [\ Y\geq k\ ]+C_{6.6}(s).\eqne

For $0<s< 1$, $C_{6.5}(s)=1$, $C_{6.6}(s)=1$, for $s>1$, $C_{6.5}(s)=\frac{1}{2^s}$, $C_{6.6}(s)=0$ and $C_{6.5}(1)=1$, $C_{6.6}(1)=0$.

Let $\nu=\log 2/\log \lambda$.
\prpb
\prpa{lowerbound}
For any $s>0$, there exist a positive constant $C_1(s)$ and $n_1 \in {\mathbb N}$ 
such that 
\[E[\ |X(n)|^s\ ] \geq C_1(s)\ n^{s \nu}, \]
for all $n>n_1$.
\prpe

\prfb
Fix $M_0 \in {\mathbb N}$ such that $C_{6.1}e^{-C_{6.2}\lambda ^{M_0}} <1/2$, where 
$C_{6.1}$ and $C_{6.2}$ are as in \prpu{shortpath}.  Take an  $n$ large enough so that 
$K(n)>M_0+2$, where $K(n)$ is as in \eqnu{enke}. 
Then  
\eqnb \eqna{1/2} P[\  |X(n)|\leq 2^{K-M_0-2}\ ]\leq P [\ D_n<K-M_0\ ]<\frac{1}{2}.\eqne
We give a proof in the case for $s>1$. 
We make use of \eqnu{Y} with  $P[\  |X(n)| >n \ ]=0$ in mind.
\begin{eqnarray*}
E[\ |X(n)|^s\ ] 
&\geq &\ 
\frac{s}{2^s} \sum _{m=0}^{\infty  } \  \sum _{k=2^m+1}^{2^{m+1} }  k^{s-1} P [\ |X(n)|\geq k\ ]\ \\ 
&\geq&
\frac{s}{2^s} \sum _{m=0}^{\infty  } \ \sum _{k=2^m+1}^{2^{m+1} } (2^m)^{s-1}\ P [\ |X(n)|>2^{m+1}\ ]\  \\ 
&\geq & 
s  \sum _{m=0}^{\infty  } \   2^{s(m-1)}\ P [\ |X(n)|>2^{m+1}\ ]\\
&\geq & 
s \ 2^{-(M_0+4)s}\  2^{Ks} \ P [\ |X(n)|>2^{K-M_0-2}\ ]\\
\\
&=&
s  \ 2^{-(M_0+4)s}\  2^{Ks}(1-P [\ |X(n)|\leq 2^{K-M_0-2}\ ])\\
\\
&\geq& s  \ 2^{-(M_0+5)s-1}\  2^{Ks}\geq C_1(s) n^{s \nu },
\end{eqnarray*}
where we used \eqnu{1/2} and set $C_1(s)=s 2^{-(M_0+5)s-1}$.
The case for $0<s\leq 1$ can be proved similarly.

\QED
\prfe

\prpb
\prpa{upperbound}
For any $s>0$, there exist a positive constant $C_2(s)$ and $n_2 \in {\mathbb N}$ 
such that 
\[E[\ |X(n)|^s\ ] \leq C_2(s)\ n^{s \nu} \]
for all $n>n_2$.
\prpe

\prfb
First note that  
\eqnb \eqna{XD} P[\ |X(n)| \geq 2^{m}\ ] \leq P[\ D_n(X) > m-1\ ].\eqne
Assume $K=K(n) \geq N_0$ as in \prpu{longpath}.  In the case of $s>1$, making use of \eqnu{Y}, we have
\begin{eqnarray*}
E[\ |X(n)|^s\ ] &\leq &
s \sum _{m=0}^{\infty  } \ 2^m \cdot  2^{(s-1)(m+1)}\ P [\ |X(n)|\geq 2^m\ ]\ \\ 
&\leq &
s 2^{s-1} \biggl(\sum _{m=0}^{K+1 } \  2^{sm}\ P [\ |X(n)|\geq 2^m\ ]\  +\sum _{m=K+2}^{\infty  } \  2^{sm}\ P [\ |X(n)|\geq 2^m\ ] \ \biggr) \\
&\leq &
s 2^{s-1} \biggl(\sum _{m=0}^{K+1 } \  2^{sm}\   + \sum _{m=K+2}^{\infty  } \  2^{sm}\ P [\ D_n(X) > m-1\ ] \ \biggr) \mbox{      \hspace{1cm} (use of \eqnu{XD})}\\
&\leq & 
 c_1 (s) 2^{Ks} + s 2^{2s-1} C_{6.3}\  2^{Ks}\ \sum _{\ell =1}^{ \infty} 2 ^{\ell s}  e^{-C_{6.4}2^{\ell } }\mbox{      \hspace{1cm} (\prpu{longpath})} \\
& \leq& C_2(s) n^{s \nu },
\end{eqnarray*}
where $c_1(s)$ and $C_2(s)$ are positive constants depending only on $s$ and we used 
the convergence of the series above. 
The case for $0<s\leq 1$ can be proved similarly. 
\QED
\prfe

\prpu{upperbound} combined with \prpu{lowerbound} gives \thmu{exponent}.

Now we go on to prove the law of the iterated logarithm. 
First we prove the upper bound: 
\prpb
\prpa{LIL1}
There exists  $C_4>0$ such that 
\[\limsup _{n \to \infty } \frac{|X(n)|}{\psi (n) }\leq C_4, \ \ \ P\mbox{-- a.s.}, \]
where $\psi (n)=n^{\nu}(\log \log n)^{1-\nu}$.
\prpe


\vspace{0.5cm}\parr
\prfb
Let $\mu _N$ be the distribution of $\lambda^{-N} T_1^{ex, N}(X)$ under $P$.
For each $x>1$ there is a unique integer $N$ such that 
$2^N\leq x<2^{N+1}$.  For $k>0$ satisfying  $ 2^{-N}k \geq C_{4.7}$,  \thmu{Tauberian} (2) implies that 
\begin{eqnarray*}
P[\ \max _{0\leq j \leq k} |X(j)| >x \ ] 
&\leq& P[\ T_1^{ex, N}(X) \leq k\ ] \\
&=&\mu _N([0, \lambda ^{-N} k]])\\
&\leq& \displaystyle C_{4.8} e^{-C_{4.9}( xk^{-\nu}/2) ^{1/(1-\nu )}}.
\end{eqnarray*}
Let $\gamma >1$ be arbitrary.  For $A>0$, let $x=A\psi (\gamma ^m)$ and $k$ be the largest integer that does not exceed $\gamma ^{m+1}$.  
The condition $2^{-N}k\geq C_{4.7}$ is satisfied for $m$ large enough.  
Thus, the above inequality leads to  
\begin{eqnarray*}
\sum _{m=1}^{\infty }P[\ \max _{\gamma ^m<j\leq \gamma ^{m+1} } |X(j)| >A\psi (\gamma ^m)\ ]  &\leq & \sum _{m=1}^{\infty }P[\ \max _{0 \leq j\leq \gamma ^{m+1} } |X(j)| >A\psi (\gamma ^m) \ ]\\
&\leq & c_{3}+C_{4.8}\sum _{m=1}^{\infty }e^{-C_{4.9}(x k^{-\nu}/2) ^{1/(1-\nu )}}\\
&\leq & \displaystyle c_{3}+c_4 C_{4.8} \sum _{m=1}^{\infty }\frac{1}{m^{\alpha} },
\end{eqnarray*}
for some constants $c_{3}, c_4>0$ 
and  $\displaystyle \alpha =C_{4.9}\biggl(\frac{A}{2  \gamma ^{\nu}}\biggr)^{1/(1-\nu)}$.
The  sequence $\dsp \sum _{m=1}^{\infty }\frac{1}{m^{\alpha} }$ converges if we take $A$ large enough so that $\alpha >1$. 
The rest is a usual Borel-Cantelli argument and the statement holds with $C_4=A$.
\QED
\prfe

Now we show the lower bound:  
\prpb
\prpa{LIL2}
There exists $C_3>0$ such that 
\[C_3 \leq \limsup _{n \to \infty } \frac{|X(n)|}{\psi (n) }, \ \ \ P\mbox{-  a.s.} \]
holds.
\prpe

The proof goes along the line of the argument used in \cite{HamblyKumagai}, but we need to show how the ELLF constraction enables us to make use of a `Markov structure' to obtain the result.  We use the following lemma:

\lemb (A version of the second Borel-Cantelli Lemma used in \cite{HHH}) 
\lema{BC}
Let $B_1, B_2,  \cdots $ be a sequence of events and assume  
\[P[B_m\mid B_{m+1}^c, B_{m+2}^c, \cdots , B_{m+k}^c]=P[B_m\mid B_{m+1}^c],\]
for all $m, k \in {\mathbb N}$.
Then 
\[\sum _{m=1}^{\infty }P[B_m\mid B_{m+1}^c] =\infty\]
implies
\[P[\limsup _{m \to \infty } B_m]=1.\]

\leme

\parr
{\it Proof of \prpu{LIL2}.}

Let $\beta =(1-\nu)/\nu$,  $0<b<1$ and 
\[A_M:=\{\lambda ^{-M}T_1^{ex, M}(X)\leq (b\log M)^{-\beta } \}.\]
We want to show that for an appropriate choice of $b$, 
\[P[\limsup_{M \to \infty } A_M]=1\]
holds. 
Let $S^{M-1}_j(X)= T_j ^{ex, M-1}(X) - T_{j-1} ^{ex, M-1}(X)$, $j \in \{1,2, 3\}$, then 
\eqnb \eqna{TS}
T_1^{ex, M}(X)=\sum _{j=1}^{|\sigma _{M-1}|}S_j^{M-1}(X),
\eqne
where $|\sigma _{M-1}|$ denotes the number of $2^{M-1}$--triangles in the 
$2^{M-1}$--skeleton of $X$ stopped at $T_{1}^{ex, M}(X)$, 
which is either $2$ or $3$.  
Since the right-hand side contains only $S_j^{M-1}$, $j \in \{1,2, 3\}$, we have
\eqnb \eqna{aac} P[A_M\mid A_{M-1}^c A_{M-2}^c, \cdots , A_{M-k}^c]=P[A_M\mid A_{M-1}^c],\eqne 
and by repeated use of the definition of conditional probability combined with \eqnu{aac}, we have 
 \[P[A_M\mid A_{M+1}^c, A_{M+2}^c, \cdots , A_{M+k}^c]=P[A_M\mid A_{M+1}^c].\]
In order to show $\displaystyle \sum_{M=1}^{\infty } P[A_M\mid A_{M+1}^c]=\infty$, 
it is sufficient to show $\displaystyle \sum_{M=1}^{\infty } P[A_M\cap A_{M+1}^c]=\infty$.
Let $x_M=(b\log M)^{-\beta }$, then since $S_1^M (X)=T_1^{ex, M} (X)$,
\[ \begin{array}{l} \dsp
P[A_M \cap A_{M+1}^c] =
\tilde{P}_{M+1}[\ T_1^{ex, M}(w)\leq \lambda ^Mx_M,\ 
T_1^{ex, M+1}(w)>  \lambda ^{M+1}x_{M+1}\  ]\\
\\
\dsp\phantom{\hat{P}_{N+1}^{(i)}}
=\sum_{y \leq \lambda ^M x_M} 
\tilde{P}_{M+1} [\ T_1^{ex, M}(w)=y, \ \sum_{i\geq 2} S_i^M(w)+y> \lambda ^{M+1}x_{M+1}\  ]\\
\\
\dsp\phantom{\hat{P}_{N+1}^{(i)}}
\geq \sum_{y \leq \lambda ^M x_M} 
\tilde{P}_{M+1} [\ T_1^{ex, M} (w)=y, \  S_2^M+y> \lambda ^{M+1}x_{M+1}\  ].
\end{array}\]

Let $\tilde{Z}_N$ be the simple random walk $\tilde{P}_N$ defines on $F_0\cap \triangle Oa_Nb_N$. Recall the procedure for loop erasure: after erasing largest-scale loops from $\tilde{Z}_{M+1}$ in Step (2),  we get $\hat{Q}_M\tilde{Z}_{M+1}$ and the law of  $2^{-M}\hat{Q}_M\tilde{Z}_{M+1}$ is equal to $\tilde{P}_1$.  We restore the original fine structures 
to these remaining parts and continue loop erasure. 
For each $\Delta _i$ in $\sigma_M(\hat{Q}_M\tilde{Z}_{M+1})$, if 
 $\Delta _i$ is Type 1 with regard to $\hat{Q}_M\tilde{Z}_{M+1}$, the rest of the procedure is the same as loop erasure for $Z_M$ (modulo rotation and reflection), and if Type 2, the same as that 
for $Z'_M$($Z_M$ and $Z'_M$ are defined in \eqnu{Z} and \eqnu{Zprime}). 
Conditioned on $\hat{Q}_M\tilde{Z}_{M+1}$, parts in different $2^M$--triangles 
are independent.   
Classifying by the types of $\Delta_1$ and $\Delta _2$ in $\sigma _M (\hat{Q}_M \tilde{Z}_{M+1})$, we have
\[ \begin{array}{l} \dsp
 \tilde{P}_{M+1} [\ T_1^{ex, M}(w)=y, \  S_2^M(w)+y> \lambda ^{M+1}x_{M+1}\  ]\\
\dsp\phantom{\hat{P}_{N+1}^{(i)}}
=\hat{P}_{M} [\ T_1^{ex, M} (w)=y\ ] \ \hat{P}_M[\ T_1^{ex, M}(w)+y> \lambda ^{M+1}x_{M+1}\  ]
\ \tilde{P}_1[\{w^*_1, w^*_5, w^*_7\}]\\
\dsp\phantom{\hat{P}_{N+1}^{(i)}=}
+\hat{P}'_{M} [\ T_1^{ex, M}(w)=y\ ] \ \hat{P}_M[\ T_1^{ex, M}(w)+y> \lambda ^{M+1}x_{M+1}\  ]
\ \tilde{P}_1[\{w^*_2, w^*_6, w^*_{10}\}]\\
\dsp\phantom{\hat{P}_{N+1}^{(i)}=}
+\hat{P}_{M} [\ T_1^{ex, M}(w)=y\ ] \ \hat{P}'_M[\ T_1^{ex, M}(w)+y> \lambda ^{M+1}x_{M+1}\  ]
\ \tilde{P}_1[\{w^*_3, w^*_8, w^*_{9}\}]\\
\dsp\phantom{\hat{P}_{N+1}^{(i)}=}
+\hat{P}'_{M} [\ T_1^{ex, M}(w)=y\ ] \ \hat{P}'_M[\ T_1^{ex, M}(w)+y> \lambda ^{M+1}x_{M+1}\  ]
\ \tilde{P}_1[\{w^*_4 \}].
\end{array}\]

Since $x_{M+1}/x_M \to 1$ as $M \to \infty$, we can take 
$0<c<1$ such that $c\lambda ^{M+1}x_{M+1} <\lambda ^Mx_M$ for all large enough $M$.
Then we have  
\[ \begin{array}{l} \dsp
\sum _{y\leq \lambda ^Mx_M}
\hat{P}_{M} [\ T_1^{ex, M} (w)=y\ ] \ \hat{P}_M[\ T_1^{ex, M}(w)+y> \lambda ^{M+1}x_{M+1}\  ]\\
\dsp\phantom{\hat{P}_{M+1}^{(i)}}
\geq 
\sum _{c\lambda ^{M+1} x_{M+1}\leq  y\leq \lambda ^Mx_M}
\hat{P}_{M} [\ T_1^{ex, M} (w)=y\ ] \ \hat{P}_M[\ T_1^{ex, M}(w)+y> \lambda ^{M+1}x_{M+1}\  ]\\
\\
\dsp\phantom{\hat{P}_{N+1}^{(i)}}
\geq
\hat{P}_{M} [\ T_1^{ex, M} (w)
\in [\  c\lambda ^{M+1}x_{M+1}, \lambda ^Mx_M] \ ] \ \hat{P}_M[\  T_1^{ex, M} (w)> \lambda ^{M+1} (1-c)x_{M+1}\  ].
\end{array}\]

By \prpu{supbr} (2),  $\lambda ^{-M}T_1^{ex, M}$ under $\hat{P}_M$ and under  $\hat{P}'_M$ converge in law to $T_1^*$ and $T_2^{*}$, respectively, as $M \to \infty$, 
which combined with the fact that $x_M \to 0$ as $M\to \infty $ leads to 
\[
\hat{P}_M[\  T_1^{ex, M} (w)> \lambda ^{M+1} (1-c)x_{M+1}\  ]
=\hat{P}_M[\  \lambda ^{-M} T_1^{ex, M} (w)> \lambda  (1-c)x_{M+1}\  ]\]
\[
\geq \hat{P}_M[\  \lambda ^{-M} T_1^{ex, M} (w)>1] >\frac{1}{2} P[T_1^*>1], \]
for $M$ large enough.

With similar argument for the other terms, we have 

\[ \begin{array}{l} \dsp
P[A_M \cap A_{M+1}^c] \\
\dsp\phantom{\hat{P}_{M+1}^{(i)}}
> a \hat{P}_{M} [\ \ T_1^{ex, M} (w)
\in [\  c\lambda ^{M+1}x_{M+1}, \lambda ^Mx_M\ ]  \ ] \ 
\tilde{P}_1[\{w^*_1, w^*_3, w_5^*, w^*_7, w_8^*, w_9^* \}]\\
\dsp\phantom{\hat{P}_{N+1}^{(i)}=}
+  a \hat{P}'_{M} [\ \ T_1^{ex, M} (w)
\in [\  c\lambda ^{M+1}x_{M+1}, \lambda ^Mx_M\ ]  \ ] \ 
\tilde{P}_1[\{w^*_2, w^*_4, w_6^*, w^*_{10} \}]\\
\dsp\phantom{\hat{P}_{M+1}^{(i)}}
=a\tilde{P}_{M} [\ \ T_1^{ex, M} (w)
\in [\  c\lambda ^{M+1}x_{M+1}, \lambda ^Mx_M\ ] \ ] ,
\end{array}\]
where $\displaystyle a =\frac{1}{2}(P[T_1^*>1]\wedge P[T_2^*>1])>0$.
Moreover, 
\[ \begin{array}{l} \dsp
\tilde{P}_{M} [\ \ T_1^{ex, M} (w)\in [\  c\lambda ^{M+1}x_{M+1}, \lambda ^Mx_M\ ] \ ] \\
\dsp\phantom{\hat{P}_{M+1}^{(i)}}
=\tilde{P}_{M} [\ \ \lambda ^{-M}T_1^{ex, M} (w)  \in [ 0, x_M ] \ ] 
\left(1-\frac{\tilde{P}_{M} [\ \ \lambda ^{-M}T_1^{ex, M} (w) \in [0, \lambda c x_{M+1}]}{\tilde{P}_{M} [\ \ \lambda ^{-M}T_1^{ex, M} (w) \in [0, x_{M}]\ ]}\right).
\end{array}\]

Since $x_M\to 0$ and $2^{M(1-\nu)/\nu}x_M \to \infty$ as $M \to \infty$,  
\thmu{Tauberian}(1) with $\alpha _N=x_N$ implies that for large enough $M$ 
\[
\tilde{P}_{M} [\ \ \lambda ^{-M}T_1^{ex, M} (w)  \in [ 0, x_M ] \ ] 
\geq \exp (-C_{4.5} x_M^{-1/\beta}) =\exp (-C_{4.5} b\log M)=\frac{1}{M^{C_{4.5} b}}.
\]
We use \thmu{Tauberian} (1) again to have  
\[
1-\frac{\tilde{P}_{M} [\ \ \lambda ^{-M}T_1^{ex, M} (w) \in [0, \lambda c x_{M+1}]}{\tilde{P}_{M} [\ \ \lambda ^{-M}T_1^{ex, M} (w) \in [0, x_{M}]]}
\to 1\]
as $M \to \infty $, thus 
this factor is greater than $1/2$ for large enough $M$.
Choose $0<b<1$ so that $C_{4.5} b <1$, then for some $M_0>0$ and $c_5>0$, it holds that 
\[\sum _{M=1}^{\infty } P[A_M \cap A_{M+1}^c] \geq c_5+ \frac{1}{2}\sum_{M\geq M_0}  
\frac{1}{M^{C_{4.5}b}}  = +\infty.\]
\lemu{BC} implies that for almost all $\omega \in \Omega$, there exists an increasing  sequence  $\{M_k (\omega )\}$, $k=1,2, \ldots $ such that 
\eqnb \eqna{BCMk} \lambda ^{-M_k}T_1^{ex, M_k}(X)\leq (b \log M_k)^{-\beta }.\eqne
It follows that for $M_k\geq 3$ 
\[M_k\geq \frac{\log T_1^{ex, M_k}(X)+\beta \log b }{\log \lambda } +
\frac{\beta \log \log M_k}{ \log \lambda} 
\geq 
\frac{\log T_1^{ex, M_k}(X)+\beta \log b }{\log \lambda },\]
and for any small $\varepsilon >0$, there exists a $k_0 \in {\mathbb N}$ such that 
\eqnb \eqna{logMk} \log M_k \geq (1-\varepsilon ) \log \log T_1^{ex, M_k}(X) \eqne
holds for all $k\geq k_0$.  On the other hand, \eqnu{BCMk} implies
\[|X(T_1^{ex, M_k} (X))| =2^{M_k}\geq (b \log M_k)^{1-\nu } \ M_k^{\nu} .\]  
This combined with \eqnu{logMk} leads to 
\[\limsup _{n\to \infty}\frac{|X(n)|}{\psi (n)}
 \geq b^{1-\nu}(1-\varepsilon)^{1-\nu}. \]
Since $\varepsilon $ is arbitrary, we have proved the proposition with $C_3=b^{1-\nu}$.

\QED

\prpu{LIL1} combined with \prpu{LIL2} gives \thmu{loglog}.


\vspace{0.5cm}\parr
{\Large\bf Acknowledgments}
\vspace{0.2cm}\parr
This work is supported by JSPS KAKENHI Grant Number 16K05210.



\end{document}

%% file: Fig1.tex
\unitlength 0.1in
\begin{picture}( 42.8000, 20.2400)(  1.9500,-22.7200)
%
{\color[named]{Black}{%
\special{pn 8}%
\special{pa 2168 2042}%
\special{pa 1218 2042}%
\special{pa 1694 1220}%
\special{pa 2168 2042}%
\special{pa 1218 2042}%
\special{fp}%
}}%
%
{\color[named]{Black}{%
\special{pn 8}%
\special{pa 242 2038}%
\special{pa 480 1628}%
\special{pa 718 2038}%
\special{pa 242 2038}%
\special{pa 480 1628}%
\special{fp}%
}}%
%
{\color[named]{Black}{%
\special{pn 8}%
\special{pa 1694 2042}%
\special{pa 1932 1630}%
\special{pa 2170 2042}%
\special{pa 1694 2042}%
\special{pa 1932 1630}%
\special{fp}%
}}%
%
{\color[named]{Black}{%
\special{pn 8}%
\special{pa 1694 2042}%
\special{pa 1932 1630}%
\special{pa 2170 2042}%
\special{pa 1694 2042}%
\special{pa 1932 1630}%
\special{fp}%
}}%
%
{\color[named]{Black}{%
\special{pn 8}%
\special{pa 1694 2042}%
\special{pa 1932 1630}%
\special{pa 2170 2042}%
\special{pa 1694 2042}%
\special{pa 1932 1630}%
\special{fp}%
}}%
%
{\color[named]{Black}{%
\special{pn 8}%
\special{pa 1694 2042}%
\special{pa 1932 1630}%
\special{pa 2170 2042}%
\special{pa 1694 2042}%
\special{pa 1932 1630}%
\special{fp}%
}}%
%
{\color[named]{Black}{%
\special{pn 8}%
\special{pa 1218 2042}%
\special{pa 1456 1630}%
\special{pa 1694 2042}%
\special{pa 1218 2042}%
\special{pa 1456 1630}%
\special{fp}%
}}%
%
{\color[named]{Black}{%
\special{pn 8}%
\special{pa 1218 2042}%
\special{pa 1456 1630}%
\special{pa 1694 2042}%
\special{pa 1218 2042}%
\special{pa 1456 1630}%
\special{fp}%
}}%
%
{\color[named]{Black}{%
\special{pn 8}%
\special{pa 1218 2042}%
\special{pa 1456 1630}%
\special{pa 1694 2042}%
\special{pa 1218 2042}%
\special{pa 1456 1630}%
\special{fp}%
}}%
%
{\color[named]{Black}{%
\special{pn 8}%
\special{pa 1218 2042}%
\special{pa 1456 1630}%
\special{pa 1694 2042}%
\special{pa 1218 2042}%
\special{pa 1456 1630}%
\special{fp}%
}}%
%
{\color[named]{Black}{%
\special{pn 8}%
\special{pa 1456 1630}%
\special{pa 1694 1218}%
\special{pa 1932 1630}%
\special{pa 1456 1630}%
\special{pa 1694 1218}%
\special{fp}%
}}%
\put(21.5600,-20.8800){\makebox(0,0)[lt]{$b_1$}}%
\put(14.5200,-16.2900){\makebox(0,0)[rb]{$a_0$}}%
\put(16.5000,-12.1500){\makebox(0,0)[lb]{$a_1$}}%
\put(16.9300,-21.4500){\makebox(0,0){$b_0$}}%
\put(11.7100,-20.8800){\makebox(0,0)[lt]{$O$}}%
\put(1.9500,-20.8500){\makebox(0,0)[lt]{$O$}}%
\put(7.4300,-21.4700){\makebox(0,0){$b_0$}}%
\put(5.5700,-16.1400){\makebox(0,0)[rb]{$a_0$}}%
\put(16.9300,-23.3700){\makebox(0,0){$F'_1$}}%
%
{\color[named]{Black}{%
\special{pn 8}%
\special{pa 3526 2042}%
\special{pa 2576 2042}%
\special{pa 3050 1220}%
\special{pa 3526 2042}%
\special{pa 2576 2042}%
\special{fp}%
}}%
%
{\color[named]{Black}{%
\special{pn 8}%
\special{pa 3050 2042}%
\special{pa 3288 1630}%
\special{pa 3526 2042}%
\special{pa 3050 2042}%
\special{pa 3288 1630}%
\special{fp}%
}}%
%
{\color[named]{Black}{%
\special{pn 8}%
\special{pa 3050 2042}%
\special{pa 3288 1630}%
\special{pa 3526 2042}%
\special{pa 3050 2042}%
\special{pa 3288 1630}%
\special{fp}%
}}%
%
{\color[named]{Black}{%
\special{pn 8}%
\special{pa 3050 2042}%
\special{pa 3288 1630}%
\special{pa 3526 2042}%
\special{pa 3050 2042}%
\special{pa 3288 1630}%
\special{fp}%
}}%
%
{\color[named]{Black}{%
\special{pn 8}%
\special{pa 3050 2042}%
\special{pa 3288 1630}%
\special{pa 3526 2042}%
\special{pa 3050 2042}%
\special{pa 3288 1630}%
\special{fp}%
}}%
%
{\color[named]{Black}{%
\special{pn 8}%
\special{pa 2576 2042}%
\special{pa 2814 1630}%
\special{pa 3052 2042}%
\special{pa 2576 2042}%
\special{pa 2814 1630}%
\special{fp}%
}}%
%
{\color[named]{Black}{%
\special{pn 8}%
\special{pa 2576 2042}%
\special{pa 2814 1630}%
\special{pa 3052 2042}%
\special{pa 2576 2042}%
\special{pa 2814 1630}%
\special{fp}%
}}%
%
{\color[named]{Black}{%
\special{pn 8}%
\special{pa 2576 2042}%
\special{pa 2814 1630}%
\special{pa 3052 2042}%
\special{pa 2576 2042}%
\special{pa 2814 1630}%
\special{fp}%
}}%
%
{\color[named]{Black}{%
\special{pn 8}%
\special{pa 2576 2042}%
\special{pa 2814 1630}%
\special{pa 3052 2042}%
\special{pa 2576 2042}%
\special{pa 2814 1630}%
\special{fp}%
}}%
%
{\color[named]{Black}{%
\special{pn 8}%
\special{pa 2814 1630}%
\special{pa 3050 1218}%
\special{pa 3288 1630}%
\special{pa 2814 1630}%
\special{pa 3050 1218}%
\special{fp}%
}}%
\put(34.9400,-20.8800){\makebox(0,0)[lt]{$b_1$}}%
\put(28.0900,-16.2900){\makebox(0,0)[rb]{$a_0$}}%
\put(29.2600,-12.1000){\makebox(0,0)[lb]{$a_1$}}%
\put(30.6500,-21.4500){\makebox(0,0){$b_0$}}%
\put(25.2800,-20.8800){\makebox(0,0)[lt]{$O$}}%
\put(35.2800,-23.3700){\makebox(0,0){$F'_2$}}%
%
{\color[named]{Black}{%
\special{pn 8}%
\special{pa 4476 2040}%
\special{pa 3526 2040}%
\special{pa 4000 1216}%
\special{pa 4476 2040}%
\special{pa 3526 2040}%
\special{fp}%
}}%
%
{\color[named]{Black}{%
\special{pn 8}%
\special{pa 4000 2040}%
\special{pa 4238 1628}%
\special{pa 4476 2040}%
\special{pa 4000 2040}%
\special{pa 4238 1628}%
\special{fp}%
}}%
%
{\color[named]{Black}{%
\special{pn 8}%
\special{pa 4000 2040}%
\special{pa 4238 1628}%
\special{pa 4476 2040}%
\special{pa 4000 2040}%
\special{pa 4238 1628}%
\special{fp}%
}}%
%
{\color[named]{Black}{%
\special{pn 8}%
\special{pa 4000 2040}%
\special{pa 4238 1628}%
\special{pa 4476 2040}%
\special{pa 4000 2040}%
\special{pa 4238 1628}%
\special{fp}%
}}%
%
{\color[named]{Black}{%
\special{pn 8}%
\special{pa 4000 2040}%
\special{pa 4238 1628}%
\special{pa 4476 2040}%
\special{pa 4000 2040}%
\special{pa 4238 1628}%
\special{fp}%
}}%
%
{\color[named]{Black}{%
\special{pn 8}%
\special{pa 3526 2040}%
\special{pa 3762 1628}%
\special{pa 4000 2040}%
\special{pa 3526 2040}%
\special{pa 3762 1628}%
\special{fp}%
}}%
%
{\color[named]{Black}{%
\special{pn 8}%
\special{pa 3526 2040}%
\special{pa 3762 1628}%
\special{pa 4000 2040}%
\special{pa 3526 2040}%
\special{pa 3762 1628}%
\special{fp}%
}}%
%
{\color[named]{Black}{%
\special{pn 8}%
\special{pa 3526 2040}%
\special{pa 3762 1628}%
\special{pa 4000 2040}%
\special{pa 3526 2040}%
\special{pa 3762 1628}%
\special{fp}%
}}%
%
{\color[named]{Black}{%
\special{pn 8}%
\special{pa 3526 2040}%
\special{pa 3762 1628}%
\special{pa 4000 2040}%
\special{pa 3526 2040}%
\special{pa 3762 1628}%
\special{fp}%
}}%
%
{\color[named]{Black}{%
\special{pn 8}%
\special{pa 3762 1628}%
\special{pa 4000 1216}%
\special{pa 4238 1628}%
\special{pa 3762 1628}%
\special{pa 4000 1216}%
\special{fp}%
}}%
\put(44.6300,-20.8700){\makebox(0,0)[lt]{$b_2$}}%
%
{\color[named]{Black}{%
\special{pn 8}%
\special{pa 4000 1216}%
\special{pa 3050 1216}%
\special{pa 3526 394}%
\special{pa 4000 1216}%
\special{pa 3050 1216}%
\special{fp}%
}}%
%
{\color[named]{Black}{%
\special{pn 8}%
\special{pa 3050 1216}%
\special{pa 3288 804}%
\special{pa 3526 1216}%
\special{pa 3050 1216}%
\special{pa 3288 804}%
\special{fp}%
}}%
%
{\color[named]{Black}{%
\special{pn 8}%
\special{pa 3524 1216}%
\special{pa 3762 804}%
\special{pa 4000 1216}%
\special{pa 3524 1216}%
\special{pa 3762 804}%
\special{fp}%
}}%
%
{\color[named]{Black}{%
\special{pn 8}%
\special{pa 3286 804}%
\special{pa 3526 390}%
\special{pa 3762 804}%
\special{pa 3286 804}%
\special{pa 3526 390}%
\special{fp}%
}}%
\put(34.4300,-3.7800){\makebox(0,0)[lb]{$a_2$}}%
\put(4.7600,-23.3600){\makebox(0,0){$F'_0$}}%
\end{picture}%

%% file: Fig2.tex
\unitlength 0.1in
\begin{picture}( 43.9500, 22.9500)(  0.7500,-26.4000)
%
{\color[named]{Black}{%
\special{pn 8}%
\special{pa 2760 2376}%
\special{pa 2268 2376}%
\special{pa 2514 1950}%
\special{pa 2760 2376}%
\special{pa 2268 2376}%
\special{fp}%
}}%
%
{\color[named]{Black}{%
\special{pn 8}%
\special{pa 2514 2376}%
\special{pa 2638 2160}%
\special{pa 2760 2376}%
\special{pa 2514 2376}%
\special{pa 2638 2160}%
\special{fp}%
}}%
%
{\color[named]{Black}{%
\special{pn 8}%
\special{pa 2514 2376}%
\special{pa 2638 2160}%
\special{pa 2760 2376}%
\special{pa 2514 2376}%
\special{pa 2638 2160}%
\special{fp}%
}}%
%
{\color[named]{Black}{%
\special{pn 8}%
\special{pa 2514 2376}%
\special{pa 2638 2160}%
\special{pa 2760 2376}%
\special{pa 2514 2376}%
\special{pa 2638 2160}%
\special{fp}%
}}%
%
{\color[named]{Black}{%
\special{pn 8}%
\special{pa 2514 2376}%
\special{pa 2638 2160}%
\special{pa 2760 2376}%
\special{pa 2514 2376}%
\special{pa 2638 2160}%
\special{fp}%
}}%
%
{\color[named]{Black}{%
\special{pn 8}%
\special{pa 2268 2376}%
\special{pa 2392 2160}%
\special{pa 2514 2376}%
\special{pa 2268 2376}%
\special{pa 2392 2160}%
\special{fp}%
}}%
%
{\color[named]{Black}{%
\special{pn 8}%
\special{pa 2268 2376}%
\special{pa 2392 2160}%
\special{pa 2514 2376}%
\special{pa 2268 2376}%
\special{pa 2392 2160}%
\special{fp}%
}}%
%
{\color[named]{Black}{%
\special{pn 8}%
\special{pa 2268 2376}%
\special{pa 2392 2160}%
\special{pa 2514 2376}%
\special{pa 2268 2376}%
\special{pa 2392 2160}%
\special{fp}%
}}%
%
{\color[named]{Black}{%
\special{pn 8}%
\special{pa 2268 2376}%
\special{pa 2392 2160}%
\special{pa 2514 2376}%
\special{pa 2268 2376}%
\special{pa 2392 2160}%
\special{fp}%
}}%
%
{\color[named]{Black}{%
\special{pn 8}%
\special{pa 2392 2160}%
\special{pa 2514 1948}%
\special{pa 2638 2160}%
\special{pa 2392 2160}%
\special{pa 2514 1948}%
\special{fp}%
}}%
%
{\color[named]{Black}{%
\special{pn 8}%
\special{pa 3250 2376}%
\special{pa 2760 2376}%
\special{pa 3004 1948}%
\special{pa 3250 2376}%
\special{pa 2760 2376}%
\special{fp}%
}}%
%
{\color[named]{Black}{%
\special{pn 8}%
\special{pa 3004 2376}%
\special{pa 3128 2162}%
\special{pa 3250 2376}%
\special{pa 3004 2376}%
\special{pa 3128 2162}%
\special{fp}%
}}%
%
{\color[named]{Black}{%
\special{pn 8}%
\special{pa 3004 2376}%
\special{pa 3128 2162}%
\special{pa 3250 2376}%
\special{pa 3004 2376}%
\special{pa 3128 2162}%
\special{fp}%
}}%
%
{\color[named]{Black}{%
\special{pn 8}%
\special{pa 3004 2376}%
\special{pa 3128 2162}%
\special{pa 3250 2376}%
\special{pa 3004 2376}%
\special{pa 3128 2162}%
\special{fp}%
}}%
%
{\color[named]{Black}{%
\special{pn 8}%
\special{pa 3004 2376}%
\special{pa 3128 2162}%
\special{pa 3250 2376}%
\special{pa 3004 2376}%
\special{pa 3128 2162}%
\special{fp}%
}}%
%
{\color[named]{Black}{%
\special{pn 8}%
\special{pa 2760 2376}%
\special{pa 2882 2162}%
\special{pa 3006 2376}%
\special{pa 2760 2376}%
\special{pa 2882 2162}%
\special{fp}%
}}%
%
{\color[named]{Black}{%
\special{pn 8}%
\special{pa 2760 2376}%
\special{pa 2882 2162}%
\special{pa 3006 2376}%
\special{pa 2760 2376}%
\special{pa 2882 2162}%
\special{fp}%
}}%
%
{\color[named]{Black}{%
\special{pn 8}%
\special{pa 2760 2376}%
\special{pa 2882 2162}%
\special{pa 3006 2376}%
\special{pa 2760 2376}%
\special{pa 2882 2162}%
\special{fp}%
}}%
%
{\color[named]{Black}{%
\special{pn 8}%
\special{pa 2760 2376}%
\special{pa 2882 2162}%
\special{pa 3006 2376}%
\special{pa 2760 2376}%
\special{pa 2882 2162}%
\special{fp}%
}}%
%
{\color[named]{Black}{%
\special{pn 8}%
\special{pa 2882 2162}%
\special{pa 3004 1948}%
\special{pa 3128 2162}%
\special{pa 2882 2162}%
\special{pa 3004 1948}%
\special{fp}%
}}%
%
{\color[named]{Black}{%
\special{pn 8}%
\special{pa 3004 1948}%
\special{pa 2514 1948}%
\special{pa 2760 1524}%
\special{pa 3004 1948}%
\special{pa 2514 1948}%
\special{fp}%
}}%
%
{\color[named]{Black}{%
\special{pn 8}%
\special{pa 2514 1948}%
\special{pa 2636 1736}%
\special{pa 2760 1948}%
\special{pa 2514 1948}%
\special{pa 2636 1736}%
\special{fp}%
}}%
%
{\color[named]{Black}{%
\special{pn 8}%
\special{pa 2760 1948}%
\special{pa 2882 1736}%
\special{pa 3004 1948}%
\special{pa 2760 1948}%
\special{pa 2882 1736}%
\special{fp}%
}}%
%
{\color[named]{Black}{%
\special{pn 8}%
\special{pa 2636 1736}%
\special{pa 2760 1522}%
\special{pa 2882 1736}%
\special{pa 2636 1736}%
\special{pa 2760 1522}%
\special{fp}%
}}%
%
{\color[named]{Black}{%
\special{pn 8}%
\special{pa 3742 2376}%
\special{pa 3250 2376}%
\special{pa 3494 1950}%
\special{pa 3742 2376}%
\special{pa 3250 2376}%
\special{fp}%
}}%
%
{\color[named]{Black}{%
\special{pn 8}%
\special{pa 3494 2376}%
\special{pa 3618 2160}%
\special{pa 3742 2376}%
\special{pa 3494 2376}%
\special{pa 3618 2160}%
\special{fp}%
}}%
%
{\color[named]{Black}{%
\special{pn 8}%
\special{pa 3494 2376}%
\special{pa 3618 2160}%
\special{pa 3742 2376}%
\special{pa 3494 2376}%
\special{pa 3618 2160}%
\special{fp}%
}}%
%
{\color[named]{Black}{%
\special{pn 8}%
\special{pa 3494 2376}%
\special{pa 3618 2160}%
\special{pa 3742 2376}%
\special{pa 3494 2376}%
\special{pa 3618 2160}%
\special{fp}%
}}%
%
{\color[named]{Black}{%
\special{pn 8}%
\special{pa 3494 2376}%
\special{pa 3618 2160}%
\special{pa 3742 2376}%
\special{pa 3494 2376}%
\special{pa 3618 2160}%
\special{fp}%
}}%
%
{\color[named]{Black}{%
\special{pn 8}%
\special{pa 3250 2376}%
\special{pa 3372 2160}%
\special{pa 3496 2376}%
\special{pa 3250 2376}%
\special{pa 3372 2160}%
\special{fp}%
}}%
%
{\color[named]{Black}{%
\special{pn 8}%
\special{pa 3250 2376}%
\special{pa 3372 2160}%
\special{pa 3496 2376}%
\special{pa 3250 2376}%
\special{pa 3372 2160}%
\special{fp}%
}}%
%
{\color[named]{Black}{%
\special{pn 8}%
\special{pa 3250 2376}%
\special{pa 3372 2160}%
\special{pa 3496 2376}%
\special{pa 3250 2376}%
\special{pa 3372 2160}%
\special{fp}%
}}%
%
{\color[named]{Black}{%
\special{pn 8}%
\special{pa 3250 2376}%
\special{pa 3372 2160}%
\special{pa 3496 2376}%
\special{pa 3250 2376}%
\special{pa 3372 2160}%
\special{fp}%
}}%
%
{\color[named]{Black}{%
\special{pn 8}%
\special{pa 3372 2160}%
\special{pa 3494 1948}%
\special{pa 3618 2160}%
\special{pa 3372 2160}%
\special{pa 3494 1948}%
\special{fp}%
}}%
%
{\color[named]{Black}{%
\special{pn 8}%
\special{pa 4232 2376}%
\special{pa 3740 2376}%
\special{pa 3986 1948}%
\special{pa 4232 2376}%
\special{pa 3740 2376}%
\special{fp}%
}}%
%
{\color[named]{Black}{%
\special{pn 8}%
\special{pa 3986 2376}%
\special{pa 4108 2162}%
\special{pa 4232 2376}%
\special{pa 3986 2376}%
\special{pa 4108 2162}%
\special{fp}%
}}%
%
{\color[named]{Black}{%
\special{pn 8}%
\special{pa 3986 2376}%
\special{pa 4108 2162}%
\special{pa 4232 2376}%
\special{pa 3986 2376}%
\special{pa 4108 2162}%
\special{fp}%
}}%
%
{\color[named]{Black}{%
\special{pn 8}%
\special{pa 3986 2376}%
\special{pa 4108 2162}%
\special{pa 4232 2376}%
\special{pa 3986 2376}%
\special{pa 4108 2162}%
\special{fp}%
}}%
%
{\color[named]{Black}{%
\special{pn 8}%
\special{pa 3986 2376}%
\special{pa 4108 2162}%
\special{pa 4232 2376}%
\special{pa 3986 2376}%
\special{pa 4108 2162}%
\special{fp}%
}}%
%
{\color[named]{Black}{%
\special{pn 8}%
\special{pa 3740 2376}%
\special{pa 3864 2162}%
\special{pa 3986 2376}%
\special{pa 3740 2376}%
\special{pa 3864 2162}%
\special{fp}%
}}%
%
{\color[named]{Black}{%
\special{pn 8}%
\special{pa 3740 2376}%
\special{pa 3864 2162}%
\special{pa 3986 2376}%
\special{pa 3740 2376}%
\special{pa 3864 2162}%
\special{fp}%
}}%
%
{\color[named]{Black}{%
\special{pn 8}%
\special{pa 3740 2376}%
\special{pa 3864 2162}%
\special{pa 3986 2376}%
\special{pa 3740 2376}%
\special{pa 3864 2162}%
\special{fp}%
}}%
%
{\color[named]{Black}{%
\special{pn 8}%
\special{pa 3740 2376}%
\special{pa 3864 2162}%
\special{pa 3986 2376}%
\special{pa 3740 2376}%
\special{pa 3864 2162}%
\special{fp}%
}}%
%
{\color[named]{Black}{%
\special{pn 8}%
\special{pa 3864 2162}%
\special{pa 3986 1948}%
\special{pa 4108 2162}%
\special{pa 3864 2162}%
\special{pa 3986 1948}%
\special{fp}%
}}%
%
{\color[named]{Black}{%
\special{pn 8}%
\special{pa 3986 1948}%
\special{pa 3494 1948}%
\special{pa 3740 1524}%
\special{pa 3986 1948}%
\special{pa 3494 1948}%
\special{fp}%
}}%
%
{\color[named]{Black}{%
\special{pn 8}%
\special{pa 3494 1948}%
\special{pa 3618 1736}%
\special{pa 3740 1948}%
\special{pa 3494 1948}%
\special{pa 3618 1736}%
\special{fp}%
}}%
%
{\color[named]{Black}{%
\special{pn 8}%
\special{pa 3740 1948}%
\special{pa 3864 1736}%
\special{pa 3986 1948}%
\special{pa 3740 1948}%
\special{pa 3864 1736}%
\special{fp}%
}}%
%
{\color[named]{Black}{%
\special{pn 8}%
\special{pa 3616 1736}%
\special{pa 3740 1522}%
\special{pa 3864 1736}%
\special{pa 3616 1736}%
\special{pa 3740 1522}%
\special{fp}%
}}%
%
{\color[named]{Black}{%
\special{pn 8}%
\special{pa 3250 1522}%
\special{pa 2758 1522}%
\special{pa 3002 1096}%
\special{pa 3250 1522}%
\special{pa 2758 1522}%
\special{fp}%
}}%
%
{\color[named]{Black}{%
\special{pn 8}%
\special{pa 3002 1522}%
\special{pa 3126 1306}%
\special{pa 3250 1522}%
\special{pa 3002 1522}%
\special{pa 3126 1306}%
\special{fp}%
}}%
%
{\color[named]{Black}{%
\special{pn 8}%
\special{pa 3002 1522}%
\special{pa 3126 1306}%
\special{pa 3250 1522}%
\special{pa 3002 1522}%
\special{pa 3126 1306}%
\special{fp}%
}}%
%
{\color[named]{Black}{%
\special{pn 8}%
\special{pa 3002 1522}%
\special{pa 3126 1306}%
\special{pa 3250 1522}%
\special{pa 3002 1522}%
\special{pa 3126 1306}%
\special{fp}%
}}%
%
{\color[named]{Black}{%
\special{pn 8}%
\special{pa 3002 1522}%
\special{pa 3126 1306}%
\special{pa 3250 1522}%
\special{pa 3002 1522}%
\special{pa 3126 1306}%
\special{fp}%
}}%
%
{\color[named]{Black}{%
\special{pn 8}%
\special{pa 2758 1522}%
\special{pa 2880 1306}%
\special{pa 3004 1522}%
\special{pa 2758 1522}%
\special{pa 2880 1306}%
\special{fp}%
}}%
%
{\color[named]{Black}{%
\special{pn 8}%
\special{pa 2758 1522}%
\special{pa 2880 1306}%
\special{pa 3004 1522}%
\special{pa 2758 1522}%
\special{pa 2880 1306}%
\special{fp}%
}}%
%
{\color[named]{Black}{%
\special{pn 8}%
\special{pa 2758 1522}%
\special{pa 2880 1306}%
\special{pa 3004 1522}%
\special{pa 2758 1522}%
\special{pa 2880 1306}%
\special{fp}%
}}%
%
{\color[named]{Black}{%
\special{pn 8}%
\special{pa 2758 1522}%
\special{pa 2880 1306}%
\special{pa 3004 1522}%
\special{pa 2758 1522}%
\special{pa 2880 1306}%
\special{fp}%
}}%
%
{\color[named]{Black}{%
\special{pn 8}%
\special{pa 2880 1306}%
\special{pa 3002 1094}%
\special{pa 3126 1306}%
\special{pa 2880 1306}%
\special{pa 3002 1094}%
\special{fp}%
}}%
%
{\color[named]{Black}{%
\special{pn 8}%
\special{pa 3740 1522}%
\special{pa 3248 1522}%
\special{pa 3492 1094}%
\special{pa 3740 1522}%
\special{pa 3248 1522}%
\special{fp}%
}}%
%
{\color[named]{Black}{%
\special{pn 8}%
\special{pa 3492 1522}%
\special{pa 3616 1308}%
\special{pa 3740 1522}%
\special{pa 3492 1522}%
\special{pa 3616 1308}%
\special{fp}%
}}%
%
{\color[named]{Black}{%
\special{pn 8}%
\special{pa 3492 1522}%
\special{pa 3616 1308}%
\special{pa 3740 1522}%
\special{pa 3492 1522}%
\special{pa 3616 1308}%
\special{fp}%
}}%
%
{\color[named]{Black}{%
\special{pn 8}%
\special{pa 3492 1522}%
\special{pa 3616 1308}%
\special{pa 3740 1522}%
\special{pa 3492 1522}%
\special{pa 3616 1308}%
\special{fp}%
}}%
%
{\color[named]{Black}{%
\special{pn 8}%
\special{pa 3492 1522}%
\special{pa 3616 1308}%
\special{pa 3740 1522}%
\special{pa 3492 1522}%
\special{pa 3616 1308}%
\special{fp}%
}}%
%
{\color[named]{Black}{%
\special{pn 8}%
\special{pa 3248 1522}%
\special{pa 3370 1308}%
\special{pa 3494 1522}%
\special{pa 3248 1522}%
\special{pa 3370 1308}%
\special{fp}%
}}%
%
{\color[named]{Black}{%
\special{pn 8}%
\special{pa 3248 1522}%
\special{pa 3370 1308}%
\special{pa 3494 1522}%
\special{pa 3248 1522}%
\special{pa 3370 1308}%
\special{fp}%
}}%
%
{\color[named]{Black}{%
\special{pn 8}%
\special{pa 3248 1522}%
\special{pa 3370 1308}%
\special{pa 3494 1522}%
\special{pa 3248 1522}%
\special{pa 3370 1308}%
\special{fp}%
}}%
%
{\color[named]{Black}{%
\special{pn 8}%
\special{pa 3248 1522}%
\special{pa 3370 1308}%
\special{pa 3494 1522}%
\special{pa 3248 1522}%
\special{pa 3370 1308}%
\special{fp}%
}}%
%
{\color[named]{Black}{%
\special{pn 8}%
\special{pa 3370 1308}%
\special{pa 3492 1094}%
\special{pa 3616 1308}%
\special{pa 3370 1308}%
\special{pa 3492 1094}%
\special{fp}%
}}%
%
{\color[named]{Black}{%
\special{pn 8}%
\special{pa 3492 1094}%
\special{pa 3002 1094}%
\special{pa 3248 670}%
\special{pa 3492 1094}%
\special{pa 3002 1094}%
\special{fp}%
}}%
%
{\color[named]{Black}{%
\special{pn 8}%
\special{pa 3002 1094}%
\special{pa 3124 882}%
\special{pa 3248 1094}%
\special{pa 3002 1094}%
\special{pa 3124 882}%
\special{fp}%
}}%
%
{\color[named]{Black}{%
\special{pn 8}%
\special{pa 3248 1094}%
\special{pa 3370 882}%
\special{pa 3492 1094}%
\special{pa 3248 1094}%
\special{pa 3370 882}%
\special{fp}%
}}%
%
{\color[named]{Black}{%
\special{pn 8}%
\special{pa 3124 880}%
\special{pa 3248 666}%
\special{pa 3370 880}%
\special{pa 3124 880}%
\special{pa 3248 666}%
\special{fp}%
}}%
%
{\color[named]{Black}{%
\special{pn 8}%
\special{pa 798 2376}%
\special{pa 308 2376}%
\special{pa 552 1950}%
\special{pa 798 2376}%
\special{pa 308 2376}%
\special{fp}%
}}%
%
{\color[named]{Black}{%
\special{pn 8}%
\special{pa 552 2376}%
\special{pa 676 2160}%
\special{pa 798 2376}%
\special{pa 552 2376}%
\special{pa 676 2160}%
\special{fp}%
}}%
%
{\color[named]{Black}{%
\special{pn 8}%
\special{pa 552 2376}%
\special{pa 676 2160}%
\special{pa 798 2376}%
\special{pa 552 2376}%
\special{pa 676 2160}%
\special{fp}%
}}%
%
{\color[named]{Black}{%
\special{pn 8}%
\special{pa 552 2376}%
\special{pa 676 2160}%
\special{pa 798 2376}%
\special{pa 552 2376}%
\special{pa 676 2160}%
\special{fp}%
}}%
%
{\color[named]{Black}{%
\special{pn 8}%
\special{pa 552 2376}%
\special{pa 676 2160}%
\special{pa 798 2376}%
\special{pa 552 2376}%
\special{pa 676 2160}%
\special{fp}%
}}%
%
{\color[named]{Black}{%
\special{pn 8}%
\special{pa 308 2376}%
\special{pa 430 2160}%
\special{pa 552 2376}%
\special{pa 308 2376}%
\special{pa 430 2160}%
\special{fp}%
}}%
%
{\color[named]{Black}{%
\special{pn 8}%
\special{pa 308 2376}%
\special{pa 430 2160}%
\special{pa 552 2376}%
\special{pa 308 2376}%
\special{pa 430 2160}%
\special{fp}%
}}%
%
{\color[named]{Black}{%
\special{pn 8}%
\special{pa 308 2376}%
\special{pa 430 2160}%
\special{pa 552 2376}%
\special{pa 308 2376}%
\special{pa 430 2160}%
\special{fp}%
}}%
%
{\color[named]{Black}{%
\special{pn 8}%
\special{pa 308 2376}%
\special{pa 430 2160}%
\special{pa 552 2376}%
\special{pa 308 2376}%
\special{pa 430 2160}%
\special{fp}%
}}%
%
{\color[named]{Black}{%
\special{pn 8}%
\special{pa 430 2160}%
\special{pa 552 1948}%
\special{pa 676 2160}%
\special{pa 430 2160}%
\special{pa 552 1948}%
\special{fp}%
}}%
%
{\color[named]{Black}{%
\special{pn 8}%
\special{pa 1290 2376}%
\special{pa 798 2376}%
\special{pa 1044 1948}%
\special{pa 1290 2376}%
\special{pa 798 2376}%
\special{fp}%
}}%
%
{\color[named]{Black}{%
\special{pn 8}%
\special{pa 1044 2376}%
\special{pa 1166 2162}%
\special{pa 1290 2376}%
\special{pa 1044 2376}%
\special{pa 1166 2162}%
\special{fp}%
}}%
%
{\color[named]{Black}{%
\special{pn 8}%
\special{pa 1044 2376}%
\special{pa 1166 2162}%
\special{pa 1290 2376}%
\special{pa 1044 2376}%
\special{pa 1166 2162}%
\special{fp}%
}}%
%
{\color[named]{Black}{%
\special{pn 8}%
\special{pa 1044 2376}%
\special{pa 1166 2162}%
\special{pa 1290 2376}%
\special{pa 1044 2376}%
\special{pa 1166 2162}%
\special{fp}%
}}%
%
{\color[named]{Black}{%
\special{pn 8}%
\special{pa 1044 2376}%
\special{pa 1166 2162}%
\special{pa 1290 2376}%
\special{pa 1044 2376}%
\special{pa 1166 2162}%
\special{fp}%
}}%
%
{\color[named]{Black}{%
\special{pn 8}%
\special{pa 798 2376}%
\special{pa 922 2162}%
\special{pa 1044 2376}%
\special{pa 798 2376}%
\special{pa 922 2162}%
\special{fp}%
}}%
%
{\color[named]{Black}{%
\special{pn 8}%
\special{pa 798 2376}%
\special{pa 922 2162}%
\special{pa 1044 2376}%
\special{pa 798 2376}%
\special{pa 922 2162}%
\special{fp}%
}}%
%
{\color[named]{Black}{%
\special{pn 8}%
\special{pa 798 2376}%
\special{pa 922 2162}%
\special{pa 1044 2376}%
\special{pa 798 2376}%
\special{pa 922 2162}%
\special{fp}%
}}%
%
{\color[named]{Black}{%
\special{pn 8}%
\special{pa 798 2376}%
\special{pa 922 2162}%
\special{pa 1044 2376}%
\special{pa 798 2376}%
\special{pa 922 2162}%
\special{fp}%
}}%
%
{\color[named]{Black}{%
\special{pn 8}%
\special{pa 922 2162}%
\special{pa 1044 1948}%
\special{pa 1166 2162}%
\special{pa 922 2162}%
\special{pa 1044 1948}%
\special{fp}%
}}%
%
{\color[named]{Black}{%
\special{pn 8}%
\special{pa 1044 1948}%
\special{pa 552 1948}%
\special{pa 798 1524}%
\special{pa 1044 1948}%
\special{pa 552 1948}%
\special{fp}%
}}%
%
{\color[named]{Black}{%
\special{pn 8}%
\special{pa 552 1948}%
\special{pa 676 1736}%
\special{pa 798 1948}%
\special{pa 552 1948}%
\special{pa 676 1736}%
\special{fp}%
}}%
%
{\color[named]{Black}{%
\special{pn 8}%
\special{pa 798 1948}%
\special{pa 922 1736}%
\special{pa 1044 1948}%
\special{pa 798 1948}%
\special{pa 922 1736}%
\special{fp}%
}}%
%
{\color[named]{Black}{%
\special{pn 8}%
\special{pa 674 1736}%
\special{pa 798 1522}%
\special{pa 922 1736}%
\special{pa 674 1736}%
\special{pa 798 1522}%
\special{fp}%
}}%
%
{\color[named]{Black}{%
\special{pn 8}%
\special{pa 1780 2376}%
\special{pa 1288 2376}%
\special{pa 1534 1950}%
\special{pa 1780 2376}%
\special{pa 1288 2376}%
\special{fp}%
}}%
%
{\color[named]{Black}{%
\special{pn 8}%
\special{pa 1534 2376}%
\special{pa 1656 2160}%
\special{pa 1780 2376}%
\special{pa 1534 2376}%
\special{pa 1656 2160}%
\special{fp}%
}}%
%
{\color[named]{Black}{%
\special{pn 8}%
\special{pa 1534 2376}%
\special{pa 1656 2160}%
\special{pa 1780 2376}%
\special{pa 1534 2376}%
\special{pa 1656 2160}%
\special{fp}%
}}%
%
{\color[named]{Black}{%
\special{pn 8}%
\special{pa 1534 2376}%
\special{pa 1656 2160}%
\special{pa 1780 2376}%
\special{pa 1534 2376}%
\special{pa 1656 2160}%
\special{fp}%
}}%
%
{\color[named]{Black}{%
\special{pn 8}%
\special{pa 1534 2376}%
\special{pa 1656 2160}%
\special{pa 1780 2376}%
\special{pa 1534 2376}%
\special{pa 1656 2160}%
\special{fp}%
}}%
%
{\color[named]{Black}{%
\special{pn 8}%
\special{pa 1288 2376}%
\special{pa 1410 2160}%
\special{pa 1534 2376}%
\special{pa 1288 2376}%
\special{pa 1410 2160}%
\special{fp}%
}}%
%
{\color[named]{Black}{%
\special{pn 8}%
\special{pa 1288 2376}%
\special{pa 1410 2160}%
\special{pa 1534 2376}%
\special{pa 1288 2376}%
\special{pa 1410 2160}%
\special{fp}%
}}%
%
{\color[named]{Black}{%
\special{pn 8}%
\special{pa 1288 2376}%
\special{pa 1410 2160}%
\special{pa 1534 2376}%
\special{pa 1288 2376}%
\special{pa 1410 2160}%
\special{fp}%
}}%
%
{\color[named]{Black}{%
\special{pn 8}%
\special{pa 1288 2376}%
\special{pa 1410 2160}%
\special{pa 1534 2376}%
\special{pa 1288 2376}%
\special{pa 1410 2160}%
\special{fp}%
}}%
%
{\color[named]{Black}{%
\special{pn 8}%
\special{pa 1410 2160}%
\special{pa 1534 1948}%
\special{pa 1656 2160}%
\special{pa 1410 2160}%
\special{pa 1534 1948}%
\special{fp}%
}}%
%
{\color[named]{Black}{%
\special{pn 8}%
\special{pa 2270 2376}%
\special{pa 1778 2376}%
\special{pa 2024 1948}%
\special{pa 2270 2376}%
\special{pa 1778 2376}%
\special{fp}%
}}%
%
{\color[named]{Black}{%
\special{pn 8}%
\special{pa 2024 2376}%
\special{pa 2148 2162}%
\special{pa 2270 2376}%
\special{pa 2024 2376}%
\special{pa 2148 2162}%
\special{fp}%
}}%
%
{\color[named]{Black}{%
\special{pn 8}%
\special{pa 2024 2376}%
\special{pa 2148 2162}%
\special{pa 2270 2376}%
\special{pa 2024 2376}%
\special{pa 2148 2162}%
\special{fp}%
}}%
%
{\color[named]{Black}{%
\special{pn 8}%
\special{pa 2024 2376}%
\special{pa 2148 2162}%
\special{pa 2270 2376}%
\special{pa 2024 2376}%
\special{pa 2148 2162}%
\special{fp}%
}}%
%
{\color[named]{Black}{%
\special{pn 8}%
\special{pa 2024 2376}%
\special{pa 2148 2162}%
\special{pa 2270 2376}%
\special{pa 2024 2376}%
\special{pa 2148 2162}%
\special{fp}%
}}%
%
{\color[named]{Black}{%
\special{pn 8}%
\special{pa 1778 2376}%
\special{pa 1902 2162}%
\special{pa 2024 2376}%
\special{pa 1778 2376}%
\special{pa 1902 2162}%
\special{fp}%
}}%
%
{\color[named]{Black}{%
\special{pn 8}%
\special{pa 1778 2376}%
\special{pa 1902 2162}%
\special{pa 2024 2376}%
\special{pa 1778 2376}%
\special{pa 1902 2162}%
\special{fp}%
}}%
%
{\color[named]{Black}{%
\special{pn 8}%
\special{pa 1778 2376}%
\special{pa 1902 2162}%
\special{pa 2024 2376}%
\special{pa 1778 2376}%
\special{pa 1902 2162}%
\special{fp}%
}}%
%
{\color[named]{Black}{%
\special{pn 8}%
\special{pa 1778 2376}%
\special{pa 1902 2162}%
\special{pa 2024 2376}%
\special{pa 1778 2376}%
\special{pa 1902 2162}%
\special{fp}%
}}%
%
{\color[named]{Black}{%
\special{pn 8}%
\special{pa 1902 2162}%
\special{pa 2024 1948}%
\special{pa 2148 2162}%
\special{pa 1902 2162}%
\special{pa 2024 1948}%
\special{fp}%
}}%
%
{\color[named]{Black}{%
\special{pn 8}%
\special{pa 2024 1948}%
\special{pa 1534 1948}%
\special{pa 1778 1524}%
\special{pa 2024 1948}%
\special{pa 1534 1948}%
\special{fp}%
}}%
%
{\color[named]{Black}{%
\special{pn 8}%
\special{pa 1534 1948}%
\special{pa 1656 1736}%
\special{pa 1778 1948}%
\special{pa 1534 1948}%
\special{pa 1656 1736}%
\special{fp}%
}}%
%
{\color[named]{Black}{%
\special{pn 8}%
\special{pa 1778 1948}%
\special{pa 1902 1736}%
\special{pa 2024 1948}%
\special{pa 1778 1948}%
\special{pa 1902 1736}%
\special{fp}%
}}%
%
{\color[named]{Black}{%
\special{pn 8}%
\special{pa 1656 1736}%
\special{pa 1778 1522}%
\special{pa 1902 1736}%
\special{pa 1656 1736}%
\special{pa 1778 1522}%
\special{fp}%
}}%
%
{\color[named]{Black}{%
\special{pn 8}%
\special{pa 1288 1522}%
\special{pa 796 1522}%
\special{pa 1040 1096}%
\special{pa 1288 1522}%
\special{pa 796 1522}%
\special{fp}%
}}%
%
{\color[named]{Black}{%
\special{pn 8}%
\special{pa 1040 1522}%
\special{pa 1164 1306}%
\special{pa 1288 1522}%
\special{pa 1040 1522}%
\special{pa 1164 1306}%
\special{fp}%
}}%
%
{\color[named]{Black}{%
\special{pn 8}%
\special{pa 1040 1522}%
\special{pa 1164 1306}%
\special{pa 1288 1522}%
\special{pa 1040 1522}%
\special{pa 1164 1306}%
\special{fp}%
}}%
%
{\color[named]{Black}{%
\special{pn 8}%
\special{pa 1040 1522}%
\special{pa 1164 1306}%
\special{pa 1288 1522}%
\special{pa 1040 1522}%
\special{pa 1164 1306}%
\special{fp}%
}}%
%
{\color[named]{Black}{%
\special{pn 8}%
\special{pa 1040 1522}%
\special{pa 1164 1306}%
\special{pa 1288 1522}%
\special{pa 1040 1522}%
\special{pa 1164 1306}%
\special{fp}%
}}%
%
{\color[named]{Black}{%
\special{pn 8}%
\special{pa 796 1522}%
\special{pa 920 1306}%
\special{pa 1042 1522}%
\special{pa 796 1522}%
\special{pa 920 1306}%
\special{fp}%
}}%
%
{\color[named]{Black}{%
\special{pn 8}%
\special{pa 796 1522}%
\special{pa 920 1306}%
\special{pa 1042 1522}%
\special{pa 796 1522}%
\special{pa 920 1306}%
\special{fp}%
}}%
%
{\color[named]{Black}{%
\special{pn 8}%
\special{pa 796 1522}%
\special{pa 920 1306}%
\special{pa 1042 1522}%
\special{pa 796 1522}%
\special{pa 920 1306}%
\special{fp}%
}}%
%
{\color[named]{Black}{%
\special{pn 8}%
\special{pa 796 1522}%
\special{pa 920 1306}%
\special{pa 1042 1522}%
\special{pa 796 1522}%
\special{pa 920 1306}%
\special{fp}%
}}%
%
{\color[named]{Black}{%
\special{pn 8}%
\special{pa 920 1306}%
\special{pa 1040 1094}%
\special{pa 1164 1306}%
\special{pa 920 1306}%
\special{pa 1040 1094}%
\special{fp}%
}}%
%
{\color[named]{Black}{%
\special{pn 8}%
\special{pa 1778 1522}%
\special{pa 1286 1522}%
\special{pa 1532 1094}%
\special{pa 1778 1522}%
\special{pa 1286 1522}%
\special{fp}%
}}%
%
{\color[named]{Black}{%
\special{pn 8}%
\special{pa 1532 1522}%
\special{pa 1656 1308}%
\special{pa 1778 1522}%
\special{pa 1532 1522}%
\special{pa 1656 1308}%
\special{fp}%
}}%
%
{\color[named]{Black}{%
\special{pn 8}%
\special{pa 1532 1522}%
\special{pa 1656 1308}%
\special{pa 1778 1522}%
\special{pa 1532 1522}%
\special{pa 1656 1308}%
\special{fp}%
}}%
%
{\color[named]{Black}{%
\special{pn 8}%
\special{pa 1532 1522}%
\special{pa 1656 1308}%
\special{pa 1778 1522}%
\special{pa 1532 1522}%
\special{pa 1656 1308}%
\special{fp}%
}}%
%
{\color[named]{Black}{%
\special{pn 8}%
\special{pa 1532 1522}%
\special{pa 1656 1308}%
\special{pa 1778 1522}%
\special{pa 1532 1522}%
\special{pa 1656 1308}%
\special{fp}%
}}%
%
{\color[named]{Black}{%
\special{pn 8}%
\special{pa 1286 1522}%
\special{pa 1410 1308}%
\special{pa 1532 1522}%
\special{pa 1286 1522}%
\special{pa 1410 1308}%
\special{fp}%
}}%
%
{\color[named]{Black}{%
\special{pn 8}%
\special{pa 1286 1522}%
\special{pa 1410 1308}%
\special{pa 1532 1522}%
\special{pa 1286 1522}%
\special{pa 1410 1308}%
\special{fp}%
}}%
%
{\color[named]{Black}{%
\special{pn 8}%
\special{pa 1286 1522}%
\special{pa 1410 1308}%
\special{pa 1532 1522}%
\special{pa 1286 1522}%
\special{pa 1410 1308}%
\special{fp}%
}}%
%
{\color[named]{Black}{%
\special{pn 8}%
\special{pa 1286 1522}%
\special{pa 1410 1308}%
\special{pa 1532 1522}%
\special{pa 1286 1522}%
\special{pa 1410 1308}%
\special{fp}%
}}%
%
{\color[named]{Black}{%
\special{pn 8}%
\special{pa 1410 1308}%
\special{pa 1532 1094}%
\special{pa 1656 1308}%
\special{pa 1410 1308}%
\special{pa 1532 1094}%
\special{fp}%
}}%
%
{\color[named]{Black}{%
\special{pn 8}%
\special{pa 1532 1094}%
\special{pa 1040 1094}%
\special{pa 1286 670}%
\special{pa 1532 1094}%
\special{pa 1040 1094}%
\special{fp}%
}}%
%
{\color[named]{Black}{%
\special{pn 8}%
\special{pa 1040 1094}%
\special{pa 1164 882}%
\special{pa 1286 1094}%
\special{pa 1040 1094}%
\special{pa 1164 882}%
\special{fp}%
}}%
%
{\color[named]{Black}{%
\special{pn 8}%
\special{pa 1286 1094}%
\special{pa 1410 882}%
\special{pa 1532 1094}%
\special{pa 1286 1094}%
\special{pa 1410 882}%
\special{fp}%
}}%
%
{\color[named]{Black}{%
\special{pn 8}%
\special{pa 1164 880}%
\special{pa 1286 666}%
\special{pa 1410 880}%
\special{pa 1164 880}%
\special{pa 1286 666}%
\special{fp}%
}}%
%
{\color[named]{Black}{%
\special{pn 8}%
\special{pa 4226 2372}%
\special{pa 4470 2372}%
\special{fp}%
\special{pa 320 2372}%
\special{pa 76 2372}%
\special{fp}%
\special{pa 1284 664}%
\special{pa 1040 664}%
\special{fp}%
\special{pa 3246 664}%
\special{pa 3490 664}%
\special{fp}%
}}%
%
{\color[named]{Black}{%
\special{pn 8}%
\special{pa 308 2372}%
\special{pa 124 2054}%
\special{fp}%
}}%
%
{\color[named]{Black}{%
\special{pn 8}%
\special{pa 1284 664}%
\special{pa 1100 346}%
\special{fp}%
}}%
%
{\color[named]{Black}{%
\special{pn 8}%
\special{pa 4230 2372}%
\special{pa 4414 2054}%
\special{fp}%
}}%
%
{\color[named]{Black}{%
\special{pn 8}%
\special{pa 3246 664}%
\special{pa 3428 346}%
\special{fp}%
}}%
\put(22.1300,-24.2500){\makebox(0,0)[lt]{$O$}}%
\put(12.2000,-27.7000){\makebox(0,0)[lb]{$F_0$}}%
\end{picture}%

%% file: Fig4.tex
\unitlength 0.1in
\begin{picture}( 42.0200, 19.0900)(  6.0000,-25.8900)
%
{\color[named]{Black}{%
\special{pn 4}%
\special{pa 1454 1472}%
\special{pa 728 1472}%
\special{pa 1092 844}%
\special{pa 1454 1472}%
\special{pa 728 1472}%
\special{fp}%
}}%
%
{\color[named]{Black}{%
\special{pn 4}%
\special{pa 1092 1472}%
\special{pa 1272 1158}%
\special{pa 1454 1472}%
\special{pa 1092 1472}%
\special{pa 1272 1158}%
\special{fp}%
}}%
%
{\color[named]{Black}{%
\special{pn 4}%
\special{pa 1092 1472}%
\special{pa 1272 1158}%
\special{pa 1454 1472}%
\special{pa 1092 1472}%
\special{pa 1272 1158}%
\special{fp}%
}}%
%
{\color[named]{Black}{%
\special{pn 4}%
\special{pa 1092 1472}%
\special{pa 1272 1158}%
\special{pa 1454 1472}%
\special{pa 1092 1472}%
\special{pa 1272 1158}%
\special{fp}%
}}%
%
{\color[named]{Black}{%
\special{pn 4}%
\special{pa 1092 1472}%
\special{pa 1272 1158}%
\special{pa 1454 1472}%
\special{pa 1092 1472}%
\special{pa 1272 1158}%
\special{fp}%
}}%
%
{\color[named]{Black}{%
\special{pn 4}%
\special{pa 1092 1472}%
\special{pa 1272 1158}%
\special{pa 1454 1472}%
\special{pa 1092 1472}%
\special{pa 1272 1158}%
\special{fp}%
}}%
%
{\color[named]{Black}{%
\special{pn 4}%
\special{pa 728 1472}%
\special{pa 910 1158}%
\special{pa 1092 1472}%
\special{pa 728 1472}%
\special{pa 910 1158}%
\special{fp}%
}}%
%
{\color[named]{Black}{%
\special{pn 4}%
\special{pa 728 1472}%
\special{pa 910 1158}%
\special{pa 1092 1472}%
\special{pa 728 1472}%
\special{pa 910 1158}%
\special{fp}%
}}%
%
{\color[named]{Black}{%
\special{pn 4}%
\special{pa 728 1472}%
\special{pa 910 1158}%
\special{pa 1092 1472}%
\special{pa 728 1472}%
\special{pa 910 1158}%
\special{fp}%
}}%
%
{\color[named]{Black}{%
\special{pn 4}%
\special{pa 728 1472}%
\special{pa 910 1158}%
\special{pa 1092 1472}%
\special{pa 728 1472}%
\special{pa 910 1158}%
\special{fp}%
}}%
%
{\color[named]{Black}{%
\special{pn 4}%
\special{pa 910 1158}%
\special{pa 1092 844}%
\special{pa 1272 1158}%
\special{pa 910 1158}%
\special{pa 1092 844}%
\special{fp}%
}}%
\put(10.9000,-8.1000){\makebox(0,0)[rb]{$a_1$}}%
\put(14.1000,-15.8000){\makebox(0,0){$b_1$}}%
%
{\color[named]{Black}{%
\special{pn 4}%
\special{pa 2296 1472}%
\special{pa 1570 1472}%
\special{pa 1934 844}%
\special{pa 2296 1472}%
\special{pa 1570 1472}%
\special{fp}%
}}%
%
{\color[named]{Black}{%
\special{pn 4}%
\special{pa 1934 1472}%
\special{pa 2116 1158}%
\special{pa 2296 1472}%
\special{pa 1934 1472}%
\special{pa 2116 1158}%
\special{fp}%
}}%
%
{\color[named]{Black}{%
\special{pn 4}%
\special{pa 1934 1472}%
\special{pa 2116 1158}%
\special{pa 2296 1472}%
\special{pa 1934 1472}%
\special{pa 2116 1158}%
\special{fp}%
}}%
%
{\color[named]{Black}{%
\special{pn 4}%
\special{pa 1934 1472}%
\special{pa 2116 1158}%
\special{pa 2296 1472}%
\special{pa 1934 1472}%
\special{pa 2116 1158}%
\special{fp}%
}}%
%
{\color[named]{Black}{%
\special{pn 4}%
\special{pa 1934 1472}%
\special{pa 2116 1158}%
\special{pa 2296 1472}%
\special{pa 1934 1472}%
\special{pa 2116 1158}%
\special{fp}%
}}%
%
{\color[named]{Black}{%
\special{pn 4}%
\special{pa 1934 1472}%
\special{pa 2116 1158}%
\special{pa 2296 1472}%
\special{pa 1934 1472}%
\special{pa 2116 1158}%
\special{fp}%
}}%
%
{\color[named]{Black}{%
\special{pn 4}%
\special{pa 1570 1472}%
\special{pa 1752 1158}%
\special{pa 1934 1472}%
\special{pa 1570 1472}%
\special{pa 1752 1158}%
\special{fp}%
}}%
%
{\color[named]{Black}{%
\special{pn 4}%
\special{pa 1570 1472}%
\special{pa 1752 1158}%
\special{pa 1934 1472}%
\special{pa 1570 1472}%
\special{pa 1752 1158}%
\special{fp}%
}}%
%
{\color[named]{Black}{%
\special{pn 4}%
\special{pa 1570 1472}%
\special{pa 1752 1158}%
\special{pa 1934 1472}%
\special{pa 1570 1472}%
\special{pa 1752 1158}%
\special{fp}%
}}%
%
{\color[named]{Black}{%
\special{pn 4}%
\special{pa 1570 1472}%
\special{pa 1752 1158}%
\special{pa 1934 1472}%
\special{pa 1570 1472}%
\special{pa 1752 1158}%
\special{fp}%
}}%
%
{\color[named]{Black}{%
\special{pn 4}%
\special{pa 1752 1158}%
\special{pa 1934 844}%
\special{pa 2116 1158}%
\special{pa 1752 1158}%
\special{pa 1934 844}%
\special{fp}%
}}%
\put(6.0000,-16.3000){\makebox(0,0)[lb]{$O$}}%
%
{\color[named]{Black}{%
\special{pn 20}%
\special{pa 728 1472}%
\special{pa 1090 844}%
\special{fp}%
}}%
%
{\color[named]{Black}{%
\special{pn 4}%
\special{pa 3132 1472}%
\special{pa 2406 1472}%
\special{pa 2768 844}%
\special{pa 3132 1472}%
\special{pa 2406 1472}%
\special{fp}%
}}%
%
{\color[named]{Black}{%
\special{pn 4}%
\special{pa 2768 1472}%
\special{pa 2950 1158}%
\special{pa 3132 1472}%
\special{pa 2768 1472}%
\special{pa 2950 1158}%
\special{fp}%
}}%
%
{\color[named]{Black}{%
\special{pn 4}%
\special{pa 2768 1472}%
\special{pa 2950 1158}%
\special{pa 3132 1472}%
\special{pa 2768 1472}%
\special{pa 2950 1158}%
\special{fp}%
}}%
%
{\color[named]{Black}{%
\special{pn 4}%
\special{pa 2768 1472}%
\special{pa 2950 1158}%
\special{pa 3132 1472}%
\special{pa 2768 1472}%
\special{pa 2950 1158}%
\special{fp}%
}}%
%
{\color[named]{Black}{%
\special{pn 4}%
\special{pa 2768 1472}%
\special{pa 2950 1158}%
\special{pa 3132 1472}%
\special{pa 2768 1472}%
\special{pa 2950 1158}%
\special{fp}%
}}%
%
{\color[named]{Black}{%
\special{pn 4}%
\special{pa 2768 1472}%
\special{pa 2950 1158}%
\special{pa 3132 1472}%
\special{pa 2768 1472}%
\special{pa 2950 1158}%
\special{fp}%
}}%
%
{\color[named]{Black}{%
\special{pn 4}%
\special{pa 2406 1472}%
\special{pa 2588 1158}%
\special{pa 2768 1472}%
\special{pa 2406 1472}%
\special{pa 2588 1158}%
\special{fp}%
}}%
%
{\color[named]{Black}{%
\special{pn 4}%
\special{pa 2406 1472}%
\special{pa 2588 1158}%
\special{pa 2768 1472}%
\special{pa 2406 1472}%
\special{pa 2588 1158}%
\special{fp}%
}}%
%
{\color[named]{Black}{%
\special{pn 4}%
\special{pa 2406 1472}%
\special{pa 2588 1158}%
\special{pa 2768 1472}%
\special{pa 2406 1472}%
\special{pa 2588 1158}%
\special{fp}%
}}%
%
{\color[named]{Black}{%
\special{pn 4}%
\special{pa 2406 1472}%
\special{pa 2588 1158}%
\special{pa 2768 1472}%
\special{pa 2406 1472}%
\special{pa 2588 1158}%
\special{fp}%
}}%
%
{\color[named]{Black}{%
\special{pn 4}%
\special{pa 2588 1158}%
\special{pa 2768 844}%
\special{pa 2950 1158}%
\special{pa 2588 1158}%
\special{pa 2768 844}%
\special{fp}%
}}%
%
{\color[named]{Black}{%
\special{pn 4}%
\special{pa 3098 2590}%
\special{pa 2374 2590}%
\special{pa 2736 1962}%
\special{pa 3098 2590}%
\special{pa 2374 2590}%
\special{fp}%
}}%
%
{\color[named]{Black}{%
\special{pn 4}%
\special{pa 2736 2590}%
\special{pa 2916 2274}%
\special{pa 3098 2590}%
\special{pa 2736 2590}%
\special{pa 2916 2274}%
\special{fp}%
}}%
%
{\color[named]{Black}{%
\special{pn 4}%
\special{pa 2736 2590}%
\special{pa 2916 2274}%
\special{pa 3098 2590}%
\special{pa 2736 2590}%
\special{pa 2916 2274}%
\special{fp}%
}}%
%
{\color[named]{Black}{%
\special{pn 4}%
\special{pa 2736 2590}%
\special{pa 2916 2274}%
\special{pa 3098 2590}%
\special{pa 2736 2590}%
\special{pa 2916 2274}%
\special{fp}%
}}%
%
{\color[named]{Black}{%
\special{pn 4}%
\special{pa 2736 2590}%
\special{pa 2916 2274}%
\special{pa 3098 2590}%
\special{pa 2736 2590}%
\special{pa 2916 2274}%
\special{fp}%
}}%
%
{\color[named]{Black}{%
\special{pn 4}%
\special{pa 2736 2590}%
\special{pa 2916 2274}%
\special{pa 3098 2590}%
\special{pa 2736 2590}%
\special{pa 2916 2274}%
\special{fp}%
}}%
%
{\color[named]{Black}{%
\special{pn 4}%
\special{pa 2374 2590}%
\special{pa 2554 2274}%
\special{pa 2736 2590}%
\special{pa 2374 2590}%
\special{pa 2554 2274}%
\special{fp}%
}}%
%
{\color[named]{Black}{%
\special{pn 4}%
\special{pa 2374 2590}%
\special{pa 2554 2274}%
\special{pa 2736 2590}%
\special{pa 2374 2590}%
\special{pa 2554 2274}%
\special{fp}%
}}%
%
{\color[named]{Black}{%
\special{pn 4}%
\special{pa 2374 2590}%
\special{pa 2554 2274}%
\special{pa 2736 2590}%
\special{pa 2374 2590}%
\special{pa 2554 2274}%
\special{fp}%
}}%
%
{\color[named]{Black}{%
\special{pn 4}%
\special{pa 2374 2590}%
\special{pa 2554 2274}%
\special{pa 2736 2590}%
\special{pa 2374 2590}%
\special{pa 2554 2274}%
\special{fp}%
}}%
%
{\color[named]{Black}{%
\special{pn 4}%
\special{pa 2554 2274}%
\special{pa 2736 1962}%
\special{pa 2916 2274}%
\special{pa 2554 2274}%
\special{pa 2736 1962}%
\special{fp}%
}}%
%
{\color[named]{Black}{%
\special{pn 4}%
\special{pa 2262 2590}%
\special{pa 1538 2590}%
\special{pa 1900 1962}%
\special{pa 2262 2590}%
\special{pa 1538 2590}%
\special{fp}%
}}%
%
{\color[named]{Black}{%
\special{pn 4}%
\special{pa 1900 2590}%
\special{pa 2082 2274}%
\special{pa 2262 2590}%
\special{pa 1900 2590}%
\special{pa 2082 2274}%
\special{fp}%
}}%
%
{\color[named]{Black}{%
\special{pn 4}%
\special{pa 1900 2590}%
\special{pa 2082 2274}%
\special{pa 2262 2590}%
\special{pa 1900 2590}%
\special{pa 2082 2274}%
\special{fp}%
}}%
%
{\color[named]{Black}{%
\special{pn 4}%
\special{pa 1900 2590}%
\special{pa 2082 2274}%
\special{pa 2262 2590}%
\special{pa 1900 2590}%
\special{pa 2082 2274}%
\special{fp}%
}}%
%
{\color[named]{Black}{%
\special{pn 4}%
\special{pa 1900 2590}%
\special{pa 2082 2274}%
\special{pa 2262 2590}%
\special{pa 1900 2590}%
\special{pa 2082 2274}%
\special{fp}%
}}%
%
{\color[named]{Black}{%
\special{pn 4}%
\special{pa 1900 2590}%
\special{pa 2082 2274}%
\special{pa 2262 2590}%
\special{pa 1900 2590}%
\special{pa 2082 2274}%
\special{fp}%
}}%
%
{\color[named]{Black}{%
\special{pn 4}%
\special{pa 1538 2590}%
\special{pa 1718 2274}%
\special{pa 1900 2590}%
\special{pa 1538 2590}%
\special{pa 1718 2274}%
\special{fp}%
}}%
%
{\color[named]{Black}{%
\special{pn 4}%
\special{pa 1538 2590}%
\special{pa 1718 2274}%
\special{pa 1900 2590}%
\special{pa 1538 2590}%
\special{pa 1718 2274}%
\special{fp}%
}}%
%
{\color[named]{Black}{%
\special{pn 4}%
\special{pa 1538 2590}%
\special{pa 1718 2274}%
\special{pa 1900 2590}%
\special{pa 1538 2590}%
\special{pa 1718 2274}%
\special{fp}%
}}%
%
{\color[named]{Black}{%
\special{pn 4}%
\special{pa 1538 2590}%
\special{pa 1718 2274}%
\special{pa 1900 2590}%
\special{pa 1538 2590}%
\special{pa 1718 2274}%
\special{fp}%
}}%
%
{\color[named]{Black}{%
\special{pn 4}%
\special{pa 1718 2274}%
\special{pa 1900 1962}%
\special{pa 2082 2274}%
\special{pa 1718 2274}%
\special{pa 1900 1962}%
\special{fp}%
}}%
%
{\color[named]{Black}{%
\special{pn 4}%
\special{pa 1428 2590}%
\special{pa 702 2590}%
\special{pa 1064 1962}%
\special{pa 1428 2590}%
\special{pa 702 2590}%
\special{fp}%
}}%
%
{\color[named]{Black}{%
\special{pn 4}%
\special{pa 1064 2590}%
\special{pa 1246 2274}%
\special{pa 1428 2590}%
\special{pa 1064 2590}%
\special{pa 1246 2274}%
\special{fp}%
}}%
%
{\color[named]{Black}{%
\special{pn 4}%
\special{pa 1064 2590}%
\special{pa 1246 2274}%
\special{pa 1428 2590}%
\special{pa 1064 2590}%
\special{pa 1246 2274}%
\special{fp}%
}}%
%
{\color[named]{Black}{%
\special{pn 4}%
\special{pa 1064 2590}%
\special{pa 1246 2274}%
\special{pa 1428 2590}%
\special{pa 1064 2590}%
\special{pa 1246 2274}%
\special{fp}%
}}%
%
{\color[named]{Black}{%
\special{pn 4}%
\special{pa 1064 2590}%
\special{pa 1246 2274}%
\special{pa 1428 2590}%
\special{pa 1064 2590}%
\special{pa 1246 2274}%
\special{fp}%
}}%
%
{\color[named]{Black}{%
\special{pn 4}%
\special{pa 1064 2590}%
\special{pa 1246 2274}%
\special{pa 1428 2590}%
\special{pa 1064 2590}%
\special{pa 1246 2274}%
\special{fp}%
}}%
%
{\color[named]{Black}{%
\special{pn 4}%
\special{pa 702 2590}%
\special{pa 884 2274}%
\special{pa 1066 2590}%
\special{pa 702 2590}%
\special{pa 884 2274}%
\special{fp}%
}}%
%
{\color[named]{Black}{%
\special{pn 4}%
\special{pa 702 2590}%
\special{pa 884 2274}%
\special{pa 1066 2590}%
\special{pa 702 2590}%
\special{pa 884 2274}%
\special{fp}%
}}%
%
{\color[named]{Black}{%
\special{pn 4}%
\special{pa 702 2590}%
\special{pa 884 2274}%
\special{pa 1066 2590}%
\special{pa 702 2590}%
\special{pa 884 2274}%
\special{fp}%
}}%
%
{\color[named]{Black}{%
\special{pn 4}%
\special{pa 702 2590}%
\special{pa 884 2274}%
\special{pa 1066 2590}%
\special{pa 702 2590}%
\special{pa 884 2274}%
\special{fp}%
}}%
%
{\color[named]{Black}{%
\special{pn 4}%
\special{pa 884 2274}%
\special{pa 1064 1962}%
\special{pa 1246 2274}%
\special{pa 884 2274}%
\special{pa 1064 1962}%
\special{fp}%
}}%
%
{\color[named]{Black}{%
\special{pn 4}%
\special{pa 3966 1472}%
\special{pa 3242 1472}%
\special{pa 3604 844}%
\special{pa 3966 1472}%
\special{pa 3242 1472}%
\special{fp}%
}}%
%
{\color[named]{Black}{%
\special{pn 4}%
\special{pa 3604 1472}%
\special{pa 3786 1158}%
\special{pa 3966 1472}%
\special{pa 3604 1472}%
\special{pa 3786 1158}%
\special{fp}%
}}%
%
{\color[named]{Black}{%
\special{pn 4}%
\special{pa 3604 1472}%
\special{pa 3786 1158}%
\special{pa 3966 1472}%
\special{pa 3604 1472}%
\special{pa 3786 1158}%
\special{fp}%
}}%
%
{\color[named]{Black}{%
\special{pn 4}%
\special{pa 3604 1472}%
\special{pa 3786 1158}%
\special{pa 3966 1472}%
\special{pa 3604 1472}%
\special{pa 3786 1158}%
\special{fp}%
}}%
%
{\color[named]{Black}{%
\special{pn 4}%
\special{pa 3604 1472}%
\special{pa 3786 1158}%
\special{pa 3966 1472}%
\special{pa 3604 1472}%
\special{pa 3786 1158}%
\special{fp}%
}}%
%
{\color[named]{Black}{%
\special{pn 4}%
\special{pa 3604 1472}%
\special{pa 3786 1158}%
\special{pa 3966 1472}%
\special{pa 3604 1472}%
\special{pa 3786 1158}%
\special{fp}%
}}%
%
{\color[named]{Black}{%
\special{pn 4}%
\special{pa 3242 1472}%
\special{pa 3422 1158}%
\special{pa 3604 1472}%
\special{pa 3242 1472}%
\special{pa 3422 1158}%
\special{fp}%
}}%
%
{\color[named]{Black}{%
\special{pn 4}%
\special{pa 3242 1472}%
\special{pa 3422 1158}%
\special{pa 3604 1472}%
\special{pa 3242 1472}%
\special{pa 3422 1158}%
\special{fp}%
}}%
%
{\color[named]{Black}{%
\special{pn 4}%
\special{pa 3242 1472}%
\special{pa 3422 1158}%
\special{pa 3604 1472}%
\special{pa 3242 1472}%
\special{pa 3422 1158}%
\special{fp}%
}}%
%
{\color[named]{Black}{%
\special{pn 4}%
\special{pa 3242 1472}%
\special{pa 3422 1158}%
\special{pa 3604 1472}%
\special{pa 3242 1472}%
\special{pa 3422 1158}%
\special{fp}%
}}%
%
{\color[named]{Black}{%
\special{pn 4}%
\special{pa 3422 1158}%
\special{pa 3604 844}%
\special{pa 3786 1158}%
\special{pa 3422 1158}%
\special{pa 3604 844}%
\special{fp}%
}}%
%
{\color[named]{Black}{%
\special{pn 4}%
\special{pa 4802 1472}%
\special{pa 4076 1472}%
\special{pa 4438 844}%
\special{pa 4802 1472}%
\special{pa 4076 1472}%
\special{fp}%
}}%
%
{\color[named]{Black}{%
\special{pn 4}%
\special{pa 4438 1472}%
\special{pa 4620 1158}%
\special{pa 4802 1472}%
\special{pa 4438 1472}%
\special{pa 4620 1158}%
\special{fp}%
}}%
%
{\color[named]{Black}{%
\special{pn 4}%
\special{pa 4438 1472}%
\special{pa 4620 1158}%
\special{pa 4802 1472}%
\special{pa 4438 1472}%
\special{pa 4620 1158}%
\special{fp}%
}}%
%
{\color[named]{Black}{%
\special{pn 4}%
\special{pa 4438 1472}%
\special{pa 4620 1158}%
\special{pa 4802 1472}%
\special{pa 4438 1472}%
\special{pa 4620 1158}%
\special{fp}%
}}%
%
{\color[named]{Black}{%
\special{pn 4}%
\special{pa 4438 1472}%
\special{pa 4620 1158}%
\special{pa 4802 1472}%
\special{pa 4438 1472}%
\special{pa 4620 1158}%
\special{fp}%
}}%
%
{\color[named]{Black}{%
\special{pn 4}%
\special{pa 4438 1472}%
\special{pa 4620 1158}%
\special{pa 4802 1472}%
\special{pa 4438 1472}%
\special{pa 4620 1158}%
\special{fp}%
}}%
%
{\color[named]{Black}{%
\special{pn 4}%
\special{pa 4076 1472}%
\special{pa 4258 1158}%
\special{pa 4440 1472}%
\special{pa 4076 1472}%
\special{pa 4258 1158}%
\special{fp}%
}}%
%
{\color[named]{Black}{%
\special{pn 4}%
\special{pa 4076 1472}%
\special{pa 4258 1158}%
\special{pa 4440 1472}%
\special{pa 4076 1472}%
\special{pa 4258 1158}%
\special{fp}%
}}%
%
{\color[named]{Black}{%
\special{pn 4}%
\special{pa 4076 1472}%
\special{pa 4258 1158}%
\special{pa 4440 1472}%
\special{pa 4076 1472}%
\special{pa 4258 1158}%
\special{fp}%
}}%
%
{\color[named]{Black}{%
\special{pn 4}%
\special{pa 4076 1472}%
\special{pa 4258 1158}%
\special{pa 4440 1472}%
\special{pa 4076 1472}%
\special{pa 4258 1158}%
\special{fp}%
}}%
%
{\color[named]{Black}{%
\special{pn 4}%
\special{pa 4258 1158}%
\special{pa 4438 844}%
\special{pa 4620 1158}%
\special{pa 4258 1158}%
\special{pa 4438 844}%
\special{fp}%
}}%
%
{\color[named]{Black}{%
\special{pn 4}%
\special{pa 4768 2590}%
\special{pa 4042 2590}%
\special{pa 4406 1962}%
\special{pa 4768 2590}%
\special{pa 4042 2590}%
\special{fp}%
}}%
%
{\color[named]{Black}{%
\special{pn 4}%
\special{pa 4406 2590}%
\special{pa 4588 2274}%
\special{pa 4770 2590}%
\special{pa 4406 2590}%
\special{pa 4588 2274}%
\special{fp}%
}}%
%
{\color[named]{Black}{%
\special{pn 4}%
\special{pa 4406 2590}%
\special{pa 4588 2274}%
\special{pa 4770 2590}%
\special{pa 4406 2590}%
\special{pa 4588 2274}%
\special{fp}%
}}%
%
{\color[named]{Black}{%
\special{pn 4}%
\special{pa 4406 2590}%
\special{pa 4588 2274}%
\special{pa 4770 2590}%
\special{pa 4406 2590}%
\special{pa 4588 2274}%
\special{fp}%
}}%
%
{\color[named]{Black}{%
\special{pn 4}%
\special{pa 4406 2590}%
\special{pa 4588 2274}%
\special{pa 4770 2590}%
\special{pa 4406 2590}%
\special{pa 4588 2274}%
\special{fp}%
}}%
%
{\color[named]{Black}{%
\special{pn 4}%
\special{pa 4406 2590}%
\special{pa 4588 2274}%
\special{pa 4770 2590}%
\special{pa 4406 2590}%
\special{pa 4588 2274}%
\special{fp}%
}}%
%
{\color[named]{Black}{%
\special{pn 4}%
\special{pa 4042 2590}%
\special{pa 4224 2274}%
\special{pa 4406 2590}%
\special{pa 4042 2590}%
\special{pa 4224 2274}%
\special{fp}%
}}%
%
{\color[named]{Black}{%
\special{pn 4}%
\special{pa 4042 2590}%
\special{pa 4224 2274}%
\special{pa 4406 2590}%
\special{pa 4042 2590}%
\special{pa 4224 2274}%
\special{fp}%
}}%
%
{\color[named]{Black}{%
\special{pn 4}%
\special{pa 4042 2590}%
\special{pa 4224 2274}%
\special{pa 4406 2590}%
\special{pa 4042 2590}%
\special{pa 4224 2274}%
\special{fp}%
}}%
%
{\color[named]{Black}{%
\special{pn 4}%
\special{pa 4042 2590}%
\special{pa 4224 2274}%
\special{pa 4406 2590}%
\special{pa 4042 2590}%
\special{pa 4224 2274}%
\special{fp}%
}}%
%
{\color[named]{Black}{%
\special{pn 4}%
\special{pa 4224 2274}%
\special{pa 4406 1962}%
\special{pa 4588 2274}%
\special{pa 4224 2274}%
\special{pa 4406 1962}%
\special{fp}%
}}%
%
{\color[named]{Black}{%
\special{pn 4}%
\special{pa 3934 2590}%
\special{pa 3208 2590}%
\special{pa 3572 1962}%
\special{pa 3934 2590}%
\special{pa 3208 2590}%
\special{fp}%
}}%
%
{\color[named]{Black}{%
\special{pn 4}%
\special{pa 3572 2590}%
\special{pa 3752 2274}%
\special{pa 3934 2590}%
\special{pa 3572 2590}%
\special{pa 3752 2274}%
\special{fp}%
}}%
%
{\color[named]{Black}{%
\special{pn 4}%
\special{pa 3572 2590}%
\special{pa 3752 2274}%
\special{pa 3934 2590}%
\special{pa 3572 2590}%
\special{pa 3752 2274}%
\special{fp}%
}}%
%
{\color[named]{Black}{%
\special{pn 4}%
\special{pa 3572 2590}%
\special{pa 3752 2274}%
\special{pa 3934 2590}%
\special{pa 3572 2590}%
\special{pa 3752 2274}%
\special{fp}%
}}%
%
{\color[named]{Black}{%
\special{pn 4}%
\special{pa 3572 2590}%
\special{pa 3752 2274}%
\special{pa 3934 2590}%
\special{pa 3572 2590}%
\special{pa 3752 2274}%
\special{fp}%
}}%
%
{\color[named]{Black}{%
\special{pn 4}%
\special{pa 3572 2590}%
\special{pa 3752 2274}%
\special{pa 3934 2590}%
\special{pa 3572 2590}%
\special{pa 3752 2274}%
\special{fp}%
}}%
%
{\color[named]{Black}{%
\special{pn 4}%
\special{pa 3208 2590}%
\special{pa 3388 2274}%
\special{pa 3572 2590}%
\special{pa 3208 2590}%
\special{pa 3388 2274}%
\special{fp}%
}}%
%
{\color[named]{Black}{%
\special{pn 4}%
\special{pa 3208 2590}%
\special{pa 3388 2274}%
\special{pa 3572 2590}%
\special{pa 3208 2590}%
\special{pa 3388 2274}%
\special{fp}%
}}%
%
{\color[named]{Black}{%
\special{pn 4}%
\special{pa 3208 2590}%
\special{pa 3388 2274}%
\special{pa 3572 2590}%
\special{pa 3208 2590}%
\special{pa 3388 2274}%
\special{fp}%
}}%
%
{\color[named]{Black}{%
\special{pn 4}%
\special{pa 3208 2590}%
\special{pa 3388 2274}%
\special{pa 3572 2590}%
\special{pa 3208 2590}%
\special{pa 3388 2274}%
\special{fp}%
}}%
%
{\color[named]{Black}{%
\special{pn 4}%
\special{pa 3388 2274}%
\special{pa 3572 1962}%
\special{pa 3752 2274}%
\special{pa 3388 2274}%
\special{pa 3572 1962}%
\special{fp}%
}}%
%
{\color[named]{Black}{%
\special{pn 20}%
\special{pa 1570 1472}%
\special{pa 1934 1472}%
\special{fp}%
\special{pa 1934 1472}%
\special{pa 1750 1156}%
\special{fp}%
\special{pa 1750 1156}%
\special{pa 1932 840}%
\special{fp}%
}}%
%
{\color[named]{Black}{%
\special{pn 20}%
\special{pa 2406 1472}%
\special{pa 2584 1156}%
\special{fp}%
\special{pa 2584 1156}%
\special{pa 2948 1156}%
\special{fp}%
\special{pa 2948 1156}%
\special{pa 2766 842}%
\special{fp}%
}}%
%
{\color[named]{Black}{%
\special{pn 20}%
\special{pa 3238 1472}%
\special{pa 3602 1472}%
\special{fp}%
\special{pa 3602 1472}%
\special{pa 3420 1156}%
\special{fp}%
\special{pa 3420 1156}%
\special{pa 3784 1156}%
\special{fp}%
\special{pa 3784 1156}%
\special{pa 3602 842}%
\special{fp}%
}}%
%
{\color[named]{Black}{%
\special{pn 20}%
\special{pa 4076 1472}%
\special{pa 4438 1472}%
\special{fp}%
\special{pa 4438 1472}%
\special{pa 4618 1156}%
\special{fp}%
\special{pa 4618 1156}%
\special{pa 4256 1156}%
\special{fp}%
\special{pa 4256 1156}%
\special{pa 4438 842}%
\special{fp}%
}}%
%
{\color[named]{Black}{%
\special{pn 20}%
\special{pa 700 2590}%
\special{pa 882 2274}%
\special{fp}%
\special{pa 882 2274}%
\special{pa 1066 2590}%
\special{fp}%
\special{pa 1066 2590}%
\special{pa 1242 2274}%
\special{fp}%
\special{pa 1242 2274}%
\special{pa 1062 1960}%
\special{fp}%
}}%
%
{\color[named]{Black}{%
\special{pn 20}%
\special{pa 1534 2590}%
\special{pa 1898 2590}%
\special{fp}%
\special{pa 1898 2590}%
\special{pa 2080 2274}%
\special{fp}%
\special{pa 2080 2274}%
\special{pa 1898 1960}%
\special{fp}%
}}%
%
{\color[named]{Black}{%
\special{pn 20}%
\special{pa 2370 2590}%
\special{pa 3098 2590}%
\special{fp}%
\special{pa 3098 2590}%
\special{pa 2734 1960}%
\special{fp}%
}}%
%
{\color[named]{Black}{%
\special{pn 20}%
\special{pa 3206 2590}%
\special{pa 3932 2590}%
\special{fp}%
\special{pa 3932 2590}%
\special{pa 3752 2274}%
\special{fp}%
\special{pa 3752 2274}%
\special{pa 3388 2274}%
\special{fp}%
\special{pa 3388 2274}%
\special{pa 3570 1960}%
\special{fp}%
}}%
%
{\color[named]{Black}{%
\special{pn 20}%
\special{pa 4042 2590}%
\special{pa 4222 2274}%
\special{fp}%
\special{pa 4222 2274}%
\special{pa 4406 2590}%
\special{fp}%
\special{pa 4406 2590}%
\special{pa 4768 2590}%
\special{fp}%
\special{pa 4768 2590}%
\special{pa 4404 1960}%
\special{fp}%
}}%
\put(12.7100,-9.4700){\makebox(0,0)[lb]{$w^*_1$}}%
\put(21.0600,-9.4700){\makebox(0,0)[lb]{$w^*_2$}}%
\put(29.4200,-9.4700){\makebox(0,0)[lb]{$w^*_3$}}%
\put(37.7800,-9.4700){\makebox(0,0)[lb]{$w^*_4$}}%
\put(46.1200,-9.4700){\makebox(0,0)[lb]{$w^*_5$}}%
\put(12.3800,-21.0700){\makebox(0,0)[lb]{$w^*_6$}}%
\put(20.7300,-21.0700){\makebox(0,0)[lb]{$w^*_7$}}%
\put(29.0900,-21.0700){\makebox(0,0)[lb]{$w^*_8$}}%
\put(37.4500,-21.0700){\makebox(0,0)[lb]{$w^*_9$}}%
\put(45.7900,-21.0700){\makebox(0,0)[lb]{$w^*_{10}$}}%
\end{picture}%

%% file: Fig3.tex
\unitlength 0.1in
\begin{picture}( 66.1000, 36.3700)(  4.6000,-36.5700)
%
{\color[named]{Black}{%
\special{pn 8}%
\special{pa 3292 3578}%
\special{pa 3526 3174}%
\special{pa 3758 3578}%
\special{pa 3292 3578}%
\special{pa 3526 3174}%
\special{fp}%
}}%
%
{\color[named]{Black}{%
\special{pn 8}%
\special{pa 3292 3578}%
\special{pa 3526 3174}%
\special{pa 3758 3578}%
\special{pa 3292 3578}%
\special{pa 3526 3174}%
\special{fp}%
}}%
%
{\color[named]{Black}{%
\special{pn 8}%
\special{pa 3292 3578}%
\special{pa 3526 3174}%
\special{pa 3758 3578}%
\special{pa 3292 3578}%
\special{pa 3526 3174}%
\special{fp}%
}}%
%
{\color[named]{Black}{%
\special{pn 8}%
\special{pa 2826 3578}%
\special{pa 3058 3174}%
\special{pa 3294 3578}%
\special{pa 2826 3578}%
\special{pa 3058 3174}%
\special{fp}%
}}%
%
{\color[named]{Black}{%
\special{pn 8}%
\special{pa 2826 3578}%
\special{pa 3058 3174}%
\special{pa 3294 3578}%
\special{pa 2826 3578}%
\special{pa 3058 3174}%
\special{fp}%
}}%
%
{\color[named]{Black}{%
\special{pn 8}%
\special{pa 2826 3578}%
\special{pa 3058 3174}%
\special{pa 3294 3578}%
\special{pa 2826 3578}%
\special{pa 3058 3174}%
\special{fp}%
}}%
%
{\color[named]{Black}{%
\special{pn 8}%
\special{pa 2826 3578}%
\special{pa 3058 3174}%
\special{pa 3294 3578}%
\special{pa 2826 3578}%
\special{pa 3058 3174}%
\special{fp}%
}}%
%
{\color[named]{Black}{%
\special{pn 8}%
\special{pa 3058 3174}%
\special{pa 3292 2770}%
\special{pa 3526 3174}%
\special{pa 3058 3174}%
\special{pa 3292 2770}%
\special{fp}%
}}%
\put(37.4500,-36.2300){\makebox(0,0)[lt]{$b_1$}}%
\put(29.9000,-27.0000){\makebox(0,0)[lb]{$a_1$}}%
\put(27.8000,-36.2300){\makebox(0,0)[lt]{$O$}}%
%
{\color[named]{Black}{%
\special{pn 8}%
\special{pa 1698 2058}%
\special{pa 768 2058}%
\special{pa 1230 1252}%
\special{pa 1698 2058}%
\special{pa 768 2058}%
\special{fp}%
}}%
%
{\color[named]{Black}{%
\special{pn 8}%
\special{pa 1230 2058}%
\special{pa 1464 1654}%
\special{pa 1698 2058}%
\special{pa 1230 2058}%
\special{pa 1464 1654}%
\special{fp}%
}}%
%
{\color[named]{Black}{%
\special{pn 8}%
\special{pa 1230 2058}%
\special{pa 1464 1654}%
\special{pa 1698 2058}%
\special{pa 1230 2058}%
\special{pa 1464 1654}%
\special{fp}%
}}%
%
{\color[named]{Black}{%
\special{pn 8}%
\special{pa 1230 2058}%
\special{pa 1464 1654}%
\special{pa 1698 2058}%
\special{pa 1230 2058}%
\special{pa 1464 1654}%
\special{fp}%
}}%
%
{\color[named]{Black}{%
\special{pn 8}%
\special{pa 1230 2058}%
\special{pa 1464 1654}%
\special{pa 1698 2058}%
\special{pa 1230 2058}%
\special{pa 1464 1654}%
\special{fp}%
}}%
%
{\color[named]{Black}{%
\special{pn 8}%
\special{pa 768 2058}%
\special{pa 998 1654}%
\special{pa 1232 2058}%
\special{pa 768 2058}%
\special{pa 998 1654}%
\special{fp}%
}}%
%
{\color[named]{Black}{%
\special{pn 8}%
\special{pa 768 2058}%
\special{pa 998 1654}%
\special{pa 1232 2058}%
\special{pa 768 2058}%
\special{pa 998 1654}%
\special{fp}%
}}%
%
{\color[named]{Black}{%
\special{pn 8}%
\special{pa 768 2058}%
\special{pa 998 1654}%
\special{pa 1232 2058}%
\special{pa 768 2058}%
\special{pa 998 1654}%
\special{fp}%
}}%
%
{\color[named]{Black}{%
\special{pn 8}%
\special{pa 768 2058}%
\special{pa 998 1654}%
\special{pa 1232 2058}%
\special{pa 768 2058}%
\special{pa 998 1654}%
\special{fp}%
}}%
%
{\color[named]{Black}{%
\special{pn 8}%
\special{pa 998 1654}%
\special{pa 1230 1250}%
\special{pa 1464 1654}%
\special{pa 998 1654}%
\special{pa 1230 1250}%
\special{fp}%
}}%
\put(4.6000,-20.7000){\makebox(0,0)[lt]{$O$}}%
%
{\color[named]{Black}{%
\special{pn 8}%
\special{pa 2628 2056}%
\special{pa 1698 2056}%
\special{pa 2162 1250}%
\special{pa 2628 2056}%
\special{pa 1698 2056}%
\special{fp}%
}}%
%
{\color[named]{Black}{%
\special{pn 8}%
\special{pa 2162 2056}%
\special{pa 2396 1652}%
\special{pa 2628 2056}%
\special{pa 2162 2056}%
\special{pa 2396 1652}%
\special{fp}%
}}%
%
{\color[named]{Black}{%
\special{pn 8}%
\special{pa 2162 2056}%
\special{pa 2396 1652}%
\special{pa 2628 2056}%
\special{pa 2162 2056}%
\special{pa 2396 1652}%
\special{fp}%
}}%
%
{\color[named]{Black}{%
\special{pn 8}%
\special{pa 2162 2056}%
\special{pa 2396 1652}%
\special{pa 2628 2056}%
\special{pa 2162 2056}%
\special{pa 2396 1652}%
\special{fp}%
}}%
%
{\color[named]{Black}{%
\special{pn 8}%
\special{pa 2162 2056}%
\special{pa 2396 1652}%
\special{pa 2628 2056}%
\special{pa 2162 2056}%
\special{pa 2396 1652}%
\special{fp}%
}}%
%
{\color[named]{Black}{%
\special{pn 8}%
\special{pa 1698 2056}%
\special{pa 1930 1652}%
\special{pa 2162 2056}%
\special{pa 1698 2056}%
\special{pa 1930 1652}%
\special{fp}%
}}%
%
{\color[named]{Black}{%
\special{pn 8}%
\special{pa 1698 2056}%
\special{pa 1930 1652}%
\special{pa 2162 2056}%
\special{pa 1698 2056}%
\special{pa 1930 1652}%
\special{fp}%
}}%
%
{\color[named]{Black}{%
\special{pn 8}%
\special{pa 1698 2056}%
\special{pa 1930 1652}%
\special{pa 2162 2056}%
\special{pa 1698 2056}%
\special{pa 1930 1652}%
\special{fp}%
}}%
%
{\color[named]{Black}{%
\special{pn 8}%
\special{pa 1698 2056}%
\special{pa 1930 1652}%
\special{pa 2162 2056}%
\special{pa 1698 2056}%
\special{pa 1930 1652}%
\special{fp}%
}}%
%
{\color[named]{Black}{%
\special{pn 8}%
\special{pa 1930 1652}%
\special{pa 2162 1250}%
\special{pa 2396 1652}%
\special{pa 1930 1652}%
\special{pa 2162 1250}%
\special{fp}%
}}%
\put(26.9000,-21.3000){\makebox(0,0)[lt]{$b_2$}}%
%
{\color[named]{Black}{%
\special{pn 8}%
\special{pa 2162 1250}%
\special{pa 1230 1250}%
\special{pa 1698 442}%
\special{pa 2162 1250}%
\special{pa 1230 1250}%
\special{fp}%
}}%
%
{\color[named]{Black}{%
\special{pn 8}%
\special{pa 1230 1250}%
\special{pa 1464 844}%
\special{pa 1698 1250}%
\special{pa 1230 1250}%
\special{pa 1464 844}%
\special{fp}%
}}%
%
{\color[named]{Black}{%
\special{pn 8}%
\special{pa 1696 1250}%
\special{pa 1930 844}%
\special{pa 2162 1250}%
\special{pa 1696 1250}%
\special{pa 1930 844}%
\special{fp}%
}}%
%
{\color[named]{Black}{%
\special{pn 8}%
\special{pa 1464 844}%
\special{pa 1698 438}%
\special{pa 1930 844}%
\special{pa 1464 844}%
\special{pa 1698 438}%
\special{fp}%
}}%
\put(14.1000,-3.5000){\makebox(0,0)[lb]{$a_2$}}%
%
{\color[named]{Black}{%
\special{pn 20}%
\special{pa 3490 20}%
\special{pa 3490 20}%
\special{fp}%
}}%
%
{\color[named]{Black}{%
\special{pn 20}%
\special{pa 750 2054}%
\special{pa 1230 1250}%
\special{fp}%
}}%
%
{\color[named]{Black}{%
\special{pn 20}%
\special{pa 1230 1250}%
\special{pa 1692 2054}%
\special{fp}%
}}%
%
{\color[named]{Black}{%
\special{pn 20}%
\special{pa 1692 2054}%
\special{pa 1926 1652}%
\special{fp}%
}}%
%
{\color[named]{Black}{%
\special{pn 20}%
\special{pa 1926 1652}%
\special{pa 2378 1652}%
\special{fp}%
}}%
%
{\color[named]{Black}{%
\special{pn 20}%
\special{pa 2378 1652}%
\special{pa 2152 1240}%
\special{fp}%
}}%
%
{\color[named]{Black}{%
\special{pn 20}%
\special{pa 2152 1240}%
\special{pa 1702 446}%
\special{fp}%
}}%
%
{\color[named]{Black}{%
\special{pn 8}%
\special{pa 1162 3612}%
\special{pa 1394 3208}%
\special{pa 1628 3612}%
\special{pa 1162 3612}%
\special{pa 1394 3208}%
\special{fp}%
}}%
%
{\color[named]{Black}{%
\special{pn 8}%
\special{pa 1162 3612}%
\special{pa 1394 3208}%
\special{pa 1628 3612}%
\special{pa 1162 3612}%
\special{pa 1394 3208}%
\special{fp}%
}}%
%
{\color[named]{Black}{%
\special{pn 8}%
\special{pa 1162 3612}%
\special{pa 1394 3208}%
\special{pa 1628 3612}%
\special{pa 1162 3612}%
\special{pa 1394 3208}%
\special{fp}%
}}%
%
{\color[named]{Black}{%
\special{pn 8}%
\special{pa 696 3612}%
\special{pa 928 3208}%
\special{pa 1162 3612}%
\special{pa 696 3612}%
\special{pa 928 3208}%
\special{fp}%
}}%
%
{\color[named]{Black}{%
\special{pn 8}%
\special{pa 696 3612}%
\special{pa 928 3208}%
\special{pa 1162 3612}%
\special{pa 696 3612}%
\special{pa 928 3208}%
\special{fp}%
}}%
%
{\color[named]{Black}{%
\special{pn 8}%
\special{pa 696 3612}%
\special{pa 928 3208}%
\special{pa 1162 3612}%
\special{pa 696 3612}%
\special{pa 928 3208}%
\special{fp}%
}}%
%
{\color[named]{Black}{%
\special{pn 8}%
\special{pa 696 3612}%
\special{pa 928 3208}%
\special{pa 1162 3612}%
\special{pa 696 3612}%
\special{pa 928 3208}%
\special{fp}%
}}%
%
{\color[named]{Black}{%
\special{pn 8}%
\special{pa 928 3208}%
\special{pa 1162 2804}%
\special{pa 1394 3208}%
\special{pa 928 3208}%
\special{pa 1162 2804}%
\special{fp}%
}}%
\put(16.2000,-36.3000){\makebox(0,0)[lt]{$b_1$}}%
\put(8.5000,-27.2000){\makebox(0,0)[lb]{$a_1$}}%
\put(6.5000,-36.5700){\makebox(0,0)[lt]{$O$}}%
%
{\color[named]{Black}{%
\special{pn 20}%
\special{pa 680 3606}%
\special{pa 924 3204}%
\special{fp}%
}}%
%
{\color[named]{Black}{%
\special{pn 20}%
\special{pa 924 3204}%
\special{pa 1150 3606}%
\special{fp}%
}}%
%
{\color[named]{Black}{%
\special{pn 20}%
\special{pa 1150 3606}%
\special{pa 1386 3204}%
\special{fp}%
}}%
%
{\color[named]{Black}{%
\special{pn 20}%
\special{pa 1386 3204}%
\special{pa 1160 2792}%
\special{fp}%
}}%
%
{\color[named]{Black}{%
\special{pn 20}%
\special{pa 996 1652}%
\special{pa 1446 1652}%
\special{fp}%
}}%
%
{\color[named]{Black}{%
\special{pn 20}%
\special{pa 1466 838}%
\special{pa 1916 838}%
\special{fp}%
}}%
%
{\color[named]{Black}{%
\special{pn 8}%
\special{pa 4850 3558}%
\special{pa 5084 3154}%
\special{pa 5316 3558}%
\special{pa 4850 3558}%
\special{pa 5084 3154}%
\special{fp}%
}}%
%
{\color[named]{Black}{%
\special{pn 8}%
\special{pa 4850 3558}%
\special{pa 5084 3154}%
\special{pa 5316 3558}%
\special{pa 4850 3558}%
\special{pa 5084 3154}%
\special{fp}%
}}%
%
{\color[named]{Black}{%
\special{pn 8}%
\special{pa 4850 3558}%
\special{pa 5084 3154}%
\special{pa 5316 3558}%
\special{pa 4850 3558}%
\special{pa 5084 3154}%
\special{fp}%
}}%
%
{\color[named]{Black}{%
\special{pn 8}%
\special{pa 4386 3558}%
\special{pa 4616 3154}%
\special{pa 4852 3558}%
\special{pa 4386 3558}%
\special{pa 4616 3154}%
\special{fp}%
}}%
%
{\color[named]{Black}{%
\special{pn 8}%
\special{pa 4386 3558}%
\special{pa 4616 3154}%
\special{pa 4852 3558}%
\special{pa 4386 3558}%
\special{pa 4616 3154}%
\special{fp}%
}}%
%
{\color[named]{Black}{%
\special{pn 8}%
\special{pa 4386 3558}%
\special{pa 4616 3154}%
\special{pa 4852 3558}%
\special{pa 4386 3558}%
\special{pa 4616 3154}%
\special{fp}%
}}%
%
{\color[named]{Black}{%
\special{pn 8}%
\special{pa 4386 3558}%
\special{pa 4616 3154}%
\special{pa 4852 3558}%
\special{pa 4386 3558}%
\special{pa 4616 3154}%
\special{fp}%
}}%
%
{\color[named]{Black}{%
\special{pn 8}%
\special{pa 4616 3154}%
\special{pa 4850 2752}%
\special{pa 5084 3154}%
\special{pa 4616 3154}%
\special{pa 4850 2752}%
\special{fp}%
}}%
\put(53.0300,-36.0300){\makebox(0,0)[lt]{$b_1$}}%
\put(45.8000,-26.7000){\makebox(0,0)[lb]{$a_1$}}%
\put(43.3800,-36.0300){\makebox(0,0)[lt]{$O$}}%
%
{\color[named]{Black}{%
\special{pn 8}%
\special{pa 6418 3558}%
\special{pa 6652 3154}%
\special{pa 6884 3558}%
\special{pa 6418 3558}%
\special{pa 6652 3154}%
\special{fp}%
}}%
%
{\color[named]{Black}{%
\special{pn 8}%
\special{pa 6418 3558}%
\special{pa 6652 3154}%
\special{pa 6884 3558}%
\special{pa 6418 3558}%
\special{pa 6652 3154}%
\special{fp}%
}}%
%
{\color[named]{Black}{%
\special{pn 8}%
\special{pa 6418 3558}%
\special{pa 6652 3154}%
\special{pa 6884 3558}%
\special{pa 6418 3558}%
\special{pa 6652 3154}%
\special{fp}%
}}%
%
{\color[named]{Black}{%
\special{pn 8}%
\special{pa 5954 3558}%
\special{pa 6184 3154}%
\special{pa 6420 3558}%
\special{pa 5954 3558}%
\special{pa 6184 3154}%
\special{fp}%
}}%
%
{\color[named]{Black}{%
\special{pn 8}%
\special{pa 5954 3558}%
\special{pa 6184 3154}%
\special{pa 6420 3558}%
\special{pa 5954 3558}%
\special{pa 6184 3154}%
\special{fp}%
}}%
%
{\color[named]{Black}{%
\special{pn 8}%
\special{pa 5954 3558}%
\special{pa 6184 3154}%
\special{pa 6420 3558}%
\special{pa 5954 3558}%
\special{pa 6184 3154}%
\special{fp}%
}}%
%
{\color[named]{Black}{%
\special{pn 8}%
\special{pa 5954 3558}%
\special{pa 6184 3154}%
\special{pa 6420 3558}%
\special{pa 5954 3558}%
\special{pa 6184 3154}%
\special{fp}%
}}%
%
{\color[named]{Black}{%
\special{pn 8}%
\special{pa 6184 3154}%
\special{pa 6418 2752}%
\special{pa 6652 3154}%
\special{pa 6184 3154}%
\special{pa 6418 2752}%
\special{fp}%
}}%
\put(68.7100,-36.0300){\makebox(0,0)[lt]{$b_1$}}%
\put(61.5000,-26.5000){\makebox(0,0)[lb]{$a_1$}}%
\put(59.0600,-36.0300){\makebox(0,0)[lt]{$O$}}%
%
{\color[named]{Black}{%
\special{pn 20}%
\special{pa 2820 3572}%
\special{pa 3742 3572}%
\special{fp}%
}}%
%
{\color[named]{Black}{%
\special{pn 20}%
\special{pa 3742 3572}%
\special{pa 3300 2748}%
\special{fp}%
}}%
%
{\color[named]{Black}{%
\special{pn 20}%
\special{pa 3290 3562}%
\special{pa 3526 3160}%
\special{fp}%
}}%
%
{\color[named]{Black}{%
\special{pn 20}%
\special{pa 4388 3552}%
\special{pa 4614 3150}%
\special{fp}%
}}%
%
{\color[named]{Black}{%
\special{pn 20}%
\special{pa 4614 3150}%
\special{pa 5074 3150}%
\special{fp}%
}}%
%
{\color[named]{Black}{%
\special{pn 20}%
\special{pa 5074 3150}%
\special{pa 4838 2748}%
\special{fp}%
}}%
%
{\color[named]{Black}{%
\special{pn 20}%
\special{pa 5946 3552}%
\special{pa 6406 2748}%
\special{fp}%
}}%
%
{\color[named]{Black}{%
\special{pn 20}%
\special{pa 6182 3150}%
\special{pa 6632 3150}%
\special{fp}%
}}%
\put(21.8100,-7.8000){\makebox(0,0)[lb]{$w$}}%
\put(16.3000,-30.0700){\makebox(0,0)[lb]{$\tilde{w}$}}%
\put(36.4300,-29.3400){\makebox(0,0)[lb]{$w_1$}}%
\put(52.6900,-29.0500){\makebox(0,0)[lb]{$w_2$}}%
\put(68.1800,-29.0500){\makebox(0,0)[lb]{$w_3$}}%
%
{\color[named]{Black}{%
\special{pn 8}%
\special{pa 1700 450}%
\special{pa 1920 450}%
\special{fp}%
}}%
%
{\color[named]{Black}{%
\special{pn 8}%
\special{pa 2620 2050}%
\special{pa 2920 2050}%
\special{fp}%
}}%
%
{\color[named]{Black}{%
\special{pn 8}%
\special{pa 2630 2050}%
\special{pa 2760 1790}%
\special{fp}%
}}%
%
{\color[named]{Black}{%
\special{pn 20}%
\special{pa 1690 2030}%
\special{pa 2150 2030}%
\special{fp}%
\special{pa 2150 2030}%
\special{pa 2150 2060}%
\special{fp}%
}}%
%
{\color[named]{Black}{%
\special{pn 20}%
\special{pa 2150 2060}%
\special{pa 1710 2060}%
\special{fp}%
}}%
%
{\color[named]{Black}{%
\special{pn 8}%
\special{pa 1610 3610}%
\special{pa 1760 3610}%
\special{fp}%
}}%
%
{\color[named]{Black}{%
\special{pn 8}%
\special{pa 3300 2750}%
\special{pa 3470 2750}%
\special{fp}%
}}%
%
{\color[named]{Black}{%
\special{pn 8}%
\special{pa 3300 2760}%
\special{pa 3400 2580}%
\special{fp}%
}}%
%
{\color[named]{Black}{%
\special{pn 8}%
\special{pa 3730 3570}%
\special{pa 3900 3570}%
\special{fp}%
}}%
%
{\color[named]{Black}{%
\special{pn 8}%
\special{pa 5290 3550}%
\special{pa 5440 3550}%
\special{fp}%
}}%
%
{\color[named]{Black}{%
\special{pn 8}%
\special{pa 6890 3550}%
\special{pa 7070 3550}%
\special{fp}%
}}%
%
{\color[named]{Black}{%
\special{pn 8}%
\special{pa 6880 3550}%
\special{pa 7010 3320}%
\special{fp}%
}}%
%
{\color[named]{Black}{%
\special{pn 8}%
\special{pa 5320 3550}%
\special{pa 5440 3340}%
\special{fp}%
}}%
%
{\color[named]{Black}{%
\special{pn 8}%
\special{pa 3740 3580}%
\special{pa 3820 3440}%
\special{fp}%
}}%
%
{\color[named]{Black}{%
\special{pn 8}%
\special{pa 4830 2750}%
\special{pa 5030 2750}%
\special{fp}%
}}%
%
{\color[named]{Black}{%
\special{pn 8}%
\special{pa 6390 2740}%
\special{pa 6620 2740}%
\special{fp}%
}}%
%
{\color[named]{Black}{%
\special{pn 8}%
\special{pa 6400 2740}%
\special{pa 6520 2530}%
\special{fp}%
}}%
%
{\color[named]{Black}{%
\special{pn 8}%
\special{pa 4830 2750}%
\special{pa 4830 2750}%
\special{fp}%
\special{pa 4950 2540}%
\special{pa 4950 2540}%
\special{fp}%
}}%
%
{\color[named]{Black}{%
\special{pn 20}%
\special{pa 4380 3520}%
\special{pa 4820 3520}%
\special{fp}%
\special{pa 4840 3550}%
\special{pa 4840 3520}%
\special{fp}%
}}%
%
{\color[named]{Black}{%
\special{pn 20}%
\special{pa 4830 3550}%
\special{pa 4370 3550}%
\special{fp}%
}}%
%
{\color[named]{Black}{%
\special{pn 8}%
\special{pa 1680 450}%
\special{pa 1810 220}%
\special{fp}%
}}%
%
{\color[named]{Black}{%
\special{pn 8}%
\special{pa 1160 2810}%
\special{pa 1320 2810}%
\special{fp}%
}}%
%
{\color[named]{Black}{%
\special{pn 8}%
\special{pa 1160 2800}%
\special{pa 1280 2590}%
\special{fp}%
}}%
%
{\color[named]{Black}{%
\special{pn 8}%
\special{pa 4840 2750}%
\special{pa 4970 2520}%
\special{fp}%
}}%
%
{\color[named]{Black}{%
\special{pn 8}%
\special{pa 1630 3610}%
\special{pa 1740 3420}%
\special{fp}%
}}%
\end{picture}%

%% file: Fig5.tex
\unitlength 0.1in
\begin{picture}( 41.1000, 38.0700)(  7.4000,-42.0700)
%
{\color[named]{Black}{%
\special{pn 4}%
\special{pa 2562 1904}%
\special{pa 890 1904}%
\special{pa 1726 518}%
\special{pa 2562 1904}%
\special{pa 890 1904}%
\special{fp}%
}}%
%
{\color[named]{Black}{%
\special{pn 4}%
\special{pa 1726 1904}%
\special{pa 2144 1212}%
\special{pa 2562 1904}%
\special{pa 1726 1904}%
\special{pa 2144 1212}%
\special{fp}%
}}%
%
{\color[named]{Black}{%
\special{pn 4}%
\special{pa 1726 1904}%
\special{pa 2144 1212}%
\special{pa 2562 1904}%
\special{pa 1726 1904}%
\special{pa 2144 1212}%
\special{fp}%
}}%
%
{\color[named]{Black}{%
\special{pn 4}%
\special{pa 1726 1904}%
\special{pa 2144 1212}%
\special{pa 2562 1904}%
\special{pa 1726 1904}%
\special{pa 2144 1212}%
\special{fp}%
}}%
%
{\color[named]{Black}{%
\special{pn 4}%
\special{pa 1726 1904}%
\special{pa 2144 1212}%
\special{pa 2562 1904}%
\special{pa 1726 1904}%
\special{pa 2144 1212}%
\special{fp}%
}}%
%
{\color[named]{Black}{%
\special{pn 4}%
\special{pa 1726 1904}%
\special{pa 2144 1212}%
\special{pa 2562 1904}%
\special{pa 1726 1904}%
\special{pa 2144 1212}%
\special{fp}%
}}%
%
{\color[named]{Black}{%
\special{pn 4}%
\special{pa 890 1904}%
\special{pa 1308 1212}%
\special{pa 1726 1904}%
\special{pa 890 1904}%
\special{pa 1308 1212}%
\special{fp}%
}}%
%
{\color[named]{Black}{%
\special{pn 4}%
\special{pa 890 1904}%
\special{pa 1308 1212}%
\special{pa 1726 1904}%
\special{pa 890 1904}%
\special{pa 1308 1212}%
\special{fp}%
}}%
%
{\color[named]{Black}{%
\special{pn 4}%
\special{pa 890 1904}%
\special{pa 1308 1212}%
\special{pa 1726 1904}%
\special{pa 890 1904}%
\special{pa 1308 1212}%
\special{fp}%
}}%
%
{\color[named]{Black}{%
\special{pn 4}%
\special{pa 890 1904}%
\special{pa 1308 1212}%
\special{pa 1726 1904}%
\special{pa 890 1904}%
\special{pa 1308 1212}%
\special{fp}%
}}%
%
{\color[named]{Black}{%
\special{pn 4}%
\special{pa 1308 1212}%
\special{pa 1726 518}%
\special{pa 2144 1212}%
\special{pa 1308 1212}%
\special{pa 1726 518}%
\special{fp}%
}}%
\put(17.9000,-5.3000){\makebox(0,0)[lb]{$a_N$}}%
%
{\color[named]{Black}{%
\special{pn 4}%
\special{pa 4702 1902}%
\special{pa 3030 1902}%
\special{pa 3866 516}%
\special{pa 4702 1902}%
\special{pa 3030 1902}%
\special{fp}%
}}%
%
{\color[named]{Black}{%
\special{pn 4}%
\special{pa 3866 1902}%
\special{pa 4284 1208}%
\special{pa 4702 1902}%
\special{pa 3866 1902}%
\special{pa 4284 1208}%
\special{fp}%
}}%
%
{\color[named]{Black}{%
\special{pn 4}%
\special{pa 3866 1902}%
\special{pa 4284 1208}%
\special{pa 4702 1902}%
\special{pa 3866 1902}%
\special{pa 4284 1208}%
\special{fp}%
}}%
%
{\color[named]{Black}{%
\special{pn 4}%
\special{pa 3866 1902}%
\special{pa 4284 1208}%
\special{pa 4702 1902}%
\special{pa 3866 1902}%
\special{pa 4284 1208}%
\special{fp}%
}}%
%
{\color[named]{Black}{%
\special{pn 4}%
\special{pa 3866 1902}%
\special{pa 4284 1208}%
\special{pa 4702 1902}%
\special{pa 3866 1902}%
\special{pa 4284 1208}%
\special{fp}%
}}%
%
{\color[named]{Black}{%
\special{pn 4}%
\special{pa 3866 1902}%
\special{pa 4284 1208}%
\special{pa 4702 1902}%
\special{pa 3866 1902}%
\special{pa 4284 1208}%
\special{fp}%
}}%
%
{\color[named]{Black}{%
\special{pn 4}%
\special{pa 3030 1902}%
\special{pa 3448 1208}%
\special{pa 3868 1902}%
\special{pa 3030 1902}%
\special{pa 3448 1208}%
\special{fp}%
}}%
%
{\color[named]{Black}{%
\special{pn 4}%
\special{pa 3030 1902}%
\special{pa 3448 1208}%
\special{pa 3868 1902}%
\special{pa 3030 1902}%
\special{pa 3448 1208}%
\special{fp}%
}}%
%
{\color[named]{Black}{%
\special{pn 4}%
\special{pa 3030 1902}%
\special{pa 3448 1208}%
\special{pa 3868 1902}%
\special{pa 3030 1902}%
\special{pa 3448 1208}%
\special{fp}%
}}%
%
{\color[named]{Black}{%
\special{pn 4}%
\special{pa 3030 1902}%
\special{pa 3448 1208}%
\special{pa 3868 1902}%
\special{pa 3030 1902}%
\special{pa 3448 1208}%
\special{fp}%
}}%
%
{\color[named]{Black}{%
\special{pn 4}%
\special{pa 3448 1208}%
\special{pa 3866 514}%
\special{pa 4284 1208}%
\special{pa 3448 1208}%
\special{pa 3866 514}%
\special{fp}%
}}%
%
{\color[named]{Black}{%
\special{pn 8}%
\special{pa 2558 1902}%
\special{pa 2536 1880}%
\special{pa 2512 1858}%
\special{pa 2486 1838}%
\special{pa 2458 1826}%
\special{pa 2434 1848}%
\special{pa 2414 1870}%
\special{pa 2392 1856}%
\special{pa 2364 1826}%
\special{pa 2332 1802}%
\special{pa 2308 1808}%
\special{pa 2296 1838}%
\special{pa 2280 1868}%
\special{pa 2258 1866}%
\special{pa 2236 1838}%
\special{pa 2218 1800}%
\special{pa 2202 1772}%
\special{pa 2168 1768}%
\special{pa 2146 1792}%
\special{pa 2128 1792}%
\special{pa 2108 1756}%
\special{pa 2074 1740}%
\special{pa 2044 1752}%
\special{pa 2034 1782}%
\special{pa 2042 1814}%
\special{pa 2020 1826}%
\special{pa 1978 1818}%
\special{pa 1948 1796}%
\special{pa 1920 1806}%
\special{pa 1908 1838}%
\special{pa 1906 1856}%
\special{pa 1870 1838}%
\special{pa 1846 1818}%
\special{pa 1812 1834}%
\special{pa 1806 1866}%
\special{pa 1778 1872}%
\special{pa 1750 1886}%
\special{pa 1724 1904}%
\special{pa 1720 1908}%
\special{fp}%
}}%
\put(7.9200,-19.4900){\makebox(0,0)[lt]{$O$}}%
%
{\color[named]{Black}{%
\special{pn 4}%
\special{pa 754 1902}%
\special{pa 910 1902}%
\special{fp}%
}}%
%
{\color[named]{Black}{%
\special{pn 4}%
\special{pa 2906 1902}%
\special{pa 3062 1902}%
\special{fp}%
}}%
%
{\color[named]{Black}{%
\special{pn 4}%
\special{pa 886 1902}%
\special{pa 772 1714}%
\special{fp}%
\special{pa 788 1738}%
\special{pa 798 1738}%
\special{fp}%
}}%
%
{\color[named]{Black}{%
\special{pn 4}%
\special{pa 3028 1902}%
\special{pa 2916 1714}%
\special{fp}%
\special{pa 2932 1738}%
\special{pa 2940 1738}%
\special{fp}%
}}%
%
{\color[named]{Black}{%
\special{pn 8}%
\special{pa 2560 1902}%
\special{pa 2542 1874}%
\special{pa 2526 1848}%
\special{pa 2508 1822}%
\special{pa 2486 1798}%
\special{pa 2462 1776}%
\special{pa 2436 1758}%
\special{pa 2404 1742}%
\special{pa 2368 1730}%
\special{pa 2332 1722}%
\special{pa 2304 1710}%
\special{pa 2298 1688}%
\special{pa 2310 1656}%
\special{pa 2324 1622}%
\special{pa 2328 1588}%
\special{pa 2322 1554}%
\special{pa 2310 1524}%
\special{pa 2288 1496}%
\special{pa 2262 1474}%
\special{pa 2230 1454}%
\special{pa 2200 1438}%
\special{pa 2180 1416}%
\special{pa 2182 1386}%
\special{pa 2186 1352}%
\special{pa 2178 1322}%
\special{pa 2160 1294}%
\special{pa 2148 1264}%
\special{pa 2144 1232}%
\special{pa 2144 1208}%
\special{fp}%
}}%
%
{\color[named]{Black}{%
\special{pn 8}%
\special{pa 2138 1208}%
\special{pa 2118 1184}%
\special{pa 2098 1162}%
\special{pa 2072 1140}%
\special{pa 2042 1120}%
\special{pa 2010 1100}%
\special{pa 1990 1080}%
\special{pa 1998 1058}%
\special{pa 2026 1036}%
\special{pa 2032 1032}%
\special{fp}%
}}%
%
{\color[named]{Black}{%
\special{pn 8}%
\special{pa 2138 1208}%
\special{pa 2126 1238}%
\special{pa 2116 1268}%
\special{pa 2110 1298}%
\special{pa 2104 1330}%
\special{pa 2104 1364}%
\special{pa 2110 1398}%
\special{pa 2104 1428}%
\special{pa 2074 1442}%
\special{pa 2046 1458}%
\special{pa 2042 1490}%
\special{pa 2048 1524}%
\special{pa 2040 1550}%
\special{pa 2008 1562}%
\special{pa 1978 1576}%
\special{pa 1970 1602}%
\special{pa 1976 1638}%
\special{pa 1972 1668}%
\special{pa 1950 1686}%
\special{pa 1916 1696}%
\special{pa 1884 1710}%
\special{pa 1860 1734}%
\special{pa 1844 1760}%
\special{pa 1828 1788}%
\special{pa 1808 1814}%
\special{pa 1784 1836}%
\special{pa 1762 1860}%
\special{pa 1742 1884}%
\special{pa 1724 1902}%
\special{fp}%
}}%
%
{\color[named]{Black}{%
\special{pn 4}%
\special{pa 740 4026}%
\special{pa 896 4026}%
\special{fp}%
}}%
%
{\color[named]{Black}{%
\special{pn 4}%
\special{pa 884 4024}%
\special{pa 770 3834}%
\special{fp}%
\special{pa 788 3862}%
\special{pa 796 3862}%
\special{fp}%
}}%
\put(16.3200,-21.1900){\makebox(0,0)[lt]{(a)}}%
%
{\color[named]{Black}{%
\special{pn 20}%
\special{pa 3028 1902}%
\special{pa 4702 1902}%
\special{fp}%
\special{pa 4702 1902}%
\special{pa 4280 1208}%
\special{fp}%
\special{pa 4280 1208}%
\special{pa 3864 1902}%
\special{fp}%
\special{pa 3864 1902}%
\special{pa 3446 1208}%
\special{fp}%
\special{pa 3446 1208}%
\special{pa 3864 518}%
\special{fp}%
}}%
\put(25.9000,-19.4000){\makebox(0,0)[lt]{$b_N$}}%
%
{\color[named]{Black}{%
\special{pn 8}%
\special{pa 880 1908}%
\special{pa 898 1884}%
\special{pa 922 1862}%
\special{pa 950 1846}%
\special{pa 982 1842}%
\special{pa 1012 1844}%
\special{pa 1046 1842}%
\special{pa 1068 1858}%
\special{pa 1098 1862}%
\special{pa 1136 1872}%
\special{pa 1166 1870}%
\special{pa 1168 1838}%
\special{pa 1182 1810}%
\special{pa 1206 1786}%
\special{pa 1232 1764}%
\special{pa 1258 1746}%
\special{pa 1290 1734}%
\special{pa 1310 1752}%
\special{pa 1310 1784}%
\special{pa 1326 1810}%
\special{pa 1358 1828}%
\special{pa 1392 1834}%
\special{pa 1424 1834}%
\special{pa 1456 1842}%
\special{pa 1488 1840}%
\special{pa 1518 1832}%
\special{pa 1550 1826}%
\special{pa 1576 1844}%
\special{pa 1594 1880}%
\special{pa 1608 1906}%
\special{pa 1628 1896}%
\special{pa 1652 1868}%
\special{pa 1684 1868}%
\special{pa 1706 1888}%
\special{pa 1706 1892}%
\special{fp}%
}}%
%
{\color[named]{Black}{%
\special{pn 8}%
\special{pa 1310 1228}%
\special{pa 1310 1258}%
\special{pa 1308 1288}%
\special{pa 1298 1328}%
\special{pa 1292 1364}%
\special{pa 1314 1376}%
\special{pa 1340 1396}%
\special{pa 1352 1426}%
\special{pa 1352 1490}%
\special{pa 1346 1524}%
\special{pa 1352 1554}%
\special{pa 1370 1580}%
\special{pa 1394 1592}%
\special{pa 1422 1604}%
\special{pa 1438 1628}%
\special{pa 1448 1664}%
\special{pa 1470 1686}%
\special{pa 1498 1706}%
\special{pa 1528 1718}%
\special{pa 1558 1722}%
\special{pa 1586 1742}%
\special{pa 1602 1774}%
\special{pa 1606 1810}%
\special{pa 1638 1808}%
\special{pa 1654 1832}%
\special{pa 1684 1842}%
\special{pa 1680 1862}%
\special{pa 1656 1858}%
\special{fp}%
}}%
%
{\color[named]{Black}{%
\special{pn 8}%
\special{pa 1302 1194}%
\special{pa 1334 1186}%
\special{pa 1362 1172}%
\special{pa 1390 1156}%
\special{pa 1418 1142}%
\special{pa 1448 1130}%
\special{pa 1480 1122}%
\special{pa 1512 1112}%
\special{pa 1538 1096}%
\special{pa 1562 1074}%
\special{pa 1586 1050}%
\special{pa 1606 1024}%
\special{pa 1614 996}%
\special{pa 1608 966}%
\special{pa 1600 934}%
\special{pa 1598 902}%
\special{pa 1598 870}%
\special{pa 1606 838}%
\special{pa 1624 810}%
\special{pa 1654 796}%
\special{pa 1690 794}%
\special{pa 1710 776}%
\special{pa 1724 744}%
\special{pa 1728 716}%
\special{pa 1718 686}%
\special{pa 1710 654}%
\special{pa 1712 626}%
\special{pa 1736 596}%
\special{pa 1730 574}%
\special{pa 1724 546}%
\special{pa 1724 530}%
\special{fp}%
}}%
%
{\color[named]{Black}{%
\special{pn 8}%
\special{pa 1706 784}%
\special{pa 1706 816}%
\special{pa 1708 846}%
\special{pa 1708 910}%
\special{pa 1740 920}%
\special{pa 1764 942}%
\special{pa 1800 946}%
\special{pa 1826 930}%
\special{pa 1824 896}%
\special{pa 1834 864}%
\special{pa 1826 838}%
\special{pa 1798 816}%
\special{pa 1776 794}%
\special{pa 1748 784}%
\special{pa 1718 776}%
\special{pa 1698 784}%
\special{fp}%
}}%
%
{\color[named]{Black}{%
\special{pn 4}%
\special{pa 2558 4020}%
\special{pa 896 4020}%
\special{pa 1726 2634}%
\special{pa 2558 4020}%
\special{pa 896 4020}%
\special{fp}%
}}%
%
{\color[named]{Black}{%
\special{pn 4}%
\special{pa 1726 4020}%
\special{pa 2144 3326}%
\special{pa 2560 4020}%
\special{pa 1726 4020}%
\special{pa 2144 3326}%
\special{fp}%
}}%
%
{\color[named]{Black}{%
\special{pn 4}%
\special{pa 1726 4020}%
\special{pa 2144 3326}%
\special{pa 2560 4020}%
\special{pa 1726 4020}%
\special{pa 2144 3326}%
\special{fp}%
}}%
%
{\color[named]{Black}{%
\special{pn 4}%
\special{pa 1726 4020}%
\special{pa 2144 3326}%
\special{pa 2560 4020}%
\special{pa 1726 4020}%
\special{pa 2144 3326}%
\special{fp}%
}}%
%
{\color[named]{Black}{%
\special{pn 4}%
\special{pa 1726 4020}%
\special{pa 2144 3326}%
\special{pa 2560 4020}%
\special{pa 1726 4020}%
\special{pa 2144 3326}%
\special{fp}%
}}%
%
{\color[named]{Black}{%
\special{pn 4}%
\special{pa 1726 4020}%
\special{pa 2144 3326}%
\special{pa 2560 4020}%
\special{pa 1726 4020}%
\special{pa 2144 3326}%
\special{fp}%
}}%
%
{\color[named]{Black}{%
\special{pn 4}%
\special{pa 896 4020}%
\special{pa 1312 3326}%
\special{pa 1728 4020}%
\special{pa 896 4020}%
\special{pa 1312 3326}%
\special{fp}%
}}%
%
{\color[named]{Black}{%
\special{pn 4}%
\special{pa 896 4020}%
\special{pa 1312 3326}%
\special{pa 1728 4020}%
\special{pa 896 4020}%
\special{pa 1312 3326}%
\special{fp}%
}}%
%
{\color[named]{Black}{%
\special{pn 4}%
\special{pa 896 4020}%
\special{pa 1312 3326}%
\special{pa 1728 4020}%
\special{pa 896 4020}%
\special{pa 1312 3326}%
\special{fp}%
}}%
%
{\color[named]{Black}{%
\special{pn 4}%
\special{pa 896 4020}%
\special{pa 1312 3326}%
\special{pa 1728 4020}%
\special{pa 896 4020}%
\special{pa 1312 3326}%
\special{fp}%
}}%
%
{\color[named]{Black}{%
\special{pn 4}%
\special{pa 1312 3326}%
\special{pa 1726 2634}%
\special{pa 2144 3326}%
\special{pa 1312 3326}%
\special{pa 1726 2634}%
\special{fp}%
}}%
%
{\color[named]{Black}{%
\special{pn 20}%
\special{pa 892 4018}%
\special{pa 1726 4018}%
\special{fp}%
\special{pa 1726 4018}%
\special{pa 1312 3326}%
\special{fp}%
\special{pa 1312 3326}%
\special{pa 1726 2634}%
\special{fp}%
}}%
\put(37.8800,-22.1600){\makebox(0,0)[lb]{(b)}}%
\put(31.2200,-20.8500){\makebox(0,0)[lb]{$O$}}%
\put(39.5000,-5.3000){\makebox(0,0)[lb]{$a_N$}}%
\put(47.4000,-20.6000){\makebox(0,0)[lb]{$b_N$}}%
\put(16.6000,-43.3700){\makebox(0,0)[lb]{(c)}}%
\put(8.1000,-42.3700){\makebox(0,0)[lb]{$O$}}%
\put(18.0000,-26.3000){\makebox(0,0)[lb]{$a_N$}}%
\put(26.1000,-41.5000){\makebox(0,0)[lb]{$b_N$}}%
%
{\color[named]{Black}{%
\special{pn 8}%
\special{pa 1350 1560}%
\special{pa 1324 1540}%
\special{pa 1294 1530}%
\special{pa 1264 1532}%
\special{pa 1236 1536}%
\special{pa 1230 1570}%
\special{pa 1220 1600}%
\special{pa 1222 1630}%
\special{pa 1256 1630}%
\special{pa 1276 1642}%
\special{pa 1308 1640}%
\special{pa 1340 1636}%
\special{pa 1358 1608}%
\special{pa 1360 1580}%
\special{fp}%
}}%
%
{\color[named]{Black}{%
\special{pn 4}%
\special{pa 4798 4048}%
\special{pa 3126 4048}%
\special{pa 3962 2662}%
\special{pa 4798 4048}%
\special{pa 3126 4048}%
\special{fp}%
}}%
%
{\color[named]{Black}{%
\special{pn 4}%
\special{pa 3962 4048}%
\special{pa 4380 3354}%
\special{pa 4798 4048}%
\special{pa 3962 4048}%
\special{pa 4380 3354}%
\special{fp}%
}}%
%
{\color[named]{Black}{%
\special{pn 4}%
\special{pa 3962 4048}%
\special{pa 4380 3354}%
\special{pa 4798 4048}%
\special{pa 3962 4048}%
\special{pa 4380 3354}%
\special{fp}%
}}%
%
{\color[named]{Black}{%
\special{pn 4}%
\special{pa 3962 4048}%
\special{pa 4380 3354}%
\special{pa 4798 4048}%
\special{pa 3962 4048}%
\special{pa 4380 3354}%
\special{fp}%
}}%
%
{\color[named]{Black}{%
\special{pn 4}%
\special{pa 3962 4048}%
\special{pa 4380 3354}%
\special{pa 4798 4048}%
\special{pa 3962 4048}%
\special{pa 4380 3354}%
\special{fp}%
}}%
%
{\color[named]{Black}{%
\special{pn 4}%
\special{pa 3962 4048}%
\special{pa 4380 3354}%
\special{pa 4798 4048}%
\special{pa 3962 4048}%
\special{pa 4380 3354}%
\special{fp}%
}}%
%
{\color[named]{Black}{%
\special{pn 4}%
\special{pa 3126 4048}%
\special{pa 3544 3354}%
\special{pa 3962 4048}%
\special{pa 3126 4048}%
\special{pa 3544 3354}%
\special{fp}%
}}%
%
{\color[named]{Black}{%
\special{pn 4}%
\special{pa 3126 4048}%
\special{pa 3544 3354}%
\special{pa 3962 4048}%
\special{pa 3126 4048}%
\special{pa 3544 3354}%
\special{fp}%
}}%
%
{\color[named]{Black}{%
\special{pn 4}%
\special{pa 3126 4048}%
\special{pa 3544 3354}%
\special{pa 3962 4048}%
\special{pa 3126 4048}%
\special{pa 3544 3354}%
\special{fp}%
}}%
%
{\color[named]{Black}{%
\special{pn 4}%
\special{pa 3126 4048}%
\special{pa 3544 3354}%
\special{pa 3962 4048}%
\special{pa 3126 4048}%
\special{pa 3544 3354}%
\special{fp}%
}}%
%
{\color[named]{Black}{%
\special{pn 4}%
\special{pa 3544 3354}%
\special{pa 3962 2660}%
\special{pa 4380 3354}%
\special{pa 3544 3354}%
\special{pa 3962 2660}%
\special{fp}%
}}%
\put(40.2000,-26.6000){\makebox(0,0)[lb]{$a_N$}}%
\put(30.2800,-40.9200){\makebox(0,0)[lt]{$O$}}%
%
{\color[named]{Black}{%
\special{pn 4}%
\special{pa 2990 4046}%
\special{pa 3146 4046}%
\special{fp}%
}}%
%
{\color[named]{Black}{%
\special{pn 4}%
\special{pa 3122 4046}%
\special{pa 3008 3856}%
\special{fp}%
\special{pa 3024 3882}%
\special{pa 3034 3882}%
\special{fp}%
}}%
\put(39.2000,-41.9700){\makebox(0,0)[lt]{(d)}}%
\put(48.5000,-41.2000){\makebox(0,0)[lt]{$b_N$}}%
%
{\color[named]{Black}{%
\special{pn 8}%
\special{pa 3116 4050}%
\special{pa 3134 4026}%
\special{pa 3158 4004}%
\special{pa 3186 3990}%
\special{pa 3218 3986}%
\special{pa 3248 3988}%
\special{pa 3282 3984}%
\special{pa 3304 4002}%
\special{pa 3334 4004}%
\special{pa 3372 4014}%
\special{pa 3402 4012}%
\special{pa 3404 3980}%
\special{pa 3418 3952}%
\special{pa 3442 3928}%
\special{pa 3468 3908}%
\special{pa 3494 3890}%
\special{pa 3526 3878}%
\special{pa 3546 3894}%
\special{pa 3546 3928}%
\special{pa 3562 3952}%
\special{pa 3594 3970}%
\special{pa 3628 3978}%
\special{pa 3660 3978}%
\special{pa 3692 3984}%
\special{pa 3724 3984}%
\special{pa 3754 3974}%
\special{pa 3786 3968}%
\special{pa 3812 3986}%
\special{pa 3830 4024}%
\special{pa 3844 4050}%
\special{pa 3864 4038}%
\special{pa 3888 4012}%
\special{pa 3920 4010}%
\special{pa 3942 4030}%
\special{pa 3942 4034}%
\special{fp}%
}}%
%
{\color[named]{Black}{%
\special{pn 8}%
\special{pa 3546 3370}%
\special{pa 3546 3402}%
\special{pa 3544 3432}%
\special{pa 3534 3470}%
\special{pa 3528 3508}%
\special{pa 3550 3518}%
\special{pa 3576 3540}%
\special{pa 3588 3570}%
\special{pa 3588 3634}%
\special{pa 3582 3668}%
\special{pa 3588 3698}%
\special{pa 3606 3724}%
\special{pa 3630 3736}%
\special{pa 3658 3748}%
\special{pa 3674 3772}%
\special{pa 3684 3808}%
\special{pa 3706 3830}%
\special{pa 3734 3848}%
\special{pa 3764 3862}%
\special{pa 3794 3864}%
\special{pa 3822 3884}%
\special{pa 3838 3916}%
\special{pa 3842 3952}%
\special{pa 3874 3950}%
\special{pa 3890 3974}%
\special{pa 3920 3984}%
\special{pa 3916 4006}%
\special{pa 3892 4002}%
\special{fp}%
}}%
%
{\color[named]{Black}{%
\special{pn 8}%
\special{pa 3538 3338}%
\special{pa 3570 3330}%
\special{pa 3598 3314}%
\special{pa 3626 3300}%
\special{pa 3654 3284}%
\special{pa 3684 3272}%
\special{pa 3748 3256}%
\special{pa 3774 3240}%
\special{pa 3822 3192}%
\special{pa 3842 3168}%
\special{pa 3850 3140}%
\special{pa 3844 3108}%
\special{pa 3836 3076}%
\special{pa 3834 3044}%
\special{pa 3834 3012}%
\special{pa 3842 2980}%
\special{pa 3860 2954}%
\special{pa 3890 2940}%
\special{pa 3926 2936}%
\special{pa 3946 2918}%
\special{pa 3960 2888}%
\special{pa 3964 2860}%
\special{pa 3954 2830}%
\special{pa 3946 2798}%
\special{pa 3948 2770}%
\special{pa 3972 2738}%
\special{pa 3966 2718}%
\special{pa 3960 2690}%
\special{pa 3960 2672}%
\special{fp}%
}}%
%
{\color[named]{Black}{%
\special{pn 8}%
\special{pa 3942 2926}%
\special{pa 3942 2958}%
\special{pa 3944 2990}%
\special{pa 3944 3052}%
\special{pa 3976 3064}%
\special{pa 4000 3086}%
\special{pa 4036 3090}%
\special{pa 4062 3072}%
\special{pa 4060 3040}%
\special{pa 4070 3006}%
\special{pa 4062 2980}%
\special{pa 4034 2960}%
\special{pa 4012 2938}%
\special{pa 3984 2928}%
\special{pa 3954 2918}%
\special{pa 3934 2926}%
\special{fp}%
}}%
%
{\color[named]{Black}{%
\special{pn 8}%
\special{pa 3586 3704}%
\special{pa 3560 3682}%
\special{pa 3530 3674}%
\special{pa 3500 3674}%
\special{pa 3472 3678}%
\special{pa 3466 3714}%
\special{pa 3456 3744}%
\special{pa 3458 3774}%
\special{pa 3492 3772}%
\special{pa 3512 3784}%
\special{pa 3544 3784}%
\special{pa 3576 3780}%
\special{pa 3594 3752}%
\special{pa 3596 3724}%
\special{fp}%
}}%
%
{\color[named]{Black}{%
\special{pn 20}%
\special{pa 3110 4048}%
\special{pa 3950 4048}%
\special{fp}%
}}%
%
{\color[named]{Black}{%
\special{pn 8}%
\special{pa 3940 4038}%
\special{pa 3950 4028}%
\special{fp}%
}}%
%
{\color[named]{Black}{%
\special{pn 20}%
\special{pa 3940 4048}%
\special{pa 3930 4048}%
\special{fp}%
\special{pa 3550 3328}%
\special{pa 3550 3328}%
\special{fp}%
}}%
%
{\color[named]{Black}{%
\special{pn 20}%
\special{pa 3950 4058}%
\special{pa 3950 4058}%
\special{fp}%
\special{pa 3940 4048}%
\special{pa 3530 3348}%
\special{fp}%
}}%
%
{\color[named]{Black}{%
\special{pn 20}%
\special{pa 3530 3348}%
\special{pa 3530 3348}%
\special{fp}%
\special{pa 3540 3348}%
\special{pa 3960 2658}%
\special{fp}%
}}%
%
\put(47.4000,-20.6000){\makebox(0,0)[lb]{}}%
\end{picture}%